\documentclass[12pt,a4paper]{amsart}
\usepackage[cp1251]{inputenc}
\usepackage[T2A]{fontenc}
\usepackage{amsfonts}
\usepackage{amsthm}
\usepackage{cite}
\usepackage[russian]{babel}
\usepackage{amsmath,amsthm}
\usepackage{amssymb}
\usepackage{mathrsfs,amscd}

\renewcommand{\le}{\leqslant}
\renewcommand{\ge}{\geqslant}

\usepackage{geometry} 
\newtheoremstyle{citedth}%
  {5pt}
  {5pt}
  {\itshape}
  {}
  {\bfseries}
  {.}
  {.5em}
  {\thmname{#1} \thmnumber{#2} \thmnote{\normalfont#3}}

  \newtheoremstyle{citedexample}%
  {5pt}
  {5pt}
  {}
  {}
  {\bfseries}
  {.}
  {.5em}
  {\thmname{#1} \thmnumber{#2} \thmnote{\normalfont#3}}

\theoremstyle{citedth}

\theoremstyle{citedexample}

\theoremstyle{theorem}
\newtheorem{theorem-}{Теорема}[section]
\newtheorem{lemma-}{Лемма}[section]
\theoremstyle{definition}
\newtheorem{example-}{Пример}[section]

\begin{document}

\bigskip

\bigskip

\bigskip

\begin{center}
{\bf {EXTREMAL AND APPROXIMATIVE PROPERTIES\\ OF SIMPLE PARTIAL
FRACTIONS} }
\bigskip

{\small \sc V.\,I. Danchenko, M.\,A. Komarov and P.\,V. Chunaev}

\end{center}

\bigskip

\bigskip

{\small

\textbf{Abstract.} In approximation theory, logarithmic
derivatives of complex polynomials are called \textit{simple
partial fractions} (\textit{SPF}) as suggested by
Eu.~P.~Dolzhenko. Many solved and unsolved extremal problems
related to SPF are traced back to works of G.~Boole, A.\,J.
Macintyre, W.\,H.\,J. Fuchs, J.\,M.~Marstrand, E.~A.~Gorin,
A.~A.~Gonchar, Eu.~P.~Dol\-zhenko. At present, many authors
systematically develop methods for approximation and interpolation
by SPF and several their modifications. Simultaneously, related
problems, being of independent interest, arise for SPF: inequalities of different metrics, estimation of derivatives, separation of singularities, etc.

We systematize some of these problems which are known to us in
Introduction of this survey. In the main part, we formulate
principal results and  outline methods to prove them if possible.
}

\bigskip

\bigskip

\bigskip

\bigskip

\bigskip

\begin{center}
{\bf ЭКСТРЕМАЛЬНЫЕ И АППРОКСИМАТИВНЫЕ СВОЙСТВА НАИПРОСТЕЙШИХ ДРОБЕЙ\footnote{Работа В.\,И. Данченко выполнена при финансовой поддержке Минобрнауки России (задание №1.574.2016/1.4). Работа П.\,В. Чунаева выполнена при финансовой поддержке РФФИ в рамках научного проекта №16-31-00252 мол\_а.}}

\bigskip

\bigskip

{\small \sc В.~И. Данченко\footnote{Кафедра функционального анализа и его приложений, Владимирский государственный университет, ул. Горького, д. 87, г. Владимир, 600000, Россия; e-mail:  \textsf{vdanch2012@yandex.ru}}, М.~А. Комаров\footnote{Кафедра функционального анализа и его приложений, Владимирский государственный университет, ул. Горького, д. 87, г. Владимир, 600000, Россия; e-mail:  \textsf{kami9@yandex.ru}}, П.~В. Чунаев\footnote{Кафедра функционального анализа и его приложений, Владимирский государственный университет, ул. Горького, д. 87, г. Владимир, 600000, Россия и Departament de Matem\`{a}tiques, Universitat Aut\`{o}noma de Barcelona, 08193 Bellaterra (Barcelona), Spain; e-mail:  \textsf{chunayev@mail.ru} }}

\end{center}

\bigskip

{\small

\textbf{Аннотация.} Наипростейшими дробями (НД) по предложению Е.\,П. Долженко в теории
аппроксимаций называют логарифмические производные алгебраических
многочленов. С НД связано много решенных и нерешенных задач
экстремального характера, восходящих к работам Дж.~Буля,
А.\,Дж.~Ма\-кин\-тай\-ра, У.\,Х.\,Дж.~Фукса, Дж.\,M.~Марстранда,
Е.\,А.~Горина, А.\,А.~Гончара, Е.\,П.~Дол\-жен\-ко. В настоящее время
многими авторами систематически развиваются методы аппроксимации и
интерполяции посредством НД и некоторых их модификаций и
обобщений. Параллельно для НД возникают и смежные задачи,
представляющие самостоятельный интерес: неравенства разных метрик,
оценки производных, разделение особенностей и др.

Во вводной части обзора в какой-то мере систематизированы
известные авторам задачи такого рода, а в основной части
сформулированы основные результаты и по возможности намечены
подходы к их доказательствам.
}

\newpage

\tableofcontents

\section{ Введение. История вопросов}
\label{paragraph1}
 Наипростейшими дробями (НД) порядка не выше $n=0,1,\ldots$
называются рациональные дроби вида
\begin{equation}
\rho_0(z)\equiv 0,\qquad \rho_n(z)=
\sum_{k=1}^{n}\;\frac{1}{z-z_k},\quad n\in {\mathbb N},\quad
z,\,z_k\in \overline{\mathbb C},
   \label{(*1)}
\end{equation}
где при $z_k=\infty$ считаем $\frac{1}{z-z_k}\equiv 0$. Очевидно,
НД $\rho_n(z)$ является логарифмической производной многочлена
$Q_n(z)=\prod_{k=1}^{n}\;(z-z_k)$, где при $z_k=\infty$ считаем
$z-z_k\equiv 1$. Приведем краткую информацию об основных старых и
новых задачах, связанных с НД, и современном состоянии теории НД.

\subsection{Оценки картановского типа}
\label{paragraph1.1} Метрические свойста НД начали систематически
изучаться в работах Дж.~Буля~\cite{[BOOL]}, А.~Дж.~Макин\-тай\-ра
и У.~Х.~Дж.\,Фук\-са~\cite{[M-F]}, Дж.~M.~Марстранда~\cite{[MAR]},
А.~А.~Гончара~\cite{[GON1],[GON2],[GON3]},
Е.~П.~Долженко~\cite{[DOL1],[DOL2],[DOL3],[DOL4]}. Так или иначе,
большинство исследований было связано с методом покрытий
А.~Картана \cite{[cartan]}. Первые оценки картановского типа для
НД были получены в \cite{[M-F]}. При $\delta>0$ и $n\in {\mathbb
N}$ положим
$$
E(\rho_n;\delta)=\{z:\; |\rho_n(z)|>\delta\};\qquad
L(\rho_n;\delta)=\inf\sum {\rm diam}\,(d_j),
$$
где нижняя грань берется по всем конечным покрытиям
$E(\rho_n;\delta)\subset \cup d_j$ наборами кругов $d_j$. В
\cite{[M-F]} показано, что, независимо от расположения полюсов НД,
\begin{equation}
L(\rho_n;\delta)\le {\rm const}\cdot {A(n)}\,{\delta^{-1}},\qquad
{\rm const}>0,\quad n\in {\mathbb N},
  \label{(001)}
\end{equation}
c $A(n)=n(1+\ln n)$.

Если же рассматривать покрытия множества
$E_0(\rho_n;\delta)=E(\rho_n;\delta)\cap \mathbb R$, то оценка
(\ref{(001)}) верна с $A(n)=n$ и она точна по порядку величины $n$
(см. \cite{[BOOL],[M-F]}). Если к тому же все полюсы $z_k$ конечны
и лежат на $\mathbb R$, то $E_0(\rho_n;\delta)$ состоит из $n$
интервалов и ${\rm mes}_1 E_0(\rho_n;\delta)=2n\delta^{-1}$.
Поэтому возник вопрос о точности величины $A(n)$ в общем случае
 (\ref{(001)}): нужен ли логарифмический множитель \cite{[M-F]}? Частично на этот
вопрос дан ответ в \cite{[MAR]}, где приведен пример, когда
справедлива оценка, противоположная (\ref{(001)}) с $A(n)= n
\sqrt{\ln n/\ln\ln n}$, ${n\ge n_0}$. В~2005~г. Дж.~М.~Андерсоном
и В.~Я.~Эйдерманом \cite{[AND-EID1],[AND-EID2]} была получена
окончательная точная по порядку оценка (\ref{(001)}) с $A(n)= n
\sqrt{\ln n}$.

Разные модификации покрытий и оценок картановского типа
применялись Е.\,П. Дол\-жен\-ко \cite{[DOL1]}, Н.~В. Говоровым и
Ю.~П.~Лапенко \cite{[GOVOR-LAP]}; они сыграли важную роль в
обратных теоремах теории рациональных аппроксимаций. Наиболее
общие результаты по этой тематике недавно получены в работе
В.~Я.~Эйдермана \cite{[EID1]}, где вместо НД оценивались
потенциалы Коши $\int (z-\zeta )^{-1}d\mu (\zeta)$.

\subsection{Задача Горина}
\label{paragraph1.2}

Другой круг задач для НД связан с проблемой
Е.~А.~Горина \cite{[GOR]} об оценке величин
\begin{equation}
d_n(\mathbb R,p)=\inf\left\{Y(\rho_n):\;\|\rho_n\|_{L_{p}({\mathbb
R})}=1 \right\},\qquad 1<p\le\infty,
  \label{(*2)}
\end{equation}
где $Y(\rho_n)=\min_{k=\overline{1,n}} |{\rm Im}\,z_k|$. При
замене $\tilde\rho_n(z)=c\rho_n(c\,z)$,
$c=\|\rho_n\|^{-q}_{L_{p}({\mathbb R})}$, сохраняющей вид НД,
имеем $\|\tilde\rho_n\|_{L_{p}({\mathbb R})}=1$ и
$Y(\tilde\rho_n)=c^{-1}Y(\rho_n)$. Поэтому определение
(\ref{(*2)}) можно переписать в виде
\begin{equation}
d_n(\mathbb R,p)=\inf_{\rho_n}\left\{Y(\rho_n)
\|\rho_n\|^{q}_{L_{p}({\mathbb R})} \right\},\qquad
p^{-1}+q^{-1}=1,
  \label{(*2bis)}
\end{equation}
где инфимум берется по всем НД вида (\ref{(*1)}), не имеющим
полюсов на ${\mathbb R}$. Следовательно, задачу Горина можно
интерпретировать как задачу о величине наилучшего приближения нуля
(наименьшего уклонения от нуля) в метрике $L_{p}({\mathbb R})$ в
классе НД $\rho_n$ порядка $\le n$, при условии $Y(\rho_n)=1$,
или, что то же самое, при условии, что все они имеют общий
фиксированный полюс, например, $z_1=i$. В этом смысле задача
Горина является аналогом классической экстремальной задачи
Чебышева о наименьшем уклонении от нуля унитарного многочлена
фиксированной степени. Вследствие этого задачу об оценке
$d_n(\mathbb R,p)$ будем иногда называть задачей Горина-Чебышева.

Вопрос о принципиальной возможности оценки снизу для $d_n(\mathbb
R,\infty)$ величиной\footnote { Далее через
$A(\cdot,\cdot,\ldots)$ (с индексами и без них) обозначаются
конечные положительные величины, зависящие лишь от указанных
аргументов и индексов и, вообще говоря, различные в различных
формулах.} $A(n)>0$ был поставлен и положительно решен в 1962 г.
Е.~А.~Гориным \cite{[GOR]}. В~1965~г. Е.~Г.~Николаев \cite{[NIK]}
получил оценку
$$
d_n(\mathbb R,\infty)\ge 2\,(\sqrt{2}-1)^{n-1},\qquad n\in {\mathbb
N}.
$$
Он также поставил вопрос о точности этой оценки и следующую
проблему: \textit{вообще, верно ли, что $d_n(\mathbb R,\infty)\to 0$ при
$n\to\infty$?}

В связи с этим в \cite{[NIK]} приведен принадлежащий
А.\,Н.~Колмогорову пример целой функции с {\it бесконечным} числом
нулей $\pm\gamma+i\pi k/(2a)$, $k\in \mathbb Z$:
$$
f(z)=f(a,\gamma;z)=\left( e^{ia\,( z+i\gamma) }+e^{-ia\,
(z+i\gamma)} \right) \left( {e^{ia\,( z-i \gamma) }}+e^{-ia(
z-i\gamma)} \right),\qquad a>0,\; \gamma>0.
$$
Несложно проверить, что $\|f'/f\|_{L_{\infty}({\mathbb R})}<1$ при
$a=\gamma^{-2}$. Таким образом, с сохранением нормировки
логарифмической производной $f'/f$ расстояние $\gamma$ от ее
полюсов до ${\mathbb R}$ можно сделать сколь угодно малым.

В 1966 г. существенное уточнение оценки Е.~Г.~Николаева получено
А.~О.~Гельфондом \cite{[GEL]}:
$$
d_n(\mathbb R,\infty)>(17\ln n)^{-1},\qquad n\ge n_{0}.
$$
В.~Э.~Кацнельсон \cite{[KAC]} получил некоторое уточнение этой
оценки, но при том же логарифмическом порядке убывания миноранты.
Вопросы, сформулированные Е.~Г. Николаевым, оставались открытыми. В
1994 г. окончательный результат был установлен в \cite{[DAN1]}:
$$
d_n(\mathbb R,\infty)\asymp \frac{\ln\ln n}{\ln n}.
$$
Кроме того, оказалось \cite{[DAN1]}, что при конечных $p$ величины
$d_n(\mathbb R,p)$ к нулю не убывают и
\begin{equation}
\inf_n d_n(\mathbb R,p)\ge A(p),\qquad A(p):={2^{q/p}}{p^{-1}}
\sin^q({\pi}{p^{-1}}).
  \label{(*3)}
\end{equation}
Эта оценка была несколько уточнена в \cite{[BOR2007]} (см.
(\ref{(*9)})); вопрос о точноcти величины $A(p)$ (хотя бы по
порядку величины $p$) остается открытым.

Рассматривалась аналогичная задача для потенциала Коши $\int
(z-\zeta )^{-1}d\mu (\zeta)$ неотрицательной борелевской меры
$\mu$, ${\rm supp}\, \mu \subset \mathbb{C}^+$, и получена оценка
снизу усредненного расстояния от ${\rm supp}\, \mu$ до
действительной оси \cite{[DAN-DISS]}.

\subsection{Модификации задачи Горина-Чебышева}
\label{paragraph1.3} Представляют интерес модификации задачи об
оценке величины  (\ref{(*2)}) с заменой $\mathbb R$ на другие
множества (окружности, ляпуновские кривые, отрезки, лучи и др.). В
случае конечных $p$ такие модификации мало изучены, но имеется ряд
окончательных результатов в случае $p=\infty$. Например, в случае
единичной окружности $\gamma_1$ и отрезка $[-1,1]$ получены слабые
эквивалентности \cite{[DAN2006],{DDJ}}:
$d_n(\gamma_1,\infty)\asymp n^{-1}\ln n$,
$d_n([-1,1],\infty)\asymp n^{-2}\ln^{2}n$ ($n\ge n_0$). В
\cite{[KOS3]} получена оценка {\it снизу} такого же порядка, как
для отрезка, для так называемых спрямляемых
компактов\footnote{Компакт $K$ называется спрямляемым, если он не
разбивает плоскость и существует такая положительная величина
$a(K)<\infty$, что любые две его точки можно соединить кривой
$L\subset K$ длины $\le a(K)$.}. Во всех этих случаях, как и в
случае прямой, проблему можно трактовать как задачу
Горина-Чебышева о НД, наименее уклоняющейся от нуля в равномерной
метрике при условии закрепленности одного полюса. Такие
экстремальные и близкие к экстремальным НД в случае некоторых
весовых пространств на отрезке $[-1,1]$ найдены в
\cite{[Chu2014]}.

Представляет интерес и другая модификация
--- задача Горина-Чебышева о НД со {\it свободными полюсами},
наименее уклоняющейся от ненулевой константы (см. \S
\ref{paragraph8.3}). Вместо констант в этой задаче рассматривались
и дробно-линейные функции~\cite{[KOND-diss]}.

\subsection{Задача Гельфонда}
\label{paragraph1.4} А.~О.~Гельфонд \cite{[GEL]} рассматривал
задачу об оценке расстояний $d' _n(\mathbb R,\infty)$ типа
(\ref{(*2)}), но с нормировкой производной:
$\|\rho'_n\|_{L_{\infty}({\mathbb R})}\le 1$. В \cite{[GEL]} была
получена оценка снизу для $d'_n(\mathbb R,\infty)$, но довольно
неточная, носящая скорее качественный характер (см.
\S\ref{paragraph2.3}). В \cite{[DAN1]} доказана слабая
эквивалентность указанных расстояний с величиной $n^{-1/2}\ln n$,
но при дополнительном условии, что все полюсы НД лежат в верхней
полуплоскости. В общем же случае известно \cite{[DANChu2013+]},
что расстояния $d' _n(\mathbb R,\infty)$ ограничены снизу
величиной порядка $\sqrt{n^{-1}\ln n}$. Вопрос о точности этой
оценки по порядку остается открытым. Кроме того, не изучен вопрос
об оценках расстояний $d' _n(\mathbb R,p)$ типа (\ref{(*2)}) с
нормировкой производной: $\|\rho'_n\|_{L_{p}({\mathbb R})}\le 1$
при конечных $p>1/2$.

\subsection{$L_p$-оценки НД и их производных}
\label{paragraph1.8} Оценки $L_p$-норм НД, зависящие явно от
полюсов $z_k$, были найдены В.~Ю.~Прота\-со\-вым,
И.~Р.~Каю\-мо\-вым, А.~В.~Каюмовой и др. \cite{[PRO],[DAN2010],
[KAY2011],[KAY2012], [KAY2012+],[KAYva]}. Особый интерес
представляют неравенства разных метрик (термин
С.~М.~Николь\-ско\-го), т.е. оценки $L_p$-норм НД сверху через их
$L_r$-нормы на различных множествах при различных $p>1$ и $r>1$.
Для алгебраических и тригонометрических многочленов такого рода
оценки на ограниченных промежутках $K\subset {\mathbb R}$ при
$p>r$ (только этот случай и представляет интерес, т.к. при $p<r$
оценки получаются из неравенства Гельдера) хорошо известны
благодаря классическим работам Д.~Джек\-со\-на,
С.~М.~Николь\-ско\-го, Н.~К.~Бари, С.~Б.~Стечкина, В.~В.~Арестова,
Л.~В.~Тайкова, Г.~Сеге, А.~Зигмунда, П.~Л.~Ульянова,
М.~К.~По\-та\-по\-ва, П.~Борвейна, А.~Ф.~Ти\-ма\-на,
И.~И.~Иб\-ра\-ги\-мова и многих других авторов (см., например,
\cite{Nikolskiy,Bari,Timan}). Напомним, что первым результатом
такого рода является следующее $(p,r)$-неравенство
Джек\-со\-на-Никольского для алгебраических многочленов $P_n$
степени $\le n$:
\begin{equation}
\|P_n\|_{L_p(K)}\le A(K,p,r)
n^{2\left(\frac{1}{r}-\frac{1}{p}\right)} \|P_n\|_{L_r(K)},\qquad
p>r\ge 1.
   \label{(*4)}
\end{equation}
Сходные неравенства верны и для тригонометрических многочленов,
причем в случае вещественных многочленов и $K=[-\pi,\pi]$
множитель 2 в показателе степени можно заменить единицей и
положить $A(K,p,r)=2$ (см. \cite{Bari,Nikolskiy}).

Первое неравенство разных метрик для НД было получено на $\mathbb
R$ \cite{[DAN1]}:
\begin{equation}
 \| \rho_n \|_{L_{\infty} ({\mathbb{R}})} \le A(r)
\cdot \|\rho_n\|_{L_r ({\mathbb{R}})}^{s},\quad r>1,\quad
r^{-1}+s^{-1}=1,
  \label{(*5)}
\end{equation}
где $A(r)< 2r\,\sin^{-s}({\pi}{r^{-1}})$. Эта оценка улучшалась и
обобщалось на $(p,r)$-не\-ра\-вен\-ства на ограниченных и
неограниченных вещественных промежутках в работах
\cite{[Dodonov],[DanChu2016+],DAN-DOD}. В основе оценок лежит
метод насечки для построения квадратурных формул с переменными
узлами, который впервые был применен в работе \cite{[DAN2010]}, а
затем развивался в работах \cite{[DanChu2016+],[DAN-SEMIN]} (см.
\S\ref{paragraph4.1}).

Отметим, что в отличие от случая многочленов неравенства для НД
нелинейны относительно сравниваемых норм, содержательны и на
бесконечных промежутках и при произвольном соотношении между $p$ и
$r$ (при $p<r$ они также не являются тривиальными и не следуют из
неравенства Гельдера даже на конечных промежутках, поскольку не
завист от их длин, см., например, (\ref{(*446)})) \cite{DAN-DOD}.
Отметим еще, что оценки на ${\mathbb R}$ при $p>r$ не зависят от
порядка НД, однако такая зависимость появляется в случае $p<r$, а
также на ограниченных промежутках $K\subset {\mathbb R}$.
Неравенства разных метрик для НД получены также на окружностях
\cite{[D-D3]}.

С неравенствами разных метрик тесно связаны оценки типа
Маркова-Берн\-штей\-на для производных НД
\cite{[DAN1],[DAN2010],[DAN2006],DAN-DOD}. Например, в
\cite{DAN-DOD} показано, что если на $[-1,1]$ НД $\rho_n$
вещественнозначна и ее модуль ограничен единицей, то
\begin{equation}
 \sqrt{1-x^2}|\rho_n'(x)|\le
n(1+\varepsilon_n),\qquad x\in [-1,1],
 \label{(*6)}
\end{equation}
для некоторых положительных $\varepsilon_n$, где $\varepsilon_n\to
0$ при $n\to\infty$, причем оценка является точной по порядку
величины $n$.

\subsection{Интерполяция аналитических функций}
\label{paragraph-1-interpol}

Установлено \cite{[D-D1],[D-D2],[KOS2],DanChu2011}, что {\it НД
Паде} (НД $n$-кратной интерполяции с одним узлом) всегда
существует и единственна. Предложены различные конструкции НД Паде
и оценки погрешности. Например, в \cite{DanChu2011} применялось
представление НД Паде в виде интеграла Эрмита, благодаря чему
найдена явная формула для остаточного члена и получены его оценки.

В задаче интерполяции с {\it различными} узлами вопросы о
разрешимости и единственности значительно усложняются и, вообще
говоря, не имеют однозначного ответа. Связь особенностей этой
задачи с алгебраической структурой интерполяционных таблиц
исследовалась в \cite{[KOND],[DAN-KOND1],[DAN-KOND2]} и других
работах. Введено понятие {\it обобщенной} интерполяции таблиц,
охватывающее и обычную интерполяцию \cite{[DAN-KOND1],[KOM1]}.
Достаточно общие результаты о существовании и единственности
решения задачи интерполяции получены методом редукции к
полиномиальной интерполяции \cite{[KOM1]}.

\subsection{Аппроксимация посредством НД на компактах}
\label{paragraph1.5} Методы, разработанные для исследования
проблемы Горина (см. \S\ref{paragraph1.2}), позже использовались в
более общей задаче аппроксимации непрерывных функций на различных
множествах $K\subset \overline{\mathbb C}$, см.
\cite{[D-D1],[D-D2],[KOS1]}. Другими методами задача
приближения посредством НД решалась в работах Дж.~Кореваара,
Ч.~Чуи и К.~Шена \cite{[KOR64],[KOR],[Chui],[ChuiShen]} (например,
в интегральных пространствах Берса-Бергмана аналитических функций
на ограниченных жордановых областях $K$). В них была предложена
конструкция аппроксимирующих НД с полюсами, которые подбирались на
фиксированных множествах с определенными свойствами по отношению к
$K$. Аналогичный подход применялся в работах
\cite{[BOR2012],[BOR2016]}.

Одна из мотивировок аппроксимации посредством НД заключена в их
простом и важном физическом смысле: они задают (с точностью до
постоянных множителей и операции комплексного сопряжения) плоские
поля различной природы, создаваемые равновеликими источниками,
расположенными в точках $z_k$. В этом смысле задачу аппроксимации
посредством НД можно интерпретировать как задачу о размещении
источников $z_k$, создающих заранее заданное поле.

Систематическое изучение НД со свободными полюсами как аппарата
приближения началось в 1999 г., когда в \cite{[D-D1],[D-D2]} для них был доказан следующий аналог полиномиальной теоремы
С.~Н.~Мергеляна (ср. с \cite[Д.1]{WOLSH}):

\smallskip

{\it Для любого компакта $K$ со связным дополнением любую
непрерывную на $K$ функцию, аналитическую во внутренних точках на
$K$, можно сколь угодно точно приблизить посредством НД в
равномерной метрике}.\smallskip

При этом, как показал О.~Н.~Ко\-су\-хи\-н \cite{[KOS1]}, для
широкого класса компактов и функций скорости равномерного
приближения посредством НД и комплексных многочленов имеют
одинаковый порядок. Это позволило О.~Н.~Ко\-су\-хи\-ну получить
для НД ряд аналогов классических теорем Д.~Джексона,
С.~Н.~Бернштейна, А.~Зигмунда, В.~К.~Дзядыка, Дж.~Л.~Уолша. Из
дальнейших исследований, однако, стало ясно, что имеются и
значительные различия между аппроксимативными свойствами НД и
полиномов (см. \S\ref{paragraph8.2}). Собственно, именно эти
различия и вызывают дополнительный интерес к изучению НД.

Отметим, что хотя скорости приближения посредством НД и
многочленов на достаточно широком классе функций мало отличаются,
НД часто имеют преимущество с точки зрения их вычислительной
надежности. Это связано с тем, что в типичных случаях с ростом $n$
относительная погрешность при вычислениях НД $\rho_n$ растет
незначительно (фактически только из-за многократных сложений), в
то же время она быстро растет при вычислении многочленов $Q_n$
(пропорционально $n$ из-за многократных умножений)
\cite{[DEMIDOVICH]}. Это обстоятельство использовалось при
вычислениях рациональных функций общего вида в \cite{DanChu2011}
путем их аппроксимации посредством определенных модификаций НД.

\subsection{Аппроксимация другими аппаратами, основанными на НД}
\label{paragraph1.5++} Значительный интерес представляют разности
НД $\rho_{n_1}-\rho_{n_2}$, т.е. логарифмические производные
рациональных функций. Аппроксимативные свойства разностей НД
исследованы еще довольно мало, а результаты носят скорее
качественный характер. Примеры показывают, что в ряде случаев
скорость аппроксимации разностями НД гораздо более высокая, чем у
НД. Например, известен такой результат \cite{[D-D2]} об
аппроксимации разностями НД рациональных функций, имеющих
однозначный интеграл.

В \cite{[BOR2016]} получены результаты  о плотности в $AC(K)$
разностей НД с полюсами, расположенными на определенных
предписанных множествах (без анализа скорости аппроксимации). В
работе \cite{[Kom2017-IzVuz]} построены некоторые оценки
наилучшего весового (с весом $(x+c)^2$) приближения на
$\mathbb{R}^+$ разностями НД функций, убывающих на $+\infty$ со
скоростью $O(x^{-2})$.

Исследовались аппроксимативные свойства некоторых других
конструкций. Например, в \cite{[DAN2008],[DAN2006]} введены дроби
вида
$$
\rho_n (p,z)=\sum\nolimits_{k=1}^{n}
\frac{z_k^{-p}}{z-z_{k}},\qquad
\Theta(z)=\frac{\rho_{n_1}(z)-\rho_{n_2}(z)}
{\rho_{n_3}(z)-\rho_{n_4}(z)},
$$
где $p$ --- целое число, $z_k\ne 0,\infty$, так что
$\rho_n(z)=\rho_n(0,z)$, $\rho_{n_k}$ --- НД. Показано, что при
$p\ge 0$ для дробей $\rho_n(p,z)$ также справедлив аналог
полиномиальной теоремы Мергеляна. При $p<0$ это уже, вообще
говоря, неверно (это неверно, например, если $K$ содержит
некоторую окрестность точки $z=0$) \cite{[DAN2008]}.

Что касается дроби $\Theta$, то оказалось, что такое
незначительное усложнение вида НД приводит к значительно более
сильным аппроксимативным свойствам. Именно, в \cite{[DAN2006],DanChu2011} показано, что любая рациональная
функция $R$ степени $m$ при каждом натуральном $q$ может быть
приближена некоторой дробью $\Theta$, в которой все $n_k$ не
превосходят $qm$, так, что 
$$
|\Theta(z) - R(z)|\le 2 e^{|R(z)|}\frac{|R(z)|^{q+1}}{q!}
\qquad\hbox{при условии}\qquad |R(z)|\le \frac{q}{5}.
$$
Отсюда видно, что в общем случае скорость аппроксимации дробями
$\Theta$ гораздо выше, чем многочленами (той же степени, что и
$\Theta$), и весьма близка по порядку к скорости аппроксимации
рациональными функциями общего вида \cite{[DAN2006]}.

\subsection{Наилучшие приближения на отрезке действительно оси}
\label{paragraph-1-sement} Рассматривалась задача о наилучшем
приближении вещественных непрерывных функций на отрезках
действительной оси посредством вещественнозначных НД. В работах
\cite{[DAN-KOND2],[KOM3],[Kom2015-MZ]} показано, что НД $\rho_n$
наилучшего равномерного приближения вещественной константы при
достаточно больших $n$ единственна, имеет порядок $n$ и
характеризуется чебышевским альтернансом, состоящим из $n+1$ точек
отрезка. Этот результат вполне аналогичен критерию наилучшего
приближения констант посредством многочленов.

Однако, при аппроксимации произвольной непрерывной функции
возможны существенные различия. Так, в отличие от случая
многочленов, НД наилучшего приближения может быть не единственна;
соответствующие примеры были построены в
\cite{[DAN-KOND2],[Kom2011],[KOM2]}. Эти примеры, кроме того,
показывают, что, вообще говоря, не существует прямой связи между
альтернансом и наилучшим приближением посредством НД.

Тем не менее, для достаточно широкого класса непрерывных функций
установлен \cite{[Kom2014-1], [Kom2015-IzvRAN]} следующий аналог
теоремы Чебышева об альтернансе:

{\it НД наилучшего приближения
характеризуется альтернансом из $n+1$ точек отрезка $[-1,1]$ при
условии, что ее порядок равен $n$ и все полюсы лежат вне
замкнутого единичного круга; при этом, дробь наилучшего
приближения единственна.}

Этот критерий усиливает предшествующий
результат Я.~В.~Новака \cite{[NOV-diss]}, где все полюсы
предполагались вещественными попарно различными, а вопрос о
единственности не затрагивался.

Недавно в \cite{[Kom-2016],[Kom2017-IzvRAN]} доказана
справедливость критерия и в более общем случае --- когда порядок
НД {\it не превосходит} $n$; кроме того, аналогичный критерий с
заметно ослабленным ограничением на полюсы построен для задачи
приближения нечетных функций на симметричных отрезках.

\subsection{Аппроксимация на неограниченных множествах}
\label{paragraph1.6} В работах В.~Ю. Про\-та\-со\-ва \cite{[PRO]},
П.~А. Бородина и О.~Н. Косухина~\cite{[KOS2],[BOR-KOS]},
И.~Р.~Каюмова~\cite{[KAY2011],[KAY2012]},
В.~И.~Данченко~\cite{[DAN2010]}  изучалось приближение посредством
НД в разных метриках на неограниченных множествах: прямых, лучах и
др. Например, установлено~\cite{[BOR-KOS]}, что каждая непрерывная
на действительной оси $\mathbb{R}$ функция $f$ с нулевым значением
на бесконечности (кратко, $f\in C_0(\mathbb{R})$) в равномерной
метрике с любой точностью приближается НД. Аналогичное утверждение
становится неверным, если вместо прямой рассматривать
неразвернутый угол; примером может служить функция
$f(z)=-\frac{1}{z-a}$, где $a$ принадлежит внутренности угла
\cite{[BOR-KOS]}.

Кроме того, в работе \cite{[BOR-KOS]} доказана следующая теорема о
существовании функции $f$ с наперед заданными наилучшими
приближениями $\mathcal{R}_n(f,\mathbb{R})$ посредством НД порядка
не выше $n$ на $\mathbb{R}$:

{\it Для любой числовой последовательности
$\{d_n\}_{n=0}^{\infty}$, строго убывающей к нулю, существует
функция $f\in C_0(\mathbb{R})$ такая, что
$\mathcal{R}_n(f,\mathbb{R})=d_n$, $n=0,1,2,\dots$}.

Это аналог известной полиномиальной теоремы С.\,Н.\,Бернштейна
(см. \cite{Timan}). Однако в отличие от полиномиального случая здесь
условие строгого убывания существенно: доказано, что не существует
функции $f\in C_0(\mathbb{R})$, для которой
$\mathcal{R}_0(f,\mathbb{R})=\mathcal{R}_1(f,\mathbb{R})>0$,
$\mathcal{R}_2(f,\mathbb{R})=\ldots=
\mathcal{R}_n(f,\mathbb{R})=\ldots=0$.

Что касается интегральных пространств, то в \cite{[BOR2009]}
показано, что при $p\in [2,\infty)$ НД всюду плотны в
$L_p({\mathbb R}^+)$ на действительной полуоси $x\ge 0$ (и что для
любой функции $f\in L_p({\mathbb R}^+)$, $1<p<\infty$, при всяком
$n\ge 1$ НД порядка $\le n$ наилучшего в $L_p({\mathbb R}^+)$
приближения существует, но, вообще говоря, неединственна).

Однако в случае всей оси класс функций, аппроксимируемых
посредством НД в $L_p(\mathbb{R})$, $p>1$, резко сужается
\cite{[PRO]}, в частности, он состоит из тех и только тех функций,
которые представляются в виде сходящихся к ним в $L_p(\mathbb{R})$
рядов НД.

Все эти результаты указывают, в частности, на своеобразную
нелинейность процесса аппроксимации на неограниченных множествах и
тесную связь с их геометрией.

\subsection{Ряды НД}
\label{paragraph1.7} Работы по аппроксимации на неограниченных
множествах в немалой мере способствовали возникновению теории
рядов НД. Точнее, были найдены различные условия и критерии
сходимости в $L_p({\mathbb R})$ бесконечных НД
$\rho_{\infty}=\lim_{n\to\infty}\rho_n$ в терминах
последовательностей их полюсов $z_k$ (см.
\cite{[PRO],[DAN2010],[KAY2011],[KAY2012],[KAYva]}). Например, при
$1<p<2$ и $z_k\in \mathbb C^+$ такая сходимость эквивалентна
следующему условию из \cite{[DAN2010]} для частичных сумм $\rho_n$
вида (\ref{(*1)}) (см. \S\ref{paragraph3.2}):
\begin{equation}
\sum_{k=1}^{n}|\rho_n(\overline{z_k})|^{p-1}\le
A(\{z_j\}_{j=1}^{\infty})<\infty\qquad \forall n.
  \label{(*777)}
\end{equation}

\subsection{Обобщения НД}
\label{paragraph1.11}

Естественным обобщением НД являются так называемые
\textit{$h$-суммы} и \textit{амплитудно-частотные суммы},
соответственно имеющие вид
$$
\sum_{k=1}^n \lambda_k h(\lambda_k z) \quad\hbox{и}\quad
\sum_{k=1}^n \mu_k h(\lambda_k z), \qquad z,\lambda_k,\mu_k\in
\mathbb{C},
$$
где $h$ --- аналитическая в окрестности начала функция. В виде
$h$-суммы представляются НД, если взять, например, $h(z)=1/(z-1)$,
$\lambda_k=1/z_k$. Эта связь явилась одной из мотивировок
использования $h$-сумм в качестве аппарата аппроксимации
\cite{[DAN2008]}.  Такие суммы эффективно использовались в
численном анализе (А.~В.~Фрянцев \cite{[FRN1],[FRN2],[FRN3]},
П.~В.~Чунаев и В.~И.~Данченко \cite{Chu2010,Chu2012,DanChu2011} и
др.). Так, в работе \cite{DanChu2011} был предложен метод
экстраполяции функций $h$ их $h$-суммами. Ряд окончательных
результатов о скорости и области сходимости экстраполяционных
процессов установлен в \cite{Chu2012}; там также на примерах
продемонстрированы некоторые преимущества метода $h$-сумм перед
классическими полиномиальными методами.

Амплитудно-частотные суммы применялись в
\cite{[DanChu2013],[DanChu2014],[DanChu2013-2],[DanChu2016]} для
$2n$-кратной Паде-интерполяции (в точке $z=0$) индивидуальных
аналитических функций, а также в качестве операторов численного
дифференцирования, интерполяции, экстраполяции и др., действующих
на определенных классах функций. Оказалось, что возможность построения нужной
амплитудно-частотной суммы обусловлена разрешимостью
ассоциированной с ней задачи дискретных моментов
$$
\sum\nolimits_{k=1}^n\mu_k \lambda_k^m=\alpha_m, \qquad
m=0,\ldots,2n-1,
$$
относительно неизвестных $\mu_k $ и $\lambda_k$ с заданными
правыми частями $\alpha_m$. Напомним, что задачам дискретных
моментов посвящены классические труды Прони, Сильвестра,
Рамануджана и работы многих современных авторов
(см.~\cite{Prony,Sylvester,Ramanujan,Lyubich,Lyubich2,Kung1}). Эти
задачи тесно связаны с ганкелевыми формами, ортогональными
многочленами, квадратурными формулами Гаусса и аппроксимациями
Паде.

В нашем случае особую трудность и интерес представляет случай
несовместных задач дискретных моментов, когда соответствующая им
амплитудно-частотная сумма не может быть построена. Для
преодоления этой трудности в \cite{[DanChu2014],[DanChu2016]} предложен метод
аналитической регуляризации амплитудно-частотной суммы, состоящий
в том, что к ней добавляется определенный бином вида
$p\,z^{n-1}+q\, z^{2n-1}$. Оказывается, что  правильный выбор
параметров $p$ и $q$ приводит к новой ассоциированной задаче
дискретных моментов, которая уже регулярно разрешима, а в ряде
прикладных задач допускает весьма простой явный вид решений. В
результате регуляризации соответствующие <<подправленные>>
интерполяционные формулы с $n$ узлами $\lambda_k z$ становятся
точными на многочленах степени $2n-1$, что в два раза выше порядка
обычной $n$-узловой интерполяции на основе многочленов Лагранжа и
других сходных аппаратов.

При работе с тригонометрическими многочленами
$$
T_n(t)=\sum\nolimits_{k=1}^n\tau_k(t), \quad \tau_k(t):=a_k\cos
kt+b_k\sin kt
$$
и сходящимися тригонометрическими рядами, вместо
амплитудно-частотных сумм удобно применять их тригонометрический
аналог --- \textit{амп\-ли\-туд\-но-фа\-зо\-вые сум\-мы}
\cite{[Dan-Vas1],[Dan-Vas2],[Dan-Vas3],[Dan-Vas4]}. Слагаемые
таких сумм суть многочлены, {\it подобные}  $T_n(t)$: они
получаются из $T_n(t)$ умножением на вещественную константу $X$ и
сдвигом на вещественную фазу $\lambda$, т.е. $T_n(t)\to X\cdot
T_n(t-\lambda)$. Показано, например, что существует
амп\-ли\-туд\-но-фа\-зо\-вая сумма, выделяющая гармонику
$\tau_{\mu}$ заданного порядка $\mu$:
\begin{equation}
\tau_{\mu}(t)=\sum_{k=1}^mX_k\cdot T_n(t-\lambda_k),\qquad
m=m(n,\mu),
  \label{AFO1}
\end{equation}
причем параметры $X_k$, $\lambda_k$ найдены явно и не зависят от
многочлена $T_n$, а зависят лишь от $\mu$, $n$. Получены
аналогичные формулы для гармоник на достаточно широком классе
сходящихся тригонометрических рядов.

Благодаря вещественности параметров в формуле (\ref{AFO1}) она
имеет простой и важный физический смысл: из стационарного сигнала
$T_{n}(t)$ выделяется гармоника $\tau_{\mu}(t)$ наложением не
более $m$ подобных сигналов, отличающихся лишь амплитудами и
начальными фазами. Выделение гармоник наложением сигналов (без
использования промежуточных спектральных замеров) позволяет
эффективно применять амп\-ли\-туд\-но-фа\-зо\-вые суммы для оценок
гармоник тригонометрических полиномов.

Отметим еще, что в отличие от ряда спектральных методов в
конструкции амп\-ли\-туд\-но-фа\-зо\-вых сумм интегральные
аппараты не применяются; формула (\ref{AFO1}) имеет чисто
арифметический характер, из нее, например, при $a_0=0$ получаются
арифметические формулы для коэффициентов Фурье:
$$
a_{\mu}=\sum_{j=1}^{m} X_j\cdot T_n\left(-\lambda_j\right),\qquad
b_{\mu}=\sum_{j=1}^{m} X_j\cdot
T_n\left(\frac{\pi}{2\mu}-\lambda_j\right).
$$

\section{ Задача Горина и родственные вопросы}
\label{paragraph2}
\subsection{Оценки мнимых частей полюсов НД, нормированных в метрике
$L_p(\mathbb{R})$,  $1<p\le \infty$} \label{paragraph2.1}
Через $H_p=H_p(\mathbb{C}^{+})$ обозначим пространство Харди
аналитических на $\mathbb{C}^{+}$ функций~$f$,
$$
\|f\|_{H_p}:=\sup_{y>0}
\left(\int\nolimits_\mathbb{R}|f(x+iy)|^p\,dx\right)^{1/p}.
$$
При вещественных $x$ через $f(x)$ обозначим некасательные угловые
пределы функции $f\in H_p$ со стороны $\mathbb{C}^{+}$. Хорошо
известно \cite{[KUSIS],[GARNET]}, что $f(x)$ существуют почти
всюду на $\mathbb{R}$, $f(x)\in L_p$ и $\|f\|_{p}=\|f\|_{H_p}$.
Здесь и далее применяются обозначения $L_{p}=L_{p}(\mathbb R)$,
$\|\cdot\|_p=\|\cdot\|_{L_{p}}$.

Переформулируем задачу Е.~А.~Горина в несколько более общем виде. Через
${HL}^{p}_{n}$ обозначим класс всех функций вида $\rho_{n}+f$, где
$\rho_{n}$ --- НД вида (\ref{(*1)}), все конечные полюсы которой лежат
на $\mathbb{C}^+$, а $f\in {H}_{p}$. Ясно, что компоненты
$\rho_{n}$ и $f$ определяются однозначно по $\rho_{n}+f\in
{HL}^{p}_{n}$. Положим
\begin{equation}
d_n(\mathbb R,p)=\inf\left\{Y(\rho_n):\quad \|\rho_{n}+f\|_{p}\le
1, \quad \rho_{n}+f\in {HL}^{p}_{n}\right\},\qquad 1<p\le\infty.
  \label{(*7)}
\end{equation}
Мы сохраняем обозначение (\ref{(*2)}), хотя здесь оно имеет
несколько иной смысл: поскольку класс ${HL}^{p}_{n}$, очевидно,
содержит все НД порядка $\le n$, то величина $d_n(\mathbb R,p)$ из
(\ref{(*7)}) не превосходит одноименной величины из (\ref{(*2)}).
По существу же оценки этих величин мало отличаются и имеют
одинаковый порядок.

Определение (\ref{(*7)}), как и в \S \ref{paragraph1.2}, можно
переписать в виде
$$
d_n(\mathbb R,p)=\inf_{\rho_n+f\in {HL}^{p}_{n}
}\left\{Y(\rho_n)\|\rho_n+f\|^{q}_{p} \right\},\qquad
p^{-1}+q^{-1}=1.
$$

При конечных $p$ величины $d_n(\mathbb R,p)$ к нулю не убывают.
Это свойство вместе с оценкой снизу для $d_n(\mathbb R,p)$
получается из соотношения двойственности \cite{[KUSIS],[GARNET]}:
\begin{equation}
\inf_{f\in H_p}\|\rho_n+f\|_{p}
 =\sup_{g\in H_q, \;\|g\|_{q}\le
1}\left\{\left|\int_\mathbb{R}\rho_n (x)g(x)\,dx\right|=2\pi
\left| \sum_{k=1}^{n}g(z_k)\right|\right\}.
  \label{(*8)}
\end{equation}
Действительно, пусть $\rho_{n}$ --- экстремальная НД (на которой
достигается (\ref{(*7)})), причем ближайший к $\mathbb R$ полюс
есть $iy_1$, $y_1>0$, так что $y_1=d_n(\mathbb R,p)$. Возьмем в
качестве пробной функции $g_h(z)\cdot\|g_h\|_{q}^{-1}$, где
$g_h(z)=(z+ih)^{-1}$, $h>0$. Тогда, очевидно,
$$
\left|\sum_{k=1}^{n} g_h(z_k)\right|\ge (h+y_1)^{-1},\quad
\|g_h\|_{q}=2\pi\,A(p)\,{h}^{-1/p},\quad A(p):
=\frac{1}{2\pi}\left(\int_{\mathbb
R}\frac{dx}{(x^2+1)^{q/2}}\right)^{1/q}.
$$
Отсюда и (\ref{(*8)}) (где левая часть равна единице) имеем
$(h+y_1)^{-1}<{h}^{-1/p}A(p)$. Взяв $h=y_1 (p-1)^{-1}$ и решив
последнее неравенство относительно $y_1$, получим
\begin{equation}
y_1=d_n(\mathbb R,p)\ge \frac{p-1}{p^{~q}}\, A^{-q}(p).
 \label{(*9)}
\end{equation}
Неравенство (\ref{(*9)}) другим способом найдено в \cite{[BOR2007]}, оно
уточняет оценку (\ref{(*3)}).

При том же условии $\inf_{f\in H_p}\|\rho_n+f\|_{p}\le 1$ приведем
более общую оценку, в которой учитывается расстояние до $\mathbb
R$ не от одного из полюсов НД $\rho_n$, а от совокупности всех ее
полюсов. Взяв в качестве пробной функции голоморфную на
$\mathbb{C}^+$ ветвь
$$
g(z)=(z+ih)^{-\alpha},\quad \alpha=\varepsilon+1/q,\quad
0<\varepsilon<{1}/{p},\quad h>0,
$$
отображающую $\mathbb C^+$ в угол $0\le\arg g(z)\le\alpha\pi<\pi$,
из (\ref{(*8)}) получим \cite{[DAN2010]}:
\begin{equation}
\sum_{k=1}^{n}|z_k+ih|^{-\alpha}\le \sec \frac{\pi
\alpha}{2}\left|\sum_{k=1}^{n} g(z_k)\right|\le \frac{1}{2\pi}\sec
\frac{\pi \alpha}{2}\|g\|_q\le A_0\cdot
\frac{\varepsilon^{-1/q}}{1-p\varepsilon}\, h^{-\varepsilon},
    \label{(*440)}
\end{equation}
с величиной $A_0=A_0(p)$. Здесь показатель степени $(-1/q)$ справа
нельзя увеличить. Отсюда получается оценка
\begin{equation}
\sum_{k=1}^{n}|z_k|^{-\varepsilon-1/q}\le A_1(p)\cdot
\frac{\varepsilon^{-1/q}}{1-p\varepsilon},
    \label{(***440)}
\end{equation}
если положить $h=1$ в (\ref{(*440)}) и использовать соотношение
$|z_k|\asymp|z_k+i|$ (которое следует из (\ref{(*9)})). Недавно
А.\,Е.\,Додонов \cite{[DODONOV]} уточнил (\ref{(***440)}): сумма
слева заменена им на
$\sum_{k=1}^{n}|z_k|^{-\varepsilon}\ln^{-(1+\varepsilon)/q}(|z_k|+1)$
при сохранении мажоранты того же порядка относительно малых
$\varepsilon$.

В \cite{[PRO]} найдено условие для мнимых частей полюсов $z_k$:
при конечных $p>1$
 $$ \sum_{k=1}^{n} |y_k|^{1-p}\le  h_p^p\,
 \left(\int_{\mathbb R}\frac{dt}{(t^2+1)^p}\right)^{-1}\|{\rm
Im}\,\rho_n\|_p^p,\qquad z_{k}=x_{k}+iy_{k}\in {\mathbb C},
 $$
где $h_p$ --- норма оператора Гильберта, равная
$\tg\frac{\pi}{2p}$ при $1<p\le 2$ и ${\rm ctg}\,\frac{\pi}{2p}$
при $p\ge 2$ (см. (\ref{(*33)})).

В случае $p=\infty$ методы двойственности не работают и оценки
$d_n(\mathbb R,\infty)$ значительно усложняются. Сформулируем
соответствующий результат, разными способами полученный в
\cite{[DAN-viniti-1],[DAN-viniti-2],[DAN1],[DAN-DAN1]}.

\smallskip

{\it Верны соотношения $d_n(\mathbb R,\infty)\asymp \ln\ln n/ \ln
n$ $($двусторонние не\-ра\-вен\-ства$)$. Более точно, при
достаточно больших $n$ эти двусторонние неравенства справедливы с
константами $1/9$ и $2$.}

\smallskip

Эта теорема легко распространяется на случай ${\mathbb R}^{m}$,
$m\ge 3$ (см. \cite{[DAN-DISS]}). Именно, пусть фиксированы
некоторая прямая $L\subset {\mathbb R}^{m}$ и множество $M_n$,
$n\in {\mathbb N}$, сумм вида
$$
\varrho_{n}(x)=\varrho_{n}(\xi_1,\ldots,\xi_n;x):=\sum^{n}_{k=1}\frac{x-\xi_{k}}{|x-\xi_{k}|^{2}},\qquad
x,\,\xi_{k}\in {\mathbb R}^{m},
$$
с нормировкой $\|\varrho_{n}\|_{C(L)}\le 1$. Тогда для точной
нижней грани $D_{n}(L)$ расстояний от всевозможных точек $\xi_{k}$
до прямой $L$ имеем $D_{n}(L)\asymp \ln\ln n/\ln n$.
Действительно, оценки сверху $D_{n}(L)$ ничем не отличаются от
приведенных выше оценок в ${\mathbb R}^{2}$. Для оценки снизу
можно провести симметризацию, добавив к множеству $\{\xi_{k}\}$
множество $\{\xi^{*}_{k}\}$, симметричное ему относительно $L$, а
затем поворотом вокруг $L$ переместить каждую пару $\xi_{k}$ и
$\xi^{*}_{k}$ в пару комплексно сопряженных точек $z_{k}$,
$\overline{z}_{k}$ некоторой фиксированной комплексной плоскости
${\mathbb C}$, содержащей $L$. Тогда $sup$-норма на $\mathbb R$ НД
с полюсами $\{z_{k}\}\cup\{\overline{z}_{k}\}$ будет ограничена
сверху удвоенной $sup$-нормой суммы $\varrho_{n}$ на $L$. Остается
воспользоваться нижней оценкой из приведенного выше результата из~\cite{[DAN-viniti-1],[DAN-viniti-2],[DAN1],[DAN-DAN1]}.

\subsection{Оценки расстояний до полуоси $\mathbb{R}^+$ от полюсов НД,
нормированных в метрике $L_p(\mathbb{R}^+)$,  $1<p\le \infty$}
\label{paragraph2.22} В случае действительной положительной
полуоси ${\mathbb R}^+$ при $1<p<2$ показано \cite{[BOR2007]}, что
полюсы нормированной в $L_p({\mathbb R}^+)$ НД $\rho_n$,
расположенные на действительной отрицательной полуоси ${\mathbb
R}^-$ (если такие имеются), отделены от нуля некоторой величиной
$A(p)>0$. Также в [37] высказана гипотеза, что при $p\ge 2$ такие
полюсы могут сколь угодно близко подступать к нулю.
В~общей постановке задача Горина для полуоси остается открытой.

Добавим (см. \cite{[Kom2017-IzvRAN]}), что в случае $p=\infty$ при
нормировке $\|2\sqrt{x}\rho_n\|_{\infty}\le1$ (более жесткой, чем
$\|\rho_n\|_{\infty}\le 1$, на бесконечности, но более слабой
вблизи нуля) и достаточно больших $n$ полюсы отделены {\it от всей
полуоси} $\mathbb{R}^+$ величиной $\alpha_n^2$, где
$\alpha_n:=(\ln\ln (2n))\cdot(9\ln (2n))^{-1}$.

\subsection{Оценки расстояний от полюсов НД до некоторых компактов}
\label{paragraph2.2}
Пусть $K$ --- компакт на $\mathbb C$, $K^+$
--- наибольшая область из дополнения к $K$, содержащая бесконечно
удаленную точку. Введем класс $SP_n(K)$, состоящий из функций
$R_n=\rho_n+f$, где $\rho_n$ -- НД (\ref{(*1)}) с полюсами в
$K^+$, а $f$
--- голоморфная и ограниченная на $K^+$ функция с $f(\infty)=0$.
При $m>0$ через ${\delta}_n(K,m)$ обозначим расстояние до $K$ от
множества всевозможных полюсов функций $R_n\in SP_n(K)$ при
условии нормировки
\begin{equation}
\|R_n\|_{\infty,K}:=\limsup_{z\to K,\;z\in K^+}|R_n(z)|\le m.
  \label{(*17)}
 \end{equation}
Здесь, в принципе, можно обойтись нормировкой
$\|R_n\|_{\infty,K}\le 1$, поскольку задача (\ref{(*17)}) сводится
к этому случаю заменой $R_n$ на $m^{-1}\,R_{n}(m^{-1}\,z)$, но с
другим компактом $m\cdot K$.

Приведем некоторые оценки величин ${\delta}_n(K,m)$. В случае
окружности $\gamma_1: |z|=1$ имеем слабую эквивалентность (т.е.
двустороннюю оценку с абсолютными постоянными) \cite{[DAN1]}:
\begin{equation}
{\delta}_n(\gamma_1,m)\asymp n^{-1}{\ln n},\qquad n\ge n_0(m).
  \label{(*18)}
 \end{equation}
Отсюда легко получаем ${\delta}_n(\gamma_1,m)\ge A(m)n^{-1}{\ln
n}$ уже при всех $n\ge 2$. Отметим, что оценка сверху в
(\ref{(*18)}) получается с помощью простого примера $Q'/Q$, где
$Q(z)={z^{n}-(n+m)/m}$, $n\ge n_0(m)$. Из (\ref{(*18)}) легко
получается аналогичная слабая эквивалентность для замкнутой
обобщенно ляпуновской кривой\footnote{гладкой кривой, модуль
непрерывности $\omega(s)$ угла наклона касательной к которой как
функции длины $s$ дуги удовлетворяет условию Дини
$\int\nolimits_{0}s^{-1}\omega(s)\,ds<\infty $} $\gamma$
\cite{[DAN2006]}:
\begin{equation}
{\delta}_n(\gamma,m)\asymp n^{-1}{\ln n},\qquad n\ge n_0(m).
  \label{(*1118)}
 \end{equation}
Действительно, пусть $z=\varphi(w)$ --- какое-либо конформное
однолистное отображение внешности $\gamma^{+}_1$ единичного круга
на внешность $\gamma^{+}$ кривой $\gamma$ с условием
$\varphi(\infty)=\infty$. Из результата С.~Е.~Варшавского
\cite{[WARSH]} следует, что $ |\varphi'(w)|\asymp 1$, $w\in
\gamma^{+}_1$, где слабая эквивалентность зависит лишь от
$\gamma$. Заметим, что
$$
\varphi'(w)\cdot
R_n(\varphi(w))=\sum_{k=1}^{n}\frac{\varphi'(w)}{\varphi(w)-z_k}
+\varphi'(w)f(\varphi(w))=\sum_{k=1}^{n}\frac{1}{w-w_k}+ F(w),
\quad  F(\infty)=0,
$$
где мы использовали то, что
$$
\sum_{k=1}^{n}\frac{\varphi'(w)}{\varphi(w)-z_k}=
\sum_{k=1}^{n}\frac{1}{w-w_k}+g(w),\qquad z_k=\varphi(w_k),\quad
g(\infty)=0,
$$
с некоторой голоморфной в $\gamma^{+}_1$ функцией $g$,
$F(w)=g(w)+\varphi'(w)f(\varphi(w))$. Отсюда и из (\ref{(*18)})
получается оценка снизу в (\ref{(*1118)}), поскольку из $
|\varphi'(w)|\asymp 1$ имеем $|w_k|-1 \asymp {\rm
dist}(z_k,\gamma)$. Оценка сверху доказывается аналогично: с
использованием обратного к $\varphi(w)$ отображения приведенный
выше пример НД для случая $\gamma^{+}_1$ <<переносится>> на случай
$\gamma$ (получается функция класса $SP_n(\gamma)$).

К случаю окружности сводится и следующая оценка:
\begin{equation}
{\delta}_n([-1,1],m)\asymp n^{-2}{\ln^2 n}, \qquad n\ge n_0(m).
  \label{(*1228)}
 \end{equation}
Действительно, пусть $R_n=\rho_n+f\in SP_n([-1,1])$,
$\|R_n\|_{\infty,[-1,1]}\le m$.  При замене Жуковского
$z=(w+1/w)/2$, $z_k=(w_k+1/w_k)/2$ (для определенности считаем,
что $|w_k|>1$) непосредственной проверкой убеждаемся в равенстве
\begin{equation}
\rho_n(z)=\frac{2w^2}{w^2-1}F\,(w),\qquad F\,(w)=\sum_{k=1}^{n}
\frac{1}{w-w_k}+\sum_{k=1}^{n}
\left(\frac{w_k}{ww_k-1}-\frac{1}{w}\right).
 \label{(*20)}
\end{equation}
Здесь $F\in SP_n(\gamma_1)$, а значит, и функция
$\sigma(w):=R_n(z(w)) (1-w^2)/(2w^2)$ также принадлежит классу
$SP_n(\gamma_1)$ и ее $\sup$-норма на $\gamma_1$ не превосходит
$m$. Из (\ref{(*18)}) имеем $\min_k |w_k|-1\ge A_1(m)\, n^{-1}{\ln
n}$. Таким образом, если выполнено (\ref{(*17)}), то полюсы НД
$\rho_n$ лежат во внешности эллипса (см. \cite{[DAN2006],DDJ})
\begin{equation}
z=\frac{1}{2}\left(a+\frac{1}{a}\right)\cos t
+\frac{i}{2}\left(a-\frac{1}{a}\right)\sin t
 \label{(*21)}
\end{equation}
с $a=1+A_1(m)\, n^{-1}{\ln n}$. Отсюда получаем нужную оценку
снизу (см. \cite{DDJ}). Оценку сверху получим, взяв для простоты
изложения $m=3$. Если в (\ref{(*20)}) числа $w_k$ являются корнями
уравнения $w^n-{\omega}=0$ с некоторым ${\omega}>1$, то после
несложных преобразований получается
\begin{equation}
\rho_n({\omega};z)=
\frac{2nw}{w^2-1}\frac{{\omega}(w^{2n}-1)}{(w^n{\omega}-1)(w^n-{\omega})},\quad
z=\frac{1}{2}\left(w+\frac{1}{w}\right).
 \label{(*23)}
 \end{equation}
Множество полюсов НД $\rho_n({\omega};z)$ лежит на эллипсе
(\ref{(*21)}) с $a=\sqrt[n]{{\omega}}$, причем два полюса, пусть
$z_{1,2}$, вещественны. Возьмем $\omega=n^2$. Несложно убедиться,
что
$$
\max_{x\in [-1,1]} |\rho_n(n^2;x)|= |\rho_n(n^2;1)|<3.
$$
Значит,
$$
{\delta}_n([-1,1],3)\le ||z_1|-1|=\frac{(a-1)^2}{2a}\sim
2\,\frac{\ln^2 n}{n^2},\qquad n\to\infty.
$$
О.~Н.~Косухин \cite{[KOS3]} оценку снизу в (\ref{(*1228)})
распространил на случай спрямляемых компактов:
 $$
{\delta}_n(K,1)\ge 4\; c(K)\frac{\ln^2 n}{n^2}(1+o(1)),\qquad
n\to\infty,
 $$
где $c(K)$ --- гармоническая (логарифмическая) емкость компакта
$K$.

В случае малых $m$ зависимость оценок от $n$ и $m$ может быть
уточнена. Например, если множество всех полюсов НД $\rho_n$ лежит
во внешности единичной окружности и $m<1/10$, то, как установлено в \cite{DDJ},
\begin{equation}
\delta_n(\gamma_1,m)\ge 0.5m^{-\frac{1}{n+1}}-1,\qquad n\ge 1.
  \label{(**18)}
 \end{equation}
Оценка дает порядок скорости удаления полюсов при $m\to 0$,
причем, как показывает пример $Q(z)={z^{n}-(n+m)/m}$, этот порядок
точный.

\subsection{Задача Гельфонда об оценке мнимых частей полюсов НД с
нормированными на $\mathbb{R}$ производными}
\label{paragraph2.3}
Задача о точном порядке величин
$$ d'_n(\mathbb
R,\infty) =\inf\left\{Y(\rho_n):\;\|\rho'_{n}\|_{\infty}=1
\right\}=\inf_{\rho_n}\left\{Y(\rho_n)\sqrt{\|\rho'_{n}\|_{\infty}}
\right\}
 $$
сформулирована А.~О.~Гельфондом \cite{[GEL]} (вместе с
аналогичными задачами с нормировками производных более высоких
порядков). Здесь первое равенство --- определение, а второе
получается из того, что для произвольной НД вида (\ref{(*1)}) при
замене $\tilde\rho_n(z)=\rho_n(z/c)/c$,
$c=\sqrt{\|\rho'_n\|_{\infty}}$, сохраняющей вид НД, имеем
$\|\tilde\rho'_n\|_{\infty}=1$ и $Y(\tilde\rho_n)=cY(\rho_n)$.

Эта задача еще не решена. В \cite{[GEL]} показано, что
$d'_n(\mathbb R,\infty)\ge \lambda _{0}2^{-n/4}$ с некоторой
абсолютной постоянной $\lambda _{0}>0$. Введем сходные величины
для НД вида (\ref{(*1)}), конечные полюсы которых лежат в ${\mathbb
C}^+$: \
$$
 d^{+}_n(\mathbb R,\infty)
=\inf\left\{Y(\rho_n):\;\|\rho'_{n}\|_{\infty}=1,\;\{
z_k\}\subset{\mathbb C}^+ \right\}.
$$
В \cite{[DAN1]}  показано, что справедливы следующие неравенства:
$$
\frac{1}{32}\frac{\ln n}{\sqrt{n}}\le d^{+}_n(\mathbb R,\infty)
\le \frac{5}{4} \frac{\ln n}{\sqrt{n}},\quad n\ge n_0.
$$
Используя последнюю оценку снизу, для общего случая можно доказать
\cite{[DANChu2013+]}, что $d'_n(\mathbb R,\infty)\succ
\sqrt{n^{-1}\ln n}$. Вопрос о точности порядка миноранты остается
открытым.

\subsection{Аналоги задачи Горина-Чебышева в случае отрезка}
\label{paragraph2.4}
Из (\ref{(*23)}), как
несложно проверить, вытекает, что значения функции
$$
\sqrt{1-x^2}\rho_n({\omega};x)={2n{\omega}}\,\frac{\sin
n\varphi}{1-2{\omega}\cos n\varphi+{\omega}^2},\qquad x=\cos
\varphi,
$$
альтернируют на отрезке $[-1,1]$, принимая в расположенных в
порядке возрастания точках $x_k=\cos \varphi_k$, для которых $\cos
(n\, \varphi_k)=2{\omega}({\omega}^2+1)^{-1}$, попеременно
максимальные и минимальные значения, равные по модулю
$q_n:=2n{\omega}({\omega}^2-1)^{-1}$ (число точек альтернанса
равно $n$). Положим
\begin{equation}
\|f\|:=\max_{x\in[-1;1]}\,|f(x)|,\quad
\|f\|^*:=\max_{x\in[-1;1]}\sqrt{1-x^2}\, |f(x)|,\quad f\in
C([-1,1]).
  \label{(*25)}
\end{equation}
В связи с указанным альтернансом возникает вопрос
\cite{[DAN2006],DAN-DOD}:
{\it являются ли НД вида $\rho_n({\omega};z)$ экстремальными в
задаче Горина-Чебышева о НД порядка $\le n$, наименее уклоняющихся
от нуля по норме $\|\cdot\|^{*}$ при условии, что расстояние от
множества их полюсов до отрезка $[-1,1]$ не превосходит заданной
величины $\delta>0$?}

При условии $\delta> \sqrt{2}-1$ положительный ответ на него
получен в \cite{[Chu2014]}: указанное наименьшее уклонение равно
норме
\begin{equation}
\|{\rho}_{n}({\omega};\cdot)\|^*
=\frac{n}{\sqrt{T_n^2(1+\delta)-1}}\sim \frac{n}{T_n(1+\delta)},
\qquad n\to \infty,\quad \delta> \sqrt{2}-1,
  \label{(*26)}
\end{equation}
при связи $2\delta+1=\omega^{1/n}+\omega^{-1/n}$, где $T_n(z)=\cos
(n\arccos z)$ --- многочлены Чебышева первого рода. При
$0<\delta<\sqrt{2}-1$ вопрос остается открытым.

Пока не решена аналогичная задача о НД, наименее уклоняющихся от
нуля по норме $\|\cdot\|$. Однако в \cite{[Chu2014]} найдены НД,
близкие к экстремальным, и, в частности, показано, что наименьшее
уклонение $e_n$ от нуля в классе всех НД порядка $\le n$, имеющих
один из полюсов в точке $x=\delta+1$, при $\delta>\sqrt{2}-1$
удовлетворяет соотношению
$$
e_n \sim \frac{2n}{T_n(1+\delta)-T_{n-2}(1+\delta)}, \qquad n\to
\infty.
$$

Отметим, что из равенства в (\ref{(*26)}) легко получается
\cite{[Chu2014]} неравенство, противоположное неравенству типа
Маркова-Бернштейна для вещественнозначного многочлена $P_n$, не
обращающегося в нуль на отрезке $[-1,1]$, с фиксированным корнем в
точке $x=1+\delta$:
$$
\|P'_n\|^*\ge \frac{n\,\min_{x\in
[-1,1]}|P_n(x)|}{\sqrt{T_n^2(1+\delta)-1}},\qquad n\ge 1.
$$
Это неравенство является асимптотически точным \cite{[Chu2014]}:
при $n\to \infty$ для многочлена $P_n(x)=T_n(x)-T_n(1+\delta)$
вместо знака $\ge$ можно поставить $\sim$.

Другой аналог задачи Горина-Чебышева для НД со {\it свободными}
полюсами --- наилучшее приближение отличных от нуля констант ---
приведен в \S\ref{paragraph8.3}.

\subsection{Аналоги задачи Золотарева о разделении компактов}
\label{paragraph2.5} Пусть $K_1$ и $K_2$ --- произвольные
непересекающиеся компакты на расширенной комплексной плоскости
$\overline{\mathbb{C}}$. Скажем, что $K_1$ и $K_2$ разделяются
некоторой рациональной функцией $R_n(z)$ (степени не выше $n\ne
0$), если
\begin{equation}
m_0:=\min_{z\in K_1}{|R_n(z)|},\quad m:=\|
R_n\|_{C(K_2)},\quad\sigma(K_1,K_2,R_n):=\frac{m_0}{m}>1.
\label{(*27)}
 \end{equation}
Здесь предполагаем, что  $m_0>0$, $m<\infty$ и допускаем случай
$m_0=\infty$ или (и) $m=0$ (когда $K_1$ или (и) $K_2$ состоит из
конечного числа точек). Задачу Е.~И.~Золотарева можно
сформулировать следующим образом: получить из условия
(\ref{(*27)}) какую-либо метрическую характеристику отдаленности
компактов друг от друга. Такая характеристика должна зависеть
только от $n$, ${m}$, $m_0$. Первоначально задача формулировалась
для двух отрезков $K_1$, $K_2$ действительной оси \cite{[ZOLOT]}.
Для произвольных компактов со связными дополнениями она
исследовалась А.~А.~Гончаром \cite{[GON4]}. В работе \cite{[GON4]}
для гармонической емкости $c({K}_1,{K}_2)$ конденсатора
$({K}_1,{K}_2)$ с условиями (\ref{(*27)}) доказано неравенство
\begin{equation} c=c({K}_1,{K}_2)\le
n\ln^{-1}\;\sigma(K_1,K_2,R_n).
 \label{(*28)}
 \end{equation}
Эта оценка точна \cite{[GON4]}: при любом малом $\varepsilon>0$
существует рациональная функция $R_n$ достаточно высокой степени
такая, что $\sigma(K_1,K_2,R_n)\ge
\exp\left(n(c+\varepsilon)^{-1}\right)$. Отсюда также следует, что
любые непересекающиеся компакты $K_1$ и $K_2$ на
$\overline{\mathbb{C}}$ со связными дополнениями разделяются
некоторой рациональной функцией достаточно высокой степени.

В некоторых случаях из (\ref{(*27)}) или (\ref{(*28)}) можно
получить оценку эвклидова расстояния ${\rm dist}(K_1,K_2)$ между
$K_1$ и $K_2$. Например, в случае вещественных отрезков
${K}_{1}=[-1,-\delta]$, ${K}_{2}=[\delta,1]$, $0<\delta<1$, из
(\ref{(*27)}) имеем оценку \cite{[ZOLOT],[GON4]}:
\begin{equation}
\delta>\exp(-\pi^2 n\cdot \ln^{-1}\sigma(K_1,K_2,R_n))
  \label{(*29)}
\end{equation}
для любой рациональной функции $R_n$, разделяющей эти отрезки
\cite{[GON4]}. В общем случае из (\ref{(*27)}) оценку ${\rm
dist}(K_1,K_2)$ получить невозможно, и применяются оценки
определенных усредненных расстояний (например, типа (\ref{(*28)})).

Однако при $R_n=\rho_n$ оценки ${\rm dist}(K_1,K_2)$ часто имеют
место, причем при ${m_0=\infty}$ получаем рассматривавшуюся в
\S\ref{paragraph2.2} задачу о расстоянии от полюсов НД до компакта
$K_2$ при условии нормировки (\ref{(*27)}).

Приведем один результат. Если $K_2$
--- континуум диаметра $\ge 1$, а $K_1$ --- произвольный компакт,
то из (\ref{(*27)}) следует \cite{[DAN2013]}:
$$
 \ln\left(1+\frac{A(m)}{{\rm dist}(K_1,K_2)}\right) <
 A_1\cdot n\, \frac{\sigma^2}{\sigma^{3/2}-1}.
$$
Заметим, что при малых положительных $\sigma-1$ отсюда получается
оценка того же порядка, что и в (\ref{(*29)}).

В \cite{[DAN2006]} рассматривался еще один вариант задачи
Золотарева. Его можно сформулировать следующим образом. При
$\delta\in(0,1/2)$ положим
$$
\lambda_n(\delta)=
\sup_{\rho_n}\left\{\frac{\min\{|\rho_n(x)|:\;x\in
[-1+\delta,-\delta]\}}{\max\{|\rho_n(x)|:\;x\in
[\delta,1-\delta]\}}\right\}
$$
где $\sup$ берется по всем НД (\ref{(*1)}), имеющим степень не выше $n$.
Требуется найти точный порядок роста величин $\lambda_n(\delta)$
при $n\to\infty$. Этот вопрос остается открытым, здесь же приведем
некоторые оценки.

Пример (\ref{(*23)}) показывает, что $\lambda_n(\delta)$ при любом
фиксированном $\delta$ растет быстрее любой степени $n^{\alpha}$,
$\alpha>1$. Действительно, пусть $\omega=n^{\alpha}$. Полюсы НД
(\ref{(*23)}) лежат на эллипсе (\ref{(*21)}) с
$a=\sqrt[n]{\omega}$ и, следовательно, они лежат в
$b_n$-окрестности отрезка $[-1,1]$, где $b_n=2\alpha n^{-1}{\ln
n}$ при $n\ge n_0$. Кроме того, из (\ref{(*26)}) имеем $
\|{\rho}_{n}(\omega;\cdot)\|^*\le 3n^{1-\alpha}$.

Далее, НД $\rho(x)=2 {\rho}_{n}(\omega;2x-1)$ обладает
аналогичными свойствами по отношению к отрезку $[0,1]$: ее полюсы
лежат в $(b_n/2)$-окрестности отрезка $[0,1]$ и $|\rho(x)|\le
6n^{1-\alpha}(\delta(1-\delta))^{-1/2}$ на $[\delta,1-\delta]$. Из
первого свойства следует, что при достаточно больших $n$ функция
$|\rho(x)|$ монотонно возрастает на $[-1,-\delta]$ и ее
минимальное значение $|\rho(-1)|>n/2$. Таким образом,
$\lambda_n(\delta)\ge A(\delta)n^{\alpha}$. Для сравнения
напомним, что аналогичная величина для рациональных функций общего
вида растет со скоростью $e^{A(\delta)n}$, $A(\delta)>0$ (см. \cite{[BULANOV]}).

\section{ Неравенства разных метрик (задача Джексона-Никольского) для~НД}
\label{paragraph4}
Результаты, связанные с аппроксимацией в $L_p$ и c рядами НД,
опираются в основном на оценки $L_p$-норм НД. Такие оценки были
найдены в \cite{[PRO],[KAY2011],[KAY2012],[KAY2012+],[KAYva]}, все
они выражаются явно через полюсы НД. Возникает интерес к задачам
другого типа --- оценкам $L_p$-норм НД через их $L_r$-нормы на
разных множествах при различных $p>1$ и $r>1$. Впервые такие
задачи на определенных классах функций (алгебраические и
тригонометрические многочлены, целые функции экспоненциального
типа и др.) рассматривали Д.~Джексон (1933) и С.~М.~Николь\-ский
(1951). Термин <<неравенства разных метрик>> предложен в
\cite{Nikolskiy}.

\subsection{Оценки на действительной оси}
\label{paragraph4.1} Если все полюсы $z_k$ функции $\rho_n$ лежат
в верхней (нижней) полуплоскости $\mathbb{C}^+$ ($\mathbb{C}^-$),
то будем писать для определенности $\rho_n(z)=\rho^+_n(z)$
($\rho_n(z)=\rho^-_n(z)$).  Положим
\begin{equation}
\mu_n(z)={\rm
Im}\,\rho_n^{+}(z)=\sum_{k=1}^{n}\frac{y_k-y}{|z-z_k|^2},\quad
\nu_n(z)={\rm
Re}\,\rho_n^{+}(z)=\sum_{k=1}^{n}\frac{x-x_k}{|z-z_k|^2}.
  \label{(*32)}
\end{equation}
Как и выше будем применять обозначения $L_{p}=L_{p}(\mathbb R)$,
$\|\cdot\|_p=\|\cdot\|_{L_{p}}$, $p\in (1,\infty]$. В основе
следующих оценок лежит метод насечки, который впервые был применен
в работе \cite{[DAN2010]}. Его суть видна из следующей простой
задачи. Пусть
$$
B_n(z)=\prod_{k=1}^n\frac{z-z_k}{z-\overline z_k},\qquad z_k\in
{\mathbb C}^+.
$$
При каждом фиксированном $\varphi\in (0,2\pi)$ уравнение
$$
B_n^2(x)=e^{i\varphi}\qquad\Leftrightarrow\qquad
e^{ia(x)}=e^{i\varphi},\qquad a(x):=4\sum_{k=1}^{n}\arg (z_k-x),
$$
$x\in \mathbb R $, имеет ровно $2n$ вещественных различных корней
(поскольку $a(x)$ монотонно возрастает от $0$ до $4\pi n$ при
возрастании $x$ от $-\infty$ до $\infty$ и, значит, $e^{i a(x)}$
непрерывно пробегает единичную окружность $2n$ раз), которые
обозначим через $t_s=t_s(\varphi)$,
$-\infty<t_1<\ldots<t_{2n}<\infty$. Будем называть их точками
насечки. Отметим равенства:
$$
2\mu_n(x)=-i\frac{B_n'(x)}{B_n(x)}=({\rm arg}\, B_n(x))',\quad
\int_{t_s}^{t_{s+1}}\mu_n(x)\,dx=\frac{\pi}{2}.
$$
Последнее равенство позволяет приближенно вычислять точки насечки,
если известны значения $\mu_n(x)$ только при $x\in {\mathbb R}$.

В \cite{[DAN2010]} получены точные квадратурные формулы для норм
$\rho_n^{+}$, $\mu_n$ в $L_2$ с квадратурными узлами в точках
насечки:
\begin{equation}
2\|\mu_n\|^2_{2}=\|\rho_n^{+}\|^2_{2}=\pi\sum\nolimits_{s=1}^{2n}\frac{\nu_n^2(t_s)}{\mu_n(t_s)}
=\pi\sum\nolimits_{s=1}^{2n}\mu_n(t_s),
      \label{(*43)}
\end{equation}
$\varphi\in (0,2\pi)$ и $t_s=t_s(\varphi)$. Отсюда легко получить
двустороннюю оценку \cite{DAN-DOD}
\begin{equation}
(2n)^{-1} \left\| \rho_n^{+} \right\|_{2}^{2}\le \pi \left\|
\rho_n^{+} \right\|_{\infty}\le 2 \left\| \rho_n^{+}
\right\|_{2}^{2}.
 \label{(*44)}
\end{equation}
Действительно, из (\ref{(*43)}) имеем
$$
\pi\sum\nolimits_{k=1}^{2n}|\nu_n(t_k)|\le
\sqrt{\pi\sum\nolimits_{k=1}^{2n}{\nu_n^2(t_k)}{\mu_n^{-1}(t_k)}}\sqrt{
\pi\sum\nolimits_{k=1}^{2n}\mu_n(t_k)}=\|\rho_n^{+}\|^2_2,
$$
откуда с учетом (\ref{(*43)}) получаем
$$
\|\rho_n^{+}\|^2_2 \le
\pi\sum\nolimits_{k=1}^{2n}|\rho_n^{+}(t_k)| \le
\pi\sum\nolimits_{k=1}^{2n}(|\nu_n(t_k)|+\mu_n(t_k)) \le
2\|\rho_n^{+}\|^2_2.
$$
Отсюда сразу следует первое неравенство в (\ref{(*44)}). Далее,
выбрав $\varphi$ так, чтобы в одной из точек $t_k(\varphi)$
функция $|\rho_n(x)|$ принимала наибольшее значение, получим и
второе неравенство в (\ref{(*44)}). Аналогичные неравенства
справедливы и для $\rho_n^{-}$. Примеры показывают, что каждый из
множителей 2 и $1/2$ в (\ref{(*44)}) нельзя заменить величиной,
меньшей единицы \cite{DAN-DOD}. С применением метода насечки
второе неравенство в (\ref{(*44)}) обобщено и доведено до точного
в недавней работе \cite{[DanChu2016+]}:
$$
\|\rho_n^{\pm}\|_{\infty} \le
\left(\frac{m_r}{\pi}\right)^{s/r}\|\rho_n^{\pm}\|_{r}^{s},\qquad
r>1,\qquad m_r\in \mathbb{N}\cap [\tfrac{1}{2}r,\tfrac{1}{2}r+1),
$$
где $r^{-1}+s^{-1}=1$ (при $r=2$ множитель 2 в (\ref{(*44)})
заменяется на 1). Равенство может достигаться при $r=2$ и $r=4$.

Оценка (\ref{(*5)}) обобщалась на $(p,r)$-метрики
\cite{DAN-DOD,[DanChu2016+]}. Наиболее точная оценка в настоящее
время имеет вид \cite{[DanChu2016+]}:
\begin{equation}
\|\rho_n\|_{p}^{q} \le 2^{\,q-s} \left(\frac{m_r}{\pi}\right)^{sq
\left(\frac{1}{r}-\frac{1}{p}\right)}(1+h_r)^{s} \|\rho_n\|_{r}^s,
\quad \frac{1}{p}+\frac{1}{q}=1,\quad \frac{1}{r}+\frac{1}{s}=1,
 \label{(*45)}
\end{equation}
где $1<r<p\le \infty$, $m_r\in \mathbb{N}\cap
[\tfrac{1}{2}r,\tfrac{1}{2}r+1)$, а $h_r$ --- норма преобразования
Гильберта. Напомним еще раз, что $h_r$ равна ${\rm tg} \tfrac{\pi}{2r}$ при $1<r\le 2$ и 
${\rm ctg}  \tfrac{\pi}{2r}$ при $2\le r<\infty$.

Заметим, что в отличие от случая многочленов в оценку
(\ref{(*45)}) сравниваемые нормы входят нелинейно. Это
естественно, поскольку при преобразовании
$\tilde\rho_n(x):=a\rho_n(ax)$, $a>0$, сохраняющем вид НД,
оцениваемые выражения зависят от $a$ линейно:
$\|\tilde\rho_n\|_{p}^{q}=a \|\rho_n\|_{p}^{q}$.

Заметим еще, что для НД неравенства типа Джексона-Никольского
содержательны без каких-либо условий на соотношение между
параметрами $p$ и $r$. Пусть $1<p\le \infty$, $1<r<\infty$,
$E$~--- произвольный ограниченный или неограниченный промежуток на
$\mathbb{R}$. Тогда, как установлено в \cite{DAN-DOD},
\begin{equation}
 \left\| \rho_n \right\|_{L_p
(E)}^{q} \le A(p)\cdot n^{q/p}\, \left\| \rho_n
\right\|_{L_{\infty}(E)}, \qquad \left\| \rho_n \right\|_{p}^{q}
\le A(p,r)\; n^{q/p}\, \left\| \rho_n \right\|_{r}^{s}.
 \label{(*446)}
\end{equation}
Отметим, что $A(p)$ не зависит от $E$. Здесь в первом неравенстве
порядок множителя $n^{q/p}$ является точным. Второе неравенство не
является точным по порядку. При $p>r$ это следует из
$(\ref{(*45)})$, но и при $p<r$ оно, по-видимому, может быть
улучшено заменой показателя степени $q/p$ на $(q/p)-(s/r)$.

\subsection{Оценки на отрезке}
\label{paragraph4.2}
Пусть все полюсы НД (\ref{(*1)})
лежат вне отрезка $[-1,1]$ симметрично относительно действительной
оси (т.е. НД вещественнозначна на ${\mathbb R}$). Тогда, как показано в \cite{DAN-DOD},
 при $p>r>1$ и достаточно больших $n\ge n_0(\|\rho_n\|^{*})$ имеем 
\begin{equation}
\|\rho_n\|\le 64\; n^{2/r}\;\|\rho_n\|_{L_r[-1,1]},\qquad
\|\rho_n\|_{L_p[-1,1]}\le A(p)\;
n^{2\left(\frac{1}{r}-\frac{1}{p}\right)}\;\|\rho_n\|_{L_r[-1,1]}.
 \label{(*46)}
\end{equation}
Обозначение норм см. в (\ref{(*25)}). Можно взять, например, $n_0(x)=10^3(x^2+1)$. Первое неравенство
точно по порядку при $r>2$: существует последовательность
вещественнозначных НД $\tilde\rho_n$ с $\|\tilde\rho_n\|\asymp 1$,
для которых выполняется противоположное неравенство (с заменой 64
на некоторую величину $A(r)>0$). Оценки (\ref{(*46)}) по форме
аналогичны оценке для многочленов (см. (\ref{(*4)})). Однако здесь
дополнительно накладывается ограничение на величину $n$. Это
связано с несовпадением <<размерностей>>: отношение левой и правой
частей (при фиксированном $n$) растет к бесконечности, когда один
из полюсов приближается к отрезку. Наиболее естественной (т.е. без
каких-либо условий на $n$) в данном случае представляется оценка
вида
$$
\|\rho_n\|^q_{L_p[-1,1]}\le A(r)\;
n^{\vartheta}\;\|\rho_n\|^{s}_{L_r[-1,1]},\quad
p^{-1}+q^{-1}=1,\quad r^{-1}+s^{-1}=1,
$$
с некоторой величиной $\vartheta=\vartheta(p,r)$ (вопрос открыт).

\subsection{Оценки на окружности}
\label{paragraph4.3}
Пусть все полюсы НД (\ref{(*1)}) лежат во
внешности окружности $\gamma_r:\;|z|=r$. Тогда при $s\in{\mathbb
N}$ справедлива полученная в \cite{[D-D3]} оценка 
$$
\|\rho_n^{(s-1)}\|_{L_{\infty}(\gamma_r)}\le
\frac{\sqrt{s}}{\sqrt{2\pi r}}\|\rho_n^{(s-1)}\|_{L_{2}(\gamma_r)}
\sqrt{2\|\rho_{n}\|_{L_{\infty}(\gamma_r)}+n}.
$$
Отсюда с помощью неравенства Гельдера можно получить аналогичную
оценку с заменой $L_2$ на $L_p$, $p>2$.

\section{ Неравенства Маркова-Бернштейна для производных НД}
\label{paragraph5}

Оценки производных функций разных классов --- одна из наиболее
известных классических экстремальных задач, возникающих, например,
в теории полиномиальных и рациональных аппроксимаций. Обширная
библиография по этой тематике имеется в \cite{Bari,[DOL1],Timan}.
Для НД оценки производных получены в \cite{[DAN2006],DAN-DOD}.
Приведем несколько таких оценок.

\subsection{Оценки на действительной оси}
\label{paragraph5.1} Пусть конечные полюсы НД $\rho_n(z)$ лежат на
$\mathbb{C}^{+}$. Тогда (см. обозначение (\ref{(*32)}))
 $$
|\rho_n'(x)|\le \left(|\rho_n(x)|+\|\rho_n\|_{\infty}\right)
\,\mu_n(x),\quad |\nu'_n(x)|\le
\left(\mu_n(x)+\|\mu_n\|_{\infty}\right) \,\mu_n(x),
 $$
 \begin{equation}
|\mu'_n(x)|\le \chi(x), \quad |\rho'_n(x)|+|\mu'_n(x)| \le 2\,
\chi(x),
 \label{(*48)}
\end{equation}
где
$\chi(x):=\left(|\nu_n(x)|+\|\nu_n\|_{\infty}\right)\,\mu_n(x)$.
При этом соотношения в $(\ref{(*48)})$ обращаются в равенства для любой НД
первого порядка $($первое -- в некоторой точке $x\in \mathbb{R}$,
а второе --- всюду на $\mathbb{R})$.

\subsection{Оценки на окружностях и отрезках}
\label{paragraph5.2} Пусть $r>0$ и полюсы НД вида (\ref{(*1)}) лежат во
внешности окружности $\gamma_r:\;|z|=r$. Тогда имеем точную оценку
 $$
\|\rho_n'\|_{C(\gamma_{r})}\le \|\rho_n\|_{C(\gamma_r)}\left(n
r^{-1}+2\|\rho_n\|_{C(\gamma_r)}\right).
 $$
В дополнение к (\ref{(*6)}) приведем другую точную по порядку
оценку из \cite{[DAN2006]}: при условии $\|\rho_n\|^{*}\le 1$ (см. обозначение
(\ref{(*25)})) имеем 
 $$
(1-x^2)|\rho_n'(x)|\le
|x\rho_n(x)|+n\|\rho_n\|^{*}(1+\varepsilon_n),\qquad
0<\varepsilon_n<12 n^{-1}\ln n\to 0.
 $$

\section{ Разделение особенностей НД}
\label{paragraph6} Пусть $K$
--- компакт на $\mathbb C$. Для класса $SP_n(K)$ функций
$\rho_n+f$ с нормировкой (\ref{(*17)}) положим
\begin{equation}
\theta_n(K,m)=
\sup\{\|\rho_n\|_{\infty,K}:\quad\|\rho_n+f\|_{\infty,K}\le m,\;\;
\rho_n+f\in SP_n(K)\}.
 \label{(RAZDEL)}
\end{equation}
Задача о разделении особенностей функций в классе $SP_n(K)$
состоит в оценке сверху величины $\theta_n(K,m)$ через $m$ и $n$,
или, другими словами, в оценке $\sup$-норм компонент $\rho_n$ и
$f$ через норму всей функции $\rho_n+f$.

Хорошо известна аналогичная задача для более широкого класса
функций, имеющих вид $r_n+f$, где $r_n$ --- рациональные функции
общего вида степени $\le n$, а $f$ --- те же, что и в $SP_n(K)$. В
этом классе через $\theta^{*}_n(K,m)$ обозначим величину,
аналогичную (\ref{(RAZDEL)}) с заменой $\rho_n$ на $r_n$. В этом
случае можно взять $m=1$, поскольку $\theta^{*}_n(K,m)$ в отличие
от $\theta_n(K,m)$ зависит от $m$ линейно.

Первоначально задача о разделении формулировалась для двух
рациональных функций $r_n$ и $f$ в случае единичной окружности
$K=\gamma_1: |z|=1$, см. \cite{[KAC]}. Из результатов \cite{[KAC]}
следует, что $\theta^{*}_n(\gamma_1,1)\le A\,n^2$. В ряде
последующих работ рассматривались более общие компакты $K$, но
оценка оставалась по порядку величины $n$ той же (А.~М.~Бочтейн,
В.~Э.~Кацнельсон, С.~И.~Пореда, Е.~Б.~Сафф, Г.~С.~Шапиро,
А.~А.~Гон\-чар, Л.~Д.~Гри\-го\-рян и др.; соответствующая
библиография имеется в \cite{[GON5],[GON6],[DAN-razdel]}). Позже
Л.~Д.~Григорян \cite{[GRIGOR]} существенно понизил порядок
мажоранты в случае гладких кривых $K=\gamma$:
$\theta^{*}_n(\gamma,1)\le A(K)\,n$. Он же показал, что порядок
оценки точный. Окончательно в \cite{[DAN-razdel]} такая оценка
установлена для произвольного компакта $K$ с устранением
зависимости $A$ от $K$.

Возникает аналогичная задача для НД. Естественно ожидать, что
порядок роста $\theta_n(K,m)$ относительно $n$ должен быть
значительно меньше порядка роста $\theta^{*}_n(K,m)$. Для
произвольных компактов $K$ вопрос оценки $\theta_n(K,m)$ остается
открытым. В простейших случаях (окружности, ляпуновские замкнутые
кривые и др.) оказалось, что порядок роста логарифмический, и этот
порядок точный. Например, в случае $K=\gamma_1$ имеем
$\theta_n(\gamma_1,m)\le 3 \,m\ln n$ при $n\ge n_0(m)$, см.
\cite{[DAN2006]}.

Аналогичную задачу об оценке компонент можно ставить на
неограниченных множествах и в разных метриках; простейший случай
см. в (\ref{(*35)}) и (\ref{(*36)}) из \S\ref{paragraph3.1}.

\section{ Интерполяция посредством НД}
\label{paragraph7}
\subsection{Задача интерполяции посредством НД. Примеры}
\label{paragraph7.1} Рассмотрим задачу построения НД $\rho_n$
порядка $\le n$, удовлетворяющей равенствам
\begin{equation}\begin{array}{r}
    \rho_n^{(s)}(\xi_j)=b_{j,s},
    \quad \xi_j,b_{j,s}\in\mathbb{C}\qquad
    (j=\overline{1,k},\ s=\overline{0,m_j-1},\ m_j\ge 1, \\
    m_1+\ldots+m_k=M),
     \label{(7.5)}
\end{array}\end{equation}
где $\xi_1,\dots,\xi_k$~--- различные узлы, $m_1,\dots,m_k$~--- их
кратности, $M$ --- суммарная кратность узлов, а $b_{j,s}$ ---
заданные значения.  В случае $M<n$ либо $M>n$ будем называть
задачу (\ref{(7.5)}) соответственно {\it неполной} либо {\it
пере\-оп\-ре\-де\-ленной} задачей интерполяции. Случай $M=n$ для
нас является основным. В этом случае, если не оговорено противное,
мы будем называть (\ref{(7.5)}) просто {\it задачей интерполяции}
(иногда для уточнения будем применять термин <<{\it полная} задача
интерполяции>>).

Если все узлы в (\ref{(7.5)}) {\it простые}, то приходим к задаче
{\it простой} интерполяции:
\begin{equation}
    \rho_n(\xi_j)=b_j, \qquad j=\overline{1,M}.
     \label{(7.3)}
\end{equation}
В противном случае (\ref{(7.5)}) называется задачей {\it кратной}
интерполяции. При интерполяции регулярной функции $f$ в
(\ref{(7.5)}) и (\ref{(7.3)}) естественным образом полагают
$b_{j,s}=f^{(s)}(\xi_j)$, $b_{j}=f(\xi_j)$.

Как известно, в классе алгебраических полиномов степени $\le n-1$
полная задача интерполяции вида (\ref{(7.5)}) всегда имеет
решение, причем единственное. Напротив, интерполяционная НД может
не быть единственной или не существовать. Проиллюстрируем это
примерами.

1) Для таблицы $\{(-1;-1),(1;1)\}$ существует бесконечно много
интерполяционных НД порядка $\le 2$ (см. \cite{[NOV-diss]}, а
также \cite{[KOND]}):
$$
\frac{1}{z},\qquad\rho_2(z;\alpha):=\frac{2z+\alpha}{z^2+\alpha
z+1}, \qquad \alpha\ne \pm 2.
$$

2) Пример {\it неразрешимой} полной задачи (\ref{(7.3)}) дает
таблица (см. \cite{[DAN-KOND2],[KOM1]})
$$(-\sqrt 2,-3\sqrt 2),\quad(-1/\sqrt 3,-\sqrt 3),\quad(0,1),\quad(1/\sqrt 3,\sqrt
3),\quad(\sqrt 2,3\sqrt 2), $$ для которой не существует
интерполяционных НД порядка $\le 5$. Действительно, если
существует такая НД $\rho_5:=Q'/Q$, то необходимо $Q'(0)=Q(0)\ne
0$, и, не нарушая общности, мы можем искать $Q$ в виде $Q(z)=1+a_1
z+\ldots+a_5z^5$. Но тогда из равенств $\rho_5(\pm\sqrt 2)=\pm 3
\sqrt 2$ и $\rho_5(\pm 1/\sqrt 3)= \pm\sqrt 3$ имеем
соответственно равенства $6a_3+4a_5=-5$ и $6a_3+4a_5=0$.
Противоречие.

Тем не менее, при определенном ограничении на суммарную кратность
$M$ узлов задача (\ref{(7.5)}) всегда (т.е. независимо от значений
$\xi_j$ и $b_{j,s}$) разрешима. А именно, справедливо утверждение \cite{[KOM1]}: \smallskip

{\it Задача $(\ref{(7.5)})$ всегда разрешима, если $M\le n-k+1$. В
частности, задача простой интерполяции всегда разрешима, если
число $k$ узлов удовлетворяет условию $2k\le n+1$. Эта оценка
неулучшаема.}\smallskip

Таким образом, при $k>1$ не все полные задачи разрешимы, а при
$k=1$ полная задача (с одним $n$-кратным узлом) всегда разрешима
({\it Паде-интерполяция}). Впервые разрешимость такой задачи была
установлена иными методами в \cite{[D-D2]} (см.
\S\ref{paragraph7.7}).

В \cite{[Kom2011],[Kom2015-IzvRAN]} рассматривался вопрос о
минимальном $M>n$, гарантирующем единственность решений
переопределенных задач (\ref{(7.5)}). Показано, что такое
минимальное $M$ равно $2n-1$. О других признаках разрешимости и
единственности решения задачи интерполяции см. в \S\ref{paragraph7.4} и \S\ref{paragraph8.6}.

\subsection{Редукция к полиномиальной интерполяции}
\label{paragraph7.2} В \cite{[KOM1]} для исследования задачи
интерполяции посредством НД предложен метод редукции к
полиномиальной интерполяции. В основе редукции лежит следующая
элементарная лемма пересчета производных (см. \cite{[KOM1]};
аналогичные формулы применялись также, например, в
\cite{[KOS2],[NOV2]}).

{\it Производные аналитической функции $u(z)$ удовлетворяют
равенствам
$$u^{(s+1)}=uF_{s}(w),\qquad w=w(z)=u'(z)/u(z),\quad
s=0,1,2,\dots,$$ где $\{F_s\}$ --- последовательность
дифференциальных операторов, определяемых рекуррентно:}
\begin{equation}
 F_0(w)=w,\qquad
F_{s+1}(w)=w F_s(w)+(F_s(w))'_z,\quad s=0,1,2,\dots.
     \label{(400)}
\end{equation}

Несложно проверить, что
$F_s(w)=w^{(s)}+{\Omega}_s(w,w',\dots,w^{(s-1)})$, где
${\Omega}_s$~--- определенный алгебраический полином степени не
выше $s+1$ по каждому аргументу, $\Omega_0=0$. Таким образом, во
всех точках аналитичности дроби $\rho_n=Q_n'/Q_n$, имеем
\begin{equation}
Q_n^{(s+1)}(z)=Q_n(z)\cdot
(\rho_n^{(s)}(z)+{\Omega}_s(\rho_n(z),\ldots,\rho_n^{(s-1)}(z)),\qquad
s\ge 0.
     \label{(200)}
\end{equation}

Существует аналогичная по виду формула обращения, выражающая
производные $\rho_n$ через производные $Q_n$ во всех точках, где
НД $\rho_n(z)$ аналитична.

Задачу (\ref{(7.5)}) можем переписать в виде задачи для
многочленов:
\begin{equation}
\left\{\begin{array}{cl}
  Q_n^{(s+1)}(\xi_j)=Q_n(\xi_j)\cdot h_{j,s},\quad
  Q_n\not\equiv 0,\
    &  h_{j,s}:=b_{j,s}+{\Omega}_s(b_{j,0},\ldots,b_{j,s-1}), \\
    &  j=\overline{1,k},\quad s=\overline{0,m_j-1},\\
    & m_1+\ldots+m_k=M.
\end{array}\right.
     \label{(100)}
\end{equation}
На (\ref{(100)}) основан критерий разрешимости полной задачи
$(\ref{(7.5)})$, см. \S \ref{paragraph7.4}.

\subsection{Дифференциальное уравнение НД}
\label{paragraph7.3}
Поскольку $Q_n^{(n+1)}\equiv 0$, то из
$(\ref{(200)})$ при $s=n$ получается обыкновенное дифференциальное
уравнение порядка $n$, которому удовлетворяет любая НД $\rho_n$
порядка $\le n$ во всех точках, отличных от ее полюсов~\cite{[KOM1]}:
\begin{equation}
y^{(n)}(z)+{\Omega}_n(y(z),\ldots,y^{(n-1)}(z))=0.
     \label{(300)}
\end{equation}
Можно показать, что дроби $\rho_n$ исчерпывают всю совокупность
решений уравнения $(\ref{(300)})$. Действительно, для
$(\ref{(300)})$ выполнены условия теоремы существования и
единственности решения задачи Коши. Предположим, что $h$~---
какое-либо решение уравнения $(\ref{(300)})$ и в точке $z=a$ это
решение аналитично. Построим НД $\rho_n$ порядка не выше $n$,
которая удовлетворяет условиям $\rho_n^{(s)}(a)=h^{(s)}(a)$,
$s=\overline{0,n-1}$. Такая дробь существует и единственна (как НД
Паде) и тоже является решением уравнения $(\ref{(300)})$, как это
уже доказано. Но тогда $h\equiv \rho_n$ ввиду единственности
решения задачи Коши.

Например, общим решением уравнения $y'+y^2=0$ является множество
НД порядка $\le 1$, а уравнения $y''+3yy'+y^3=0$ --- множество НД
порядка $\le 2$.

\subsection{Задача обобщенной интерполяции}
\label{paragraph7.4} Начнем с интерполяции в случае простых узлов.
Решение $\rho_n=Q'_n/Q_n$ полной задачи (\ref{(7.3)}) необходимо
удовлетворяет системе уравнений
\begin{equation}
    Q_n'(\xi_j)=b_jQ_n(\xi_j), \qquad j=\overline{1,n}; \qquad
    Q_n\not\equiv 0,
     \label{(7.4)}
\end{equation}
которую называют {\it задачей\footnote{Эта задача по форме
аналогична классической интерполяции Паде в случае рациональных
функций общего вида. Термин {\it обобщенная интерполяция} был
введен в \cite{[KOND.Abstr2010]}.} обобщенной интерполяции с
простыми узлами}. Легко видеть, что решение $Q_n$ системы
(\ref{(7.4)}) всегда существует. Таких решений, не получающихся
друг из друга домножением на постоянные множители, может быть
бесконечно много; будем называть их и соответствующие НД
$\rho_n=Q'_n/Q_n$ решениями задачи обобщенной интерполяции
(\ref{(7.4)}). Отметим, что решение $\rho_n=Q'_n/Q_n$ задачи
(\ref{(7.4)}) может не быть решением задачи~(\ref{(7.3)}). Дело в
том, что порождающий многочлен $Q_n$ может обращаться в нуль в
некоторых узлах $\xi_j$ одновременно с производной $Q_n'$. Такие
узлы называются {\it особыми}, а остальные --- {\it регулярными}
относительно решения $\rho_n$.  В особых узлах дробь $\rho_n$
имеет полюсы, а равенства (\ref{(7.4)}) выполняются независимо от
значений $b_j$. При этом узел, особый относительно одного решения
задачи (\ref{(7.4)}), может быть регулярным относительно другого
решения. Итак, задача (\ref{(7.3)}) разрешима, если и только если
все ее узлы регулярны относительно одного из решений $\rho_n$
задачи~(\ref{(7.4)}).

Обобщенная интерполяция (\ref{(7.4)}) с простыми узлами
систематически изучалась в \cite{[DAN-KOND1]}, но при
дополнительном условии $\deg Q_n=n$. Это ограничение может
приводить к неразрешимым задачам (\ref{(7.4)}). Были получены
некоторые критерии их разрешимости с таким ограничением в
алгебраических и геометрических терминах. Рассмотрим примеры:

1) Дроби $\rho_2(z;\alpha)$, определенные в \S\ref{paragraph7.1},
являются решениями задачи обобщенной интерполяции таблицы
$\{(-1;-1),(1;1)\}$ при всех $\alpha\in\mathbb{C}$, а вот
решениями задачи интерполяции той же таблицы в обычном смысле они
являются только при $\alpha\ne \pm 2$, поскольку $\rho_2(z;2)$ и
$\rho_2(z;-2)$ имеют полюсы в узлах интерполяции $\xi_1=-1$ и
$\xi_2=1$ соответственно.

2) Задача обычной интерполяции 5-узловой таблицы из
\S\ref{paragraph7.1} решений не имеет, однако она имеет бесконечно
много решений в обобщенном смысле (\ref{(7.4)}):
$$ \rho_5=Q'/Q, \qquad Q(z)=a_2z^2+2a_3z^3-a_2z^4-3a_3z^5, \quad
|a_2|+|a_3|\ne 0.
$$

{\bf Обобщенная кратная интерполяции посредством НД.} Аналогично
пре\-ды\-дущему определим обобщенную кратную интерполяцию таблицы
$b_{j,s}$, заменив равенства $(\ref{(7.5)})$ на $(\ref{(100)})$.
Узел $\xi_j$ будем называть {\it особым} (соответственно {\it
регулярным}) {\it относительно решения $\rho_n=Q_n'/Q_n$ задачи
обобщенной интерполяции} $(\ref{(100)})$, если $Q_n(\xi_j)=0$
(соответственно $Q_n(\xi_j)\ne 0$). Отметим, что в особых узлах
$\xi_j$ многочлен $Q_n$ имеет нули кратности не ниже $m_j+1$, а
$\rho_n$ имеет полюсы. В \cite{[KOM1]} доказано следующее
предложение: \smallskip

{\it Полная задача $(\ref{(100)})$ $($т.е. при $M=n)$ всегда
разрешима.} \smallskip

Как и в случае интерполяции с простыми узлами, решение задачи
$(\ref{(100)})$ не всегда дает решение задачи~$(\ref{(7.5)})$.
Очевидным следствием предыдущего предложения является утверждение:
\smallskip

{\it Полная задача $(\ref{(7.5)})$ разрешима, если и только если
среди решений обобщенной задачи интерполяции $(\ref{(100)})$
найдется решение $\rho_n$, для которого все узлы регулярны.}\smallskip

В работе \cite{[KOM1]} с помощью указанной выше редукции получен
критерий однозначной разрешимости задачи обобщенной интерполяции
$(\ref{(100)})$, а также алгебраический критерий единственности
решения задачи (\ref{(7.5)}). Для получения условий единственности
решения задачи интерполяции применялись и другие методы (см. \S
\ref{paragraph8.6}).

\subsection{Интерполяция вещественных констант}
\label{paragraph7.6} Решение $\rho_n=Q_n'/Q_n$ задачи обобщенной
интерполяции констант при любом наборе $n$ узлов $\{\xi_k\}$
всегда существует, единственно и имеет порядок, равный $n$, см.
\cite{[KOND],[DAN-KOND2]}. Таким образом, в смысле обобщенной
интерполяции имеем
\begin{equation}
\rho_n(z)-c=-c\frac{\Pi_n(z)}{Q_n(z)},\qquad
\Pi_n(z):=\prod_{k=1}^n (z-\xi_k),\qquad
\rho_n(z)=\rho_n(c,\{\xi_k\};z).
     \label{(777.4)}
\end{equation}
Можно выписать и явный вид порождающего многочлена:
\begin{equation}
Q_n(z)=\sum\nolimits_{k=0}^{n}c^{-k}\Pi_n^{(k)}(z),
     \label{(777.4+)}
\end{equation}
что значительно упрощает оценки границ корней этого многочлена и
погрешности интерполяции. Наиболее полно изучена задача об
интерполяции вещественных констант на отрезках вещественной оси;
см. \cite{[KOND],[DAN-KOND2],[KOM3], [Kom2015-MZ]}. Не нарушая
общности, ограничимся изложением результатов в случае интерполяции
положительных констант по узлам, лежащим на отрезке $[-1,1]$.

Первый результат для случая интерполяции констант $c\in(0,15/31)$
по чебышевскому набору узлов
$t_k=\cos\left((2k-1)\pi/(2n)\right)$, $k=\overline{1,n}$, получен
в \cite{[KOND]}, где доказано, что $\rho_n(c,\{t_k\};z)$ не имеет
полюсов в круге $|z|\le 1$, и дана оценка
$$
\|\rho_n(c,\{t_k\};\cdot)-c\|_{C([-1,1])} \le
\frac{c}{2^{2n-1}n!}\cdot\frac{1-c}{1-2c},\qquad n\ge 2.
$$
К настоящему времени наиболее сильный результат получен в
\cite{[Kom2015-MZ]}: \smallskip

{\it При $c>0$ и $n>8c+1$ для любого набора узлов $\xi_k\in
[-1,1]$ имеем точную по порядку двустороннюю оценку
$$
 \frac{c^{n+1}}{2^{n-1}n!e^{2c}}\le
\|\rho_n(c,\{\xi_k\};\cdot)-c\|_{C([-1,1])}\le\frac{2^n
    c^{n+1}e^{2c}}{n!v_n(c)}
$$
с определенным $v_n(c)>0$, $\lim_{n\to\infty}v_n(c)=1$. При этом
порядок по $n$ верхней оценки достигается, например, при
интерполяции с $n$-кратным узлом $\xi_1=1$, а порядок нижней ---
при интерполяции по чебышевскому набору узлов
$\{t_k\}$.}\smallskip

Отметим еще, что явные интерполяционные формулы были получены в
\cite{[Kom2012]} методом разностных уравнений для рациональных
функций вида
 $$
az+b,\qquad
 \frac{a_1z+b_1}{(a_2z+b_2)z}, \qquad |a_1|+|b_1|\ne
0, \qquad |a_2|+|b_2|\ne 0
$$
(часть этих результатов была независимо получена Е.\,Н.~Кондаковой
\cite{[KOND-diss]}). Оценки погрешности интерполяции получены к
настоящему моменту только в исключительных случаях, к которым,
помимо интерполяции констант, относится Паде-ин\-тер\-по\-ля\-ция.

\subsection{Паде-интерполяция }
\label{paragraph7.7} Пусть фиксирована аналитическая в некоторой
окрестности начала функция
\begin{equation}
f(z)=f_0+f_1z+f_2z^2+\ldots.
 \label{7.1-f}
\end{equation}
Скажем, что НД $\rho_{\nu}$ вида (\ref{(*1)}) порядка $\nu\le n$
является \textit{НД Паде} $n$-кратной интерполяции  (с узлом
$z=0$) функции $f$, если
\begin{equation}
\label{7.1-Pade} f(z)-\rho_{\nu} (z)=O(z^n)\qquad\hbox{при}\qquad
z\to 0.
\end{equation}

Задача Паде для НД всегда имеет решение и притом единственное
\cite{[D-D2]}. В \cite{[D-D2]} предложен следующий метод
построения НД Паде. При натуральных $m$ обозначим через
\begin{equation}
S_m=S_m(z_1^{-1},\ldots,z_{n}^{-1}):=\sum\nolimits_{k=1}^n
z_k^{-m}=- \tfrac{1}{(m-1)!}\rho_{n}^{(m-1)}(0)
 \label{7.1-Sm}
\end{equation}
\textit{степенные суммы} чисел $z_k^{-1}$. Условие
(\ref{7.1-Pade}), как легко видеть из $(\ref{7.1-Sm})$,
равносильно разрешимости относительно $\lambda_k=z_k^{-1}$ системы
$S_m=-f_{m-1}$, $m=\overline{1,n}$, а последняя, как хорошо
известно \cite[Гл. 11, \S52-53]{Kurosh}, всегда имеет
(единственное) решение совпадающее с корнями многочлена
 \begin{equation}
 \label{7.1-Tn}
\mathcal T_n(\lambda)={\lambda}^n-\tau_1{\lambda}^{n-1}+
\tau_2{\lambda}^{n-2}+\ldots+(-1)^n\tau_n,
\end{equation}
коэффициенты которого находятся по  \textit{рекуррентным формулам
Ньютона}
$$
\tau_1=-f_{0},\quad \tau_m={(-1)^{m}}{m^{-1}}
\left(f_{m-1}+\sum_{j=1}^{m-1} (-1)^{j}\,f_{m-j-1}\tau_j\right),
\quad m=\overline{2,n}.
$$
Поэтому все конечные (и ненулевые) полюсы $z_k$ искомой НД Паде
суть величины, обратные ненулевым корням $\lambda_k$ (обозначим их
число через $\nu$) многочлена (\ref{7.1-Tn}). Если $\lambda_k=0$,
то полагаем $z_k=\infty$, и соответствующее слагаемое в $\rho_\nu$
исчезает. В результате получается искомая НД Паде порядка $\nu$.

О.~Н.~Косухин в \cite{[KOS2]} предложил другой подход к построению НД Паде; его 
суть в следующем. Для функции (\ref{7.1-f})
рассмотрим интеграл $\alpha(z):=\int_0^z f(\zeta)\,d\zeta$, где
интегрирование ведется по любому спрямляемому пути в некоторой
окрестности начала. Оказывается, НД Паде для функции
$f$ совпадает с логарифмической производной частичной суммы
$\sum_{k=0}^n c_kz^k$ ряда Маклорена функции $\exp(\alpha(z))$.
Благодаря такой конструкции в \cite{[KOS2]} был получен ряд
интересных результатов об области сходимости последовательности
интерполяционных НД Паде $\{\rho_n\}$ к $f$ (эта область, вообще
говоря, шире чем круг сходимости ряда Тейлора для функции $f$).
Кроме того, в \cite{[KOS2]} получены оценки остаточного члена
интерполяции НД Паде для аналитических в единичном круге функций
$f$ класса Харди $H_1$ (с конечной нормой
$\sup_{0<r<1}\int_0^{2\pi }|f(re^{i\varphi})|\,d\varphi$).

В работе \cite{DanChu2011} разработан еще один метод построения НД
Паде --- в виде интеграла Эрмита
\begin{equation}
\rho_{\nu}(z)= \frac{1}{2\pi i}
\int_{\gamma}f(\zeta)\;\frac{R_n(\zeta)-R_n(z)}{(\zeta-z)\,R_n(\zeta)}\,d\zeta,
\qquad z\in G(\gamma),
  \label{7.100-rem}
\end{equation}
где $\gamma$ --- спрямляемый жорданов контур, содержащий внутри
себя точку $z=0$, функция $f$ аналитична на замыкании области
$G(\gamma)$, ограниченной контуром $\gamma$,
$$
R_n(z)=\frac{z^n}{Q(z)},\quad
Q(z)=\prod_{k=1}^{{\nu}}(z-z_k)=z^{\nu}+q_{{\nu}-1}z^{{\nu}-1}+\ldots+q_0,
$$
а значения $\nu$, $z_k$ определены выше вокруг формулы
(\ref{7.1-Tn}). Из (\ref{7.100-rem}) имеем явный вид остаточного
члена
\begin{equation}
f(z)-\rho_{\nu}(z)=\frac{1}{Q(z)}\sum_{k=n}^{\infty}z^{k}\sum_{m=0}^{\nu}\,q_m
f_{k-m}, \quad |z|\le\min_{\zeta\in\gamma}|\zeta|.
 \label{7.1-rem}
\end{equation}
Приведем полученную в \cite{DanChu2011} оценку остаточного члена в
случае, когда коэффициенты $f_m$ по модулю ограничены членами
некоторой геометрической прогрессии.\smallskip

{\it Если $|f_{m-1}|\le a^m$ для всех $m$  при некотором $a>0$, то
\begin{equation}
|f(z)-\rho_{\nu}(z)|\le\frac{a}{1-a|z|}\frac{|z|^n}{r^{n}}
\left(\frac{1-\varepsilon_n+ar}{1-\varepsilon_n-ar}\right)^n
\ln\frac{er}{r-|z|},\quad |z|<r<\frac{1-\varepsilon_n}{a},
 \label{4.1-eps}
 \end{equation}
где $\varepsilon_n$ удовлетворяют соотношениям
$\varepsilon_n^2-(1-\varepsilon_n)^{n+1}=0$, $\varepsilon_n\sim
2n^{-1}\ln n$ при $n\to \infty$.}\smallskip

Отметим, что (\ref{4.1-eps}) опирается на следующую точную по
порядку величины $a$ оценку, представляющую и самостоятельный
интерес \cite{Chu2010}:\smallskip

{\it Если для чисел $\lambda_k$ их степенные суммы вида
$(\ref{7.1-Sm})$ удовлетворяют неравенствам $|S_m|\le a^m$,
$m=\overline{1,n}$, то $|\lambda_k|<(1-\varepsilon_n)^{-1}a$,
$k=\overline{1,n}$.}

\section{Аналог теоремы Мергеляна.
Скорость равномерной аппроксимации посредством НД} Пусть
$K\subset{\mathbb C}$ --- компакт со связным дополнением. Через
$CA(K)$ обозначим класс непрерывных на компакте $K$ функций,
аналитических в его внутренних точках. Как говорилось в
\S\ref{paragraph1.5}, такие функции сколь угодно точно
приближаются в $C(K)$ посредством НД (аналог теоремы Мергеляна)
\cite{[D-D1],[D-D2]}. Возникает естественный вопрос о скорости
аппроксимации функций $f\in CA(K)$ и ее связи со свойствами~$f$.

Обозначим через ${\mathcal R}_{n}(f,K)$ и ${E}_{n}(f,K)$
наименьшие равномерные уклонения на $K$ функции $f$ от множества
НД порядка не выше $n$ и многочленов степени не выше $n$
соответственно.

Пусть $K$ --- спрямляемый компакт (см. \S\ref{paragraph1.3}) и
$f\in CA(K)$. В \cite{[D-D1],[D-D2]} показано, что
$$
{\mathcal R}_{[n\ln (1/{E}_{n}(f,K))]}(f,K)<A(f,K)\cdot
{E}_{n}(f,K),\qquad n\ge n_0(f,K).
 $$

При фиксированном $b\in K$ положим
$\alpha(f;b,z)=\int_{b}^{z}\,f(t)\,dt$, где интеграл берется по
лежащей на $K$ спрямляемой кривой, соединяющей точки $b$ и $z\in
K$. О.~Н.~Косухин в \cite{[KOS2]} установил на спрямляемых
компактах слабую эквивалентность:
\begin{equation}\label{KOS-N}
{\mathcal R}_{n+1}(f,K)\asymp {E}_{n}(fe^{\alpha(f;b,\cdot)},K).
\end{equation}
Точнее, если каждую точку компакта $K$ можно соединить с точкой
$b$ спрямляемым лежащим в $K$ контуром длины $\le d$, то
$$
A_1(d\|f\|_{C(K)}){\mathcal R}_{n+1}(f,K)\le
{E}_{n}(fe^{\alpha(f;b,\cdot)},K)\le A_2(d\|f\|_{C(K)}){\mathcal
R}_{n+1}(f,K),
$$
где $A_1(x)=1/(2(1+x)e^x)$, $A_2(x)=(1+2xe^x)e^x$.

Этот результат и ряд разработанных О.~Н.~Косухиным методов позволили ему
получить для НД аналоги классических полиномиальных теорем
Д.~Джексона, С.~Н.~Бернштейна, А.~Зигмунда, В.~К.~Дзядыка,
Дж.~Л.~Уолша \cite{[KOS2]}. Приведем полученный в \cite{[KOS2]} аналог 
классической теоремы Джексона-Корней\-чу\-ка об оценке уклонений
$E_n(f,K)$ через модуль непрерывности $\omega
(f,K;\delta)$.\smallskip

{\it Если $K=D$ --- замкнутый круг или отрезок прямой диаметра
$2d$, то для $f\in CA(D)$ имеем}
$$
 {\mathcal R}_{n}(f,D)\le A(d\|f\|_{C(D)})\left(\omega \left(f,D;\frac {\pi
d}{n}\right)+ \left(\frac {\pi d}{n}\right)\|f\|_{C(D)}^2\right).
$$

Здесь, в отличие от полиномиального случая, нельзя отбросить второе
слагаемое в мажоранте, так как уже для $f\equiv {\rm const}\ne 0$
имеем $ ~\omega(f,D; \pi d/n)=0$ и в то же время ${\mathcal
R}_n(f,D)\ne 0$.

Этот результат в \cite{[KOS2]} обобщен на функции с $f^{(s)}\in
CA(D)$, $s=1,2,\ldots$. Приведем еще результат О.~Н.~Косухина
\cite{[KOS2]} и результат Я.~В.~Новака \cite{[NOV2]} (аналоги
теорем Бернштейна), показывающие, насколько в некоторых задачах
сходны аппроксимативные свойства НД и многочленов. \smallskip

{\it Пусть $D$ --- замкнутый круг, $\beta\in (0,1)$. Функция
$f^{(s)}\in CA(D)\cap{\rm Lip}^{\beta}(D)$ тогда и только тогда,
когда ${\mathcal R}_{n}(f,D)\le C/n^{s+\beta}$.} \smallskip

{\it Непрерывная на $[-1,1]$ функция $f$ принадлежит классу
$C^n([-1,1])$, если и только если существует функция $F\in
C([-1,1])$ такая, что
\begin{equation}
{\mathcal{R}_n(f,[x_1,x_2])}{(x_2-x_1)^{-n}} \rightarrow F(x),
\qquad -1\le x_1<x<x_2\le 1,
 \label{(*NOV)}
\end{equation}
равномерно по $x$ при $x_1, x_2 \to x$. Если последнее условие
выполняется, то имеет место тождество}
$$
n!2^{2n-1}F(x)\equiv
\left|\frac{d^n}{dt^n}\bigg(f(t)e^{\alpha(f;x,t)}\bigg)_{t=x}\right|.
$$

Из этого утверждения можно усмотреть порядок скорости наилучшего
приближения малых констант на отрезке $[-1,1]$. Действительно,
пусть $f(x)=c>0$. Тогда $F(x)\equiv c^{n+1}2^{1-2n}/n!$. При
$x_1=-\varepsilon$, $x_2=\varepsilon$ и $\varepsilon\to 0$ из
(\ref{(*NOV)}) для НД $\tilde
\rho_n(\varepsilon;x)=\varepsilon\rho_n(\varepsilon;x
\varepsilon)$, где $\rho_n(\varepsilon;\cdot)$
--- НД наилучшего приближения константы $c$ на
$[-\varepsilon,\varepsilon]$, имеем (ср. с (\ref{8.4.delta})),
$$
 \max_{x\in[-1,1]}|\tilde \rho_n(\varepsilon;x)-\tilde c| = \frac{\tilde
c^{n+1}}{2^{n-1}n!}+o(\varepsilon^{n+1}),\qquad \tilde
c=\varepsilon c.
$$

Пользуясь методами работы \cite{[KOS1]}, Я.~В.~Новак (см.
\cite{[NOV*],[NOV-diss]}) распространил слабую
эквивалентность (\ref{KOS-N}) на пространства $L_p([-1,1])$, $1\le
p<\infty$. Опираясь на этот результат, он получил, в частности,
описание класса $C^n([-1,1])$ в терминах скорости локальных
приближений посредством НД в метрике $L_p([-1,1])$ при $p>1$ (этот
результат обобщает его теорему, приведенную выше, и является
аналогом теоремы А.\,Н.\,Морозова \cite{[MOROZOV]} для
алгебраических многочленов).

\section{Равномерная аппроксимация на отрезке действительной оси}
\label{paragraph7--}

В этом параграфе рассматривается задача наилучшего равномерного
приближения непрерывных вещественных функций $f$ на отрезке
$K\subset \mathbb{R}$ {\it вещественнозначными} НД. Не нарушая
общности, можем считать $K=[-1,1]$. Для этого случая обозначим
через ${\mathcal R}_n(f)$ наименьшее уклонение функции $f$ от
множества вещественнозначных НД порядка $\le n$, а через
$\rho^{*}_n(f;x)$
--- вещественнозначную НД наилучшего приближения (несложно
показать, что эта НД существует \cite{[DAN-KOND2]}).

\subsection{Неединственность и другие особенности НД наилучшего приближения}
\label{paragraph8.2}  Как уже говорилось, ряд аппроксимативных
свойств НД принципиально отличаются от свойств многочленов. Так,
существуют непрерывные функции $f$, для которых

$\rm(a)$ НД $\rho_n^*(f;x)$ неединственна;

$\rm(b)$ разность $\rho_n^*(f;x)-f(x)$ не характеризуется
чебышевским альтернансом, состоящим из $n+1$ точек отрезка
$[-1,1]$.

Впервые это было показано в \cite{[DAN-KOND2]} на примере НД
порядка 2 вида
\begin{equation}
\frac{2x+\lambda}{x^2+\lambda x+1}, \qquad
\lambda\in[1,\lambda^*], \quad \lambda^*=1.62\dots,
  \label{(NEEDIN)}
\end{equation}
каждая из которых наименее уклоняется от функции $f(x)=x+1$ на
отрезке $[-1,1]$ среди всех НД порядка не выше $2$. Величина
наименьшего уклонения равна $1$. Только в случае $\lambda=\lambda^*$
имеется альтернанс из трех точек отрезка $[-1,1]$; для остальных
значений $\lambda$ нет даже двухточечного альтернанса.

К этому примеру можно добавить следующее. Для функции $x+1/2$ НД
наилучшего приближения уже единственна, она совпадает с НД
(\ref{(NEEDIN)}), где $\lambda=1$, причем существует альтернанс из
трех точек: $\{-1,0,1\}$. Таким образом, при приближении
посредством НД функций, отличающихся друг от друга лишь на
константу, могут проявляться совершенно разные аппроксимативные
свойства.

Позже пример неединственности, равно как и необязательности
альтернанса для наилучшего приближения, был построен и для
произвольного $n\ge 2$ (см. \cite{[Kom2012-2],[KOM2]}).
Конструкция весьма сложна и мы остановимся лишь на ее идее.

Рассмотрим таблицу интерполяции $T=\{(x_k,y_k)\}_{k=1}^{n}$ с
вещественными узлами:
\begin{equation}
\label{8.2-X}
0=x_1<x_2<\ldots<x_n,\qquad y_k=\frac{n}{x_k} \quad
(k=\overline{2,n}),
    \quad y_1=-\frac{y_2+\ldots+y_n}{n-1}.
\end{equation}
Несложно проверить, что  задача интерполяции таблицы $T$ имеет
бесконечно много решений вида
$$
    \begin{array}{l}
        \rho_{n,s}(x)=(Q_{n,s}(x))'_x/Q_{n,s}(x), \\
        Q_{n,s}(x)=x^n+s\cdot \left(x^{n-1}+\sum_{j=2}^n(-1)^{j-1}\cdot
        \frac{\sigma_{j-1}}{j}x^{n-j}\right),\quad s>0,
    \end{array}
$$
где  $\sigma_m=\sigma_m(x_2,\ldots,x_n)$ --- элементарные
симметрические многочлены
$$
    \sigma_0:= 1, \qquad \sigma_k(z_1,\ldots,z_q):=\sum_{1\le j_1<j_2<\ldots<j_k\le q}
    z_{j_1}z_{j_2}\dots z_{j_k}, \qquad k=\overline{1,q}.
$$
При этом существует отрезок $I:=[s_1,s_2]$ значений параметра $s$
такой, что $\rho_{n,s}$ не имеет полюсов на $[x_1,x_n]$. Отсюда, в
частности следует, что $\rho_{n,s}$, $s\in I$, образуют компактное
семейство на $[x_1,x_n]$. Вследствие этого, выбрав
$\varepsilon>0$, можно определить непрерывную функцию $f$ так, что
\begin{equation}\label{111-R*}
 f(x_1)=y_1+\varepsilon,\quad f(x_k)=y_k+(-1)^k
\varepsilon, \qquad k=\overline{2,n},\qquad
|f(x)-\rho_{n,s}(x)|<\varepsilon
\end{equation}
при $x\ne x_k$, $s\in I_1\subset I$. Построение завершается
доказательством существования числа $\varepsilon^*$ и таблицы
$(x_k^{*},y_k^{*})$ вида (\ref{8.2-X}) таких, что любая НД, для
которой неравенство в (\ref{111-R*}) выполнено при
$\varepsilon=\varepsilon^*$ всюду, включая узлы $x=x_k^{*}$, имеет
полюс на отрезке $[0,x_n]$.

Отметим, что величина наилучшего приближения в этом примере равна
$\varepsilon^*$ и достигается во всех $x_k^*$. При этом разность
$f(x)-\rho_{n,s}(x)$ имеет на отрезке $[0,x_n^*]$ альтернанс ровно
из $n-1$ точек. В \cite{[KOM2]} показано также, что функцию $f$
при построении примера можно выбрать так, что при одном значении
параметра $s\in I_1$ разность будет альтернировать в $n+1$ точках
отрезка, а при других его значениях, по-прежнему, --- в $n-1$
точках.

\subsection{Наилучшее приближение констант и альтернанс}
\label{paragraph8.3} Из предыдущего параграфа видно, что для НД не
существует точного аналога теоремы
П.~Л.~Че\-бы\-ше\-ва об альтернансе. Однако для функций некоторых
классов такие аналоги все же имеют место. Первые результаты в этом
направлении получены для {\it вещественных постоянных функций}
$f(x)=c$ (см.
\cite{[DAN-KOND2],[Kom2011],[Kom2013],[KOM3],[Kom2015-MZ]}).
Отметим, что задачу о наилучшем приближении констант можно считать
одним из аналогов задачи П.\,Л.~Чебышева об унитарном многочлене,
наименее уклоняющемся от константы.

В работе \cite{[DAN-KOND2]} при достаточно больших $n>n_0(c)$
доказана {\it необходимость} существования $(n+1)$-точечного
альтернанса для наилучшего приближения констант
$c\in(0,\ln\sqrt{2})$ посредством НД порядка $n$. Там же доказана
{\it единственность} вещественнозначной НД $\rho^{*}_n(c;x)$
наилучшего приближения. В работе \cite{[KOM3]} доказана {\it
достаточность} указанного альтернанса, но при весьма жестком
ограничении: $|c|\le c_n$, где последовательность $\{c_n\}$
достаточно быстро убывает к нулю. Эти предварительные результаты
объединяет и усиливает следующее утверждение из \cite{[Kom2015-MZ]}:\smallskip

{\it Для любой вещественной константы $c\ne 0$ дробь
$\rho_n^*(c;\cdot)$ при каждом $n>8|c|+1$ единственна и имеет
порядок $n$. Существование чебышевского альтернанса из $n+1$ точек
на отрезке $[-1,1]$ для разности $\rho_n-c$ необходимо и
достаточно для того, чтобы $\rho_n(x)\equiv \rho_n^*(c;x)$. Для
наименьших уклонений ${\mathcal R}_n(c)$ имеем двустороннюю
оценку:
\begin{equation}\label{8.4.delta}
    \frac{|c|^{n+1}}{2^{n-1}n!}
    \cdot\frac{e^{-|c|}}{\ch c\cdot W_n(|c|)}\le
    {\mathcal R}_n(c)
    \le\frac{|c|^{n+1}}{2^{n-1}n!}\cdot\frac{e^{2|c|}}{V_n(|c|)}
\end{equation}
с определенными $V_n(|c|)>0$, $W_n(|c|)>0$, стремящимися к $1$ при
$n\to\infty$.}\smallskip

Опустить или существенно ослабить ограничение на $n$ в утверждении
теоремы нельзя в следующем смысле. Из \cite{[M-F]} следует, что
для любой НД $\rho_n$ найдется точка $x_0\in [-1,1]$, в которой
$|\rho_n(x_0)|<A\, n$, где $A$ --- константа. Значит, константы
$c>A\, n$ заведомо не могут хорошо приближаться на $[-1,1]$
дробями $\rho_n$ (наилучшее приближение имеет порядок самой
константы). При этом, в \cite{[Kom2015-IzvRAN],[Kom2015-MZ]}
отмечено, что из (\ref{KOS-N}) и известной оценки величины
$E_{n-1}(F,[-1,1])$ для функций $F$ со знакопостоянной на $[-1,1]$
производной $F^{(n)}$ также получается двусторонняя оценка
наименьших уклонений НД от константы:
\[\frac{|c|^{n+1}}{2^{n-1}n!}\cdot\frac{e^{-2|c|}}{1+2|c|e^{|c|}}\le
\mathcal{R}_n(c)\le \frac{|c|^{n+1}}{2^{n-1}n!}
\cdot2(1+|c|)e^{2|c|}.\] В этой оценке нет ограничений на
соотношение между $c$ и $n$, но постоянные множители при
$|c|^{n+1}2^{1-n}/n!$ менее точны, чем в (\ref{8.4.delta}).

Правдоподобна гипотеза о том, что, как и в случае констант, аналог
теоремы П.~Л.~Чебышева верен для любой непрерывной вещественнозначной
функции $f$, {\it если порядок $n\ge n_0(f)$ достаточно велик}
\cite{[D-D.Abstr2010]}. Например, в \cite{[Kom2012-2]} построен
пример появления единственности и альтернанса с ростом $n$, а
именно, показано, что дробь $\rho_3^*(x^3+1;\cdot)$ единственна и
характеризуется альтернансом из 4-х точек, тогда как при $n=2$
наилучшее приближение функции $x^3+1$ неединственно и альтернанс
необязателен.

В заключение этого раздела заметим, что оценки скорости
приближения констант находят применение в оценках наилучшего
приближения более общих функций. Например, в
\cite{[Kom2017-IzvRAN]} с помощью (\ref{8.4.delta}) показано, что
$\mathcal{R}_{2n}(Ax)\asymp A^{n+1}8^{-n}(n!)^{-1}$ с любым $A>0$
при $n>2A+1$.

\subsection{Алгоритм построения НД наилучшего приближения
констант} \label{paragraph8.4}

В теории полиномиальных аппроксимаций хорошо известен алгоритм
Ремеза, с помощью которого можно построить многочлен наилучшего
приближения и соответствующий альтернанс из $n+2$ точек. Для
$\rho_n-c$ в \cite{[DAN-KOND2]} предложен численный алгоритм
построения альтернанса из $n+1$ точек отрезка $[-1,1]$ в случае
констант $c\in(0,\ln\sqrt{2})$, $n>n_0(c)$. Мы пока не имеем
строгого обоснования сходимости, но в численных экспериментах
алгоритм работает надежно. Таким образом, с учетом предыдущей
теоремы, алгоритм приводит к численному построению НД
$\rho_n^*(c;\cdot)$ наилучшего приближения. Приведем схему этого
алгоритма \cite{[DAN-KOND2]}.

Сначала берутся произвольно узлы $-1<{\xi}_1<\ldots<{\xi}_n<1$ и
для $f(x)=c$ строится интерполяционная НД $\rho_n(x)$ (один из
способов построения см. в (\ref{(777.4)}) и (\ref{(777.4+)})).
Вычисляются $\sup$-нормы $N_k$ разности $\rho_n(x)-c$ на отрезках
$[{\xi}_{k-1},{\xi}_k]$, $k=1,\ldots,n+1$, где ${\xi}_0=-1$,
${\xi}_{n+1}=1$. Если все $N_k$ равны, то искомая НД построена. В
противном случае выбирается какой-либо номер $m$, для которого
значение $N_m$ минимально. Если $2\le m \le n$, то соответствующие
узлы ${\xi}_{m-1}$ и ${\xi}_m$ слегка <<раздвигаются>>, т.е.
заменяются на ${{\xi}_{m-1}-\varepsilon},\;{{\xi}_m+\varepsilon}$,
$\varepsilon>0$. Если $m=1$, то узел ${\xi}_1$ <<сдвигается>>
вправо, если $m=n+1$, то узел ${\xi}_n$ <<сдвигается>> влево. И
все повторяется с новым набором узлов. Показано, что при
достаточно малых $\varepsilon>0$ имеем убывание норм разностей
интерполяционных НД и $c$ на $[-1,1]$. Параметр $\varepsilon$ в
этом процессе, естественно, уменьшается.

За начальный набор узлов интерполяции рекомендуется брать набор
узлов Чебышева, т.к. норма разности соответствующей
интерполяционной НД и $c$ близка по порядку к наименьшему
уклонению \cite{[Kom2015-MZ]}.

\subsection{Критерии наилучшего приближения вещественнозначных
функций на отрезке} \label{paragraph8.5} Я.~В.~Новак
\cite{[NOV-diss],[NOV2010]} доказал следующий критерий наилучшего
приближения вещественнозначными НД вещественных функций $f\in
C([-1,1])$, вполне аналогичный критерию Колмогорова
\cite{[Kolmogorov48]} для полиномов:\smallskip

{\it Имеем ${\rho}_n(x)\equiv \rho_n^*(f;x)$, если и только если для
произвольной НД $\widetilde{\rho}_n$, не имеющей полюсов на
$[-1,1]$, выполняется неравенство
\begin{equation}\label{NOVAK-100}
\min_{x\in E}
({\rho}_n(x)-\widetilde{\rho}_n(x))({{\rho}_n(x)-f(x)})\le 0,
\end{equation}
где
$E=\{x\in[-1,1]:|f(x)-{\rho}_n(x)|=\|f-{\rho}_n\|_{C([-1,1])}\}$.}
\smallskip

Для компакта $K$ общего вида Я.~В.~Новаком доказана достаточность
условия Колмогорова типа (\ref{NOVAK-100}), где надо второй
множитель заменить на комплексно сопряженный и заменить полученное
произведение его вещественной частью.

Следующий критерий  наилучшего приближения в терминах альтернанса (далее --- {\it
теорема об альтернансе} или {\it критерий}) в настоящее время
является наиболее общим результатом в этом направлении (см.
\cite{[Kom2014-1],[Kom2015-IzvRAN],[Kom-2016],[Kom2017-IzvRAN]}).\smallskip

{\it Пусть полюсы вещественнозначной НД $\rho_n$ порядка $\le n$
расположены вне круга $|z|\le 1$. Тогда $\rho_n(x)\equiv
\rho_n^*(f;x)$ в том и только том случае, когда для разности
$f-\rho_n$ на $[-1,1]$ имеется альтернанс из $n+1$ точек. При этом
$\rho^{*}_n(f;x)$ является единственной НД наилучшего приближения.
Условие на расположение полюсов ослабить нельзя.}\smallskip

Невозможность ослабления условия на полюсы означает, что если хотя
бы одна пара комплексных сопряженных полюсов $\rho_n$ попадает в
круг $|z|\le 1$, то, вообще говоря, при выполнении других условий
теоремы $\rho_n$ может не быть НД наилучшего приближения, или НД
наилучшего приближения неединственна. Так, в первом примере \S
\ref{paragraph8.2} возникает неединственность, хотя при
$\lambda=\lambda^*$ на отрезке $[-1,1]$ имеется альтернанс из трех
точек: условие критерия не выполнено, поскольку полюсы НД
(\ref{(NEEDIN)}) лежат на единичной окружности.

В построенном ниже в \S \ref{paragraph8.7} примере показана
недостаточность даже $2n-1$ точек альтернанса для наилучшего
приближения. Там условие критерия также не выполнено: все полюсы
НД лежат {\it внутри} единичного круга. Другие примеры, связанные
с невозможностью ослабления условия на полюсы, см. в
\cite{[Kom2015-IzvRAN]}.

Ранее аналог теоремы об альтернансе для случая $n$ {\it
вещественных различных} полюсов за исключением утверждения о
единственности доказал Я.~В.~Новак \cite{[NOV-diss]}. Для случая
$n$ {\it различных} полюсов достаточность
альтернанса была доказана в \cite{[Kom2013]}, а необходимость
альтернанса и единственность наилучшего приближения --- в
\cite{[Kom2015-PMA]}.

Теорему об альтернансе можно переформулировать. Пусть $X$ ---
произвольный компакт на $\mathbb C$, расположенный симметрично
относительно $\mathbb R$ и не имеющий общих точек с кругом $|z|\le
1$. Через $ \rho^{**}_n(X;f;x)$ обозначим НД наилучшего
приближения функции $f$ на $[-1,1]$ среди всех вещественнозначных
НД порядка $n$, полюсы которых лежат на $X$ (такая НД, очевидно,
существует и имеет порядок, равный $n$). Приходим к следующему
результату:\smallskip

{\it Имеем $\rho^{**}_n(X;f;x)\equiv \rho_n^*(f;x)$ в том и только том
случае, когда для разности $f(x)-\rho^{**}_n(X;f;x)$ на $[-1,1]$
имеется альтернанс из $n+1$ точек.}\smallskip

Здесь условие на компакт ослабить нельзя в следующем смысле. Если
$X$ имеет непустое пересечение с кругом $|z|\le 1$, то альтернанс,
вообще говоря, не гарантирует, что $\rho^{**}_n(X;f;x)\equiv
\rho_n^*(f;x)$ \cite{[Kom2015-IzvRAN]}.

Теорему об альтернансе дополняет аналог теоремы Валле-Пуссена
\cite[Гл. 2, \S32]{[Akhiezer]}, полученный в
\cite{[Kom2015-IzvRAN]} (см. также \cite{[Kom2014-1]}):\smallskip

{\it Если полюсы НД $\rho_n$ порядка $\le n$ лежат вне круга
$|z|\le 1$ и найдутся точки $-1\le t_1<\ldots<t_{n+1}\le 1$ такие,
что
$$
f(t_j)-\rho_n(t_j)=\pm(-1)^j a_j,\qquad a_j>0,\quad
j=\overline{1,n+1},
$$
то ${\mathcal R}_n(f)\ge \min \{a_1,\dots,a_{n+1}\}$.}\smallskip

Отметим, что для $\rho_n$ порядка, равного $n$, подобное
утверждение фактически получено уже в \cite{[Kom2013]} (см. также
\cite{[NOV-diss]}).

\subsection{ Чебышевские системы функций, связанные с НД. Идея
доказательства теоремы об альтернансе} \label{paragraph8.6}
Напомним, что система $\{f_1,\dots,f_n\}$ функций, непрерывных на
множестве $K\subset\mathbb{C}$, содержащем не менее $n+1$ точек,
удовлетворяет {\it условию Хаара}, если для всякого обобщенного
полинома вида
$$
F_n(\alpha_1,\ldots,\alpha_n;z):=\alpha_1 f_1(z)+\ldots+\alpha_n
f_n(z)\not\equiv 0
$$
с коэффициентами $\alpha_j\in\mathbb{C}$ число его различных нулей
на $K$ не превосходит $n-1$ (см., например, \cite[Гл.~1,
\S\S43-48]{[Akhiezer]}). Условие Хаара, очевидно, эквивалентно
тому, что $\det(f_j(\xi_l))_{l,j=1}^n\ne 0$ для любых различных
$n$ точек $\xi_l\in K$ и, как доказал А.\,Хаар, равносильно {\it
единственности} полинома $\alpha_1^* f_1+\ldots+\alpha_n^* f_n$
наилучшего приближения для каждой непрерывной на $K$ функции $f$.
В частности, если функции $f_j$ вещественны и удовлетворяют
условию Хаара на отрезке $[-1,1]$, то систему $\{f_1,\dots,f_n\}$
называют {\it системой Чебышева} на $[-1,1]$. На системы функций
Чебышева дословно переносится теорема П.~Л.~Чебышева об
альтернансе и некоторые другие результаты теории приближений
алгебраическими полиномами.

Будем говорить, что система $\{f_1,\dots,f_n\}$ аналитических
функций удовлетворяет {\it усиленному условию Хаара} в (замкнутой
или открытой) области $K\subset\mathbb{C}$, если для всякого
обобщенного полинома $F_n(\alpha_1,\ldots,\alpha_n;z)\not\equiv 0$
сумма кратностей его нулей на $K$ не превосходит $n-1$.

Пусть точки $z_1,\dots,z_k$ лежат вне круга $|z|\le 1$ и
расположены симметрично относительно вещественной оси ${\mathbb
R}$. Введем систему $\chi$ из $n$ функций вида
\begin{equation}\label{SystF}
    f_{j,s}(x)=(x-z_j)^{-s}, \qquad j=\overline{1,k}, \quad
    s=\overline{2,m_j+1}; \quad \sum\nolimits_{j=1}^k m_j=n,
\end{equation}
где каждому полюсу $z_j$ сопоставлено число
$m_j=m(z_j)\in\mathbb{N}$, причем $m(z_j)=m(z_p)$, если
$z_j=\overline{z_p}$ и $z_j\notin{\mathbb R}$.

Одним из важных вспомогательных инструментов доказательства
теоремы об альтернансе и теоремы типа Валле-Пуссена (см. \S
\ref{paragraph8.5}) является следующая теорема
\cite{[Kom2015-PMA],[Kom2015-IzvRAN]}, представляющая и
самостоятельный интерес.\smallskip

{\it Система $\chi$ на отрезке $[-1,1]$  удовлетворяет усиленному
условию Хаара.}\smallskip

Ограничимся схемой доказательства ослабленного варианта теоремы:
система $\chi$ удовлетворяет обычному условию Хаара в случае,
когда все $m_j=1$. Как говорилось, в этом случае условие Хаара
равносильно условию
$$
{\mathcal
D}(\{\xi_l\};\{z_j\}):=\det((\xi_l-z_j)^{-2})_{l,j=1}^n\ne 0
 $$
по всем наборам из $n$ различных точек $\xi_l$ отрезка $[-1,1]$. В
доказательстве этого неравенства ключевым является тождество
Борхарта~(см. \cite[\S1.3]{[Minc]}):
\begin{equation}\label{D=W*T}
    {\mathcal
D}(\{\xi_l\};\{z_j\})=W\cdot T,
\end{equation}
где $W$ и $T$ --- детерминант и перманент\footnote{Перманент
квадратной матрицы --- это матричная функция, вычисляемая по
правилу разложения определителя, но с тем отличием, что перед
каждым произведением элементов независимо от четности
соответствующей перестановки ставится плюс.} матрицы
$((\xi_l-z_j)^{-1})_{l,j=1}^n$ соответственно. Определитель $W$
раскладывается по известной (см. \cite[Гл. 1, \S14]{[Akhiezer]})
формуле
$$
W=\det((\xi_l-z_j)^{-1})_{l,j=1}^n=\frac{\prod_{1\le l< j\le
n}(\xi_j-\xi_l)(z_l-z_j)}{\prod_{l,j=1}^n(\xi_l-z_j)}.
$$
В рассматриваемом случае точки $\xi_l$, равно как и полюсы $z_j$,
попарно различны, поэтому $W\ne 0$. Несколько сложнее доказывается
(см. \cite{[Kom2013],[Kom2015-IzvRAN]}), что и перманент $T\ne 0$.
В результате приходим\footnote{Если все $z_j$ вещественны и лежат
за пределами отрезка $[-1,1]$, естественно говорить, что система
система $\chi$ удовлетворяет на отрезке {\it усиленному условию
Чебышева}. Это утверждение для системы функций
$\{(x-z_j)^{-2}\}_{j=1}^n$ с $n$ {\it попарно различными} полюсами
$z_j\in {\mathbb R}$ при условии справедливости тождества Борхарта
(\ref{D=W*T}) впервые доказано в \cite{[NOV-diss]} (см. также
\cite{[Kom2013],[Kom2015-PMA]}).} к искомому неравенству
${\mathcal D}(\{\xi_l\};\{z_j\})\ne 0$. \smallskip

Применение теоремы о системе Хаара к доказательству теоремы об
альтернансе осуществляется посредством следующей леммы и ее
следствия \cite{[Kom2011],[Kom2015-IzvRAN]}:\smallskip

{\it Пусть $\rho_n(x)=(\ln Q(x))'$, где $Q$
--- полином степени $n$ вида
$$
Q(x)=(x-z_1)^{m_1}\ldots(x-z_k)^{m_k} \qquad
(m_1+\ldots+m_k=n),
$$
$z_1,\dots,z_k$ --- попарно различные точки на $\mathbb{C}$, а
$m_j$ --- натуральные числа. Тогда для любой НД
$\tilde\rho_n(x)=(\ln P(x))'$ порядка не выше $n$ существуют $n$
однозначно определенных чисел $\alpha_{j,s}$ таких, что}
$$
    \tilde\rho_n(x)-\rho_n(x)\equiv\frac{Q(x)}{P(x)}\cdot\sum_{j=1}^k\sum_{s=2}^{m_j+1}
    \frac{\alpha_{j,s}}{(x-z_j)^s}.
$$

Отсюда ввиду указанного выше свойства системы $\chi$ получается утверждение:
\smallskip

{\it Если все полюсы вещественнозначной НД $\rho_n$ порядка $n$
лежат вне круга ${|z|\le 1}$, а $\tilde\rho_n\not\equiv \rho_n$
--- другая вещественнозначная НД порядка не выше $n$, то сумма
кратностей нулей разности $\tilde\rho_n-\rho_n$ на $[-1,1]$ не
превосходит $n-1$. В частности, $\tilde\rho_n-\rho_n$ имеет на
$[-1,1]$ не более $n-1$ различных нулей.}
\smallskip

Иначе это следствие можно переформулировать как теорему о
единственности решения задачи обычной интерполяции {\it
вещественных} таблиц
\cite{[Kom2015-IzvRAN],[Kom2017-IzvRAN]}:\smallskip

{\it Пусть полюсы вещественнозначной НД $\rho_n$ порядка $\le n$
лежат вне круга ${|z|\le 1}$. Если $\rho_n$ является решением полной
задачи $(\ref{(7.5)})$ с узлами $\xi_j\in [-1,1]$, то $\rho_n$ ---
единственная вещественная интерполяционная дробь.}\smallskip

Теперь легко получается достаточность альтернанса в случае, когда
порядок НД $\rho_n$ равен $n$. Действительно, пусть для разности
$f-\rho_n$ существует альтернанс из $n+1$ точек $-1\le
t_1<\ldots<t_{n+1}\le 1$. Допустим существование НД
$\tilde\rho_n$, для которой
$\|f-\rho_n\|_{C([-1,1])}>\|\tilde\rho_n-L\|_{C([-1,1])}$. Тогда
при $j=\overline{1,n+1}$ имеем
$$
{\rm sgn}[\tilde\rho_n(t_j)-\rho_n(t_j)]={\rm
sgn}[(\tilde\rho_n(t_j)-f(t_j))-(\rho_n(t_j)-f(t_j))]=\pm (-1)^j,
$$
так что функция $\tilde\rho_n-\rho_n$ имеет на интервале $(-1,1)$
не менее $n$ различных нулей. Но тогда по следствию
$\tilde\rho_n\equiv \rho_n$. Противоречие с допущением.

С помощью некоторых дополнительных соображений (см.
\cite{[Kom2017-IzvRAN]}) получается достаточность альтернанса и в
случае НД $\rho_n$ порядка $\le n$.

\subsection{Минимальное число точек альтернанса,
гарантирующее наилучшее приближение} \label{paragraph8.7} Вопрос о
минимальном числе $N$ точек альтернанса, гарантирующем наилучшее
приближение посредством НД порядка $\le n$ независимо от свойств
приближаемой непрерывной функции, рассматривался в
\cite{[Kom2012-2]}. Оказывается, что $N=2n$. Действительно,
достаточность такого числа точек альтернанса очевидна, а
недостаточность $2n-1$ точек подтверждается примером
\cite{[Kom2012-2],[Kom2015-IzvRAN]}, который мы изложим для случая
четных $n$.

Введем положительные на ${\mathbb R}$ многочлены
$$
P(x)=\varepsilon+\prod_{k=1}^{m}(x+2^{-k})^2,\qquad
Q(x)=\varepsilon+\prod_{k=1}^{m}(x-2^{-k})^2,\qquad \varepsilon>0,
$$
степени $n=2m$. Доказано \cite{[Kom2012-2]}, что при достаточно
малом $\varepsilon>0$ многочлен $P'Q-Q'P$ имеет ровно $2n-2$
простых нулей $\xi_k\in (-1,1)$. Отсюда следует, что разность НД
$\rho_n=Q'/Q$ и $\widetilde{\rho}_n=P'/P$ в точках $\xi_k$ имеет
$2n-2$ перемен знаков. Легко видеть, что любая непрерывная
вещественнозначная функция $f$, график которой проходит через
точки $(\xi_k,\rho_n(\xi_k))$ и которая удовлетворяет условию
$$
{\rm sgn}~(f(x)-\rho_n(x))={\rm sgn}~
(f(x)-\widetilde{\rho}_n(x))={\rm sgn}~
(\widetilde{\rho}_n(x)-\rho_n(x))
$$
($2n-2$ перемен знаков), приближается (на любом отрезке) дробью
$\widetilde{\rho}_n$ лучше, чем дробью ${\rho}_n$. Для того, чтобы
разность $f(x)-\rho_n(x)$ имела $2n-1$ точек альтернанса, остается
дополнительно потребовать от $f$, чтобы на всех $2n-1$ промежутках
отрезка $[-1,1]$, где знаки разности $f(x)-\rho_n(x)$ постоянны,
$\sup$-нормы этой разности были одинаковыми.

\section{ Аппроксимация посредством НД в $L_p(\mathbb{R})$. Ряды НД}
\label{paragraph3}
\subsection{Класс аппроксимируемых функций}
\label{paragraph3.1} Пусть $p\in (1,\infty)$. Через $S_p$
обозначим класс всех комплекснозначных функций $f\in
L_p=L_p(\mathbb{R})$, приближаемых сколь угодно точно посредством
НД в метрике $L_p$. Полное описание класса $S_p$ получено
В.~Ю.~Протасовым в работе \cite{[PRO]}. Им, в частности, показано,
что любая функция $f\in S_p$ является аналитической на
$\mathbb{R}$, продолжается до мероморфной на комплексной плоскости
${\mathbb C}$ функции и представляется в виде ряда
$f(x)=\rho_{\infty}(x)=\sum_{k=1}^{\infty}\;(x-z_k)^{-1}$,
сходящегося в $L_p$. При этом показатель сходимости
последовательности $\{z_k\}$ удовлетворяет неравенству
\begin{equation}
\tau(\{z_k\}) \le 1-1/p,\qquad \tau(\{z_k\}):= \inf
\left\{\gamma>0, \;\;\sum\nolimits_k
|z_k|^{-\gamma}<\infty\right\}.
  \label{(*30)}
\end{equation}
Таким образом, класс $S_p$ состоит из тех и только тех функций
$f$, которые представляются в виде сходящихся к ним в $L_p$ рядов
НД. Этот результат инициировал изучение рядов НД
\cite{[KAY2011],[KAY2012+],[KAY2012],[DanDod2013]}.

Отметим, что в случае $p=\infty$ ситуация значительно изменяется.
Как показали П.~А.~Бородин и О.~Н.~Ко\-су\-хин
\cite{[BOR-KOS],[BOR2009]}, в равномерной метрике {\it каждая}
непрерывная на $\mathbb{R}$ функция $f$ с нулевым значением на
бесконечности с любой точностью приближается посредством НД.

Предположим, что ряд $\sum\nolimits_{k}(x-\xi_k)^{-1}$ сходится в
$L_p$, $p\in (1,\infty)$, к некоторой функции $\rho\in L_p$ (т.е.
$\rho\in S_p$). Запишем это в виде
\begin{equation}
\rho(x)=\rho_{\infty}(x)=\sum^{(p)}\nolimits_{k}(x-\xi_k)^{-1},\qquad
x\in \mathbb{R}.
   \label{(*31)}
\end{equation}
Ряд (\ref{(*31)}) сходится к $\rho$ в $L_p$ безусловно
\cite{[PRO]}. Через $z_k=x_k+iy_k$, $k=1,2,\ldots$, $y_k>0$, будем
обозначать занумерованные в каком-либо порядке полюсы суммы
(\ref{(*31)}), лежащие в $\mathbb{C}^{+}$. Если число $m\ge 0$
таких полюсов конечно, то для единообразия записей будем считать
$z_{m+k}=\infty$ и $(z-z_{m+k})^{-1}\equiv 0$, $k \in \mathbb{N}$.
То же будем предполагать относительно полюсов, лежащих в
$\mathbb{C}^{-}$; будем обозначать их через $\widetilde
z_k=\widetilde x_k+i \widetilde y_k$. Введем частичные суммы
$\rho_n(z)+\widetilde\rho_n(z)$, где
$$
\rho_n(z) =\sum\nolimits_{k=1}^{n}\frac{1}{z-z_k},\quad
\widetilde\rho_n(z)=\sum\nolimits_{k=1}^{n}\frac{1}{z-\widetilde
z_k}.
$$
и положим $\mu_n={\rm Im}\,\rho_n$, $\nu_n={\rm Re}\,\rho_n$,
$\widetilde\mu_n={\rm Im}\,\widetilde\rho_n$,
$\widetilde\nu_n={\rm Re}\,\widetilde\rho_n$ (см. (\ref{(*32)})).
Доказательство основных результатов в \cite{[PRO]} базируется на
преобразовании Гильберта, которое (см., например, \cite{[GARNET]}) определяется (почти всюду) как
\begin{equation}
H(f)(x) =\frac{1}{\pi}\lim_{\varepsilon \to 0}
\int\nolimits_{|t-x|\ge \varepsilon} \frac{f(t)}{t-x}dt,\qquad
f\in L_p, \qquad x\in \mathbb{R},
  \label{(*33)}
\end{equation}
и на формулах Сохоцкого-При\-ва\-ло\-ва $H(\rho_n)=-i\rho_n$ и
$H(\widetilde\rho_n)=i\widetilde\rho_n$, т.е.
\begin{equation}
H(\nu_n)=\mu_n,\quad H(\mu_n)=-\nu_n,\quad
H(\widetilde\nu_n)=-\widetilde\mu_n,\quad
H(\widetilde\mu_n)=\widetilde\nu_n.
    \label{(*34)}
\end{equation}
Хорошо известно, что оператор $H:\;L_p\rightarrow L_p$ при $p\in
(1,\infty)$ ограничен, поэтому из (\ref{(*34)}), в частности,
получаются слабые эквивалентности $\|\nu_n \|_{p}\asymp \|\mu_n
\|_{p}\asymp \|\rho_n \|_{p}$ (не зависящие от $n$). Здесь и далее
в этом параграфе применяется обозначение
$\|\cdot\|_p:=\|\cdot\|_{L_{p}}$, $1<p\le\infty$.

Пусть выполнено (\ref{(*31)}) при $1<p<\infty$. Тогда с учетом
(\ref{(*34)}) получим (см. \cite{[PRO]})
\begin{equation}
\label{(*35)}
\begin{array}{l}
 \|\mu_n\|_p \le \|\mu_n-\widetilde\mu_n\|_p
=\|H(\nu_n+\widetilde\nu_n)\|_p\le\\
\qquad\qquad\le h_p\|\nu_n+\widetilde\nu_n\|_p \le h_p\,
\|\rho_n+\widetilde\rho_n\|_p \to h_p\, \|\rho\|_p,\quad n\to
\infty,
\end{array}
\end{equation}
где $h_p$ --- норма оператора Гильберта. Из теоремы Б.~Леви
следует, что неубывающая последовательность $\mu_n$ сходится в
$L_p$ к функции $\mu$ и $\|\mu\|_p \le h_p \,\|\rho\|_p$. Отсюда и
из $\|\nu_n \|_{p}\asymp \|\mu_n \|_{p}$ вытекает, что
последовательность $\{\nu_n\}$ является фундаментальной в $L_p$ и,
значит, сходится к некоторой функции $\nu$, причем $\|\nu\|_p\le
h_p^2 \,\|\rho\|_p$. Следовательно, последовательность $\rho_n$
сходится в $L_p$ к функции $\sigma=\mu+i\nu$ и $\|\sigma\|_p\le
h_p(1+h_p)\|\rho\|_p$. Аналогичны рассуждения и оценки для функций
со знаком~$\widetilde{\phantom{z}}$.

Здесь важно, что $p<\infty$, поскольку оценка компонент $\sigma$ и
$\widetilde\sigma$ через их сумму $\rho=\sigma+\widetilde\sigma$
при $p=\infty$ невозможна. Можно утверждать \cite{[DAN2006]} лишь то, что
\begin{equation}
\|\mu_n\|_{\infty}\le A\cdot \ln
n\cdot\|\rho_n+\widetilde\rho_n\|_{\infty}.
    \label{(*36)}
\end{equation}
Эта оценка точна по порядку.

Таким образом, для сходимости (\ref{(*31)}) необходима и
достаточна одновременная сходимость в $L_p$ сумм $\rho_n$ и
$\widetilde\rho_n$, или $\nu_n$ и $\widetilde\nu_n$, или $\mu_n$ и
$\widetilde\mu_n$ (то есть одновременная конечность $L_p$-норм
сумм $\mu(x):=\lim_{n\to\infty}\mu_n(x)$,
$\widetilde\mu(x):=\lim_{n\to\infty}\widetilde\mu_n(x)$
знакопостоянных рядов).

Итак, задача аппроксимации в $L_p=L_p(\mathbb{R})$ при конечных
$p$ фактически сводится к исследованию сходимости рядов НД в
$L_p$, а последнее значительно облегчается тем, что достаточно
изучать ряды, {\it все полюсы которых лежат в $\mathbb{C}^{+}$}.
Кроме того, в этом случае сходимость (\ref{(*31)}) равносильна
конечности $L_p$-нормы $\|\mu\|_p$ знакоположительного ряда
$\mu(x)=\lim_{n\to\infty}\mu_n(x)$. В связи с этим в работе
\cite{[PRO]} сформулирована задача: {\it каков критерий сходимости
(\ref{(*31)}) в терминах последовательности полюсов
$\{z_k\}\subset\mathbb{C}^{+}$, или, что то же самое, каков
критерий для $\|\mu\|_p<\infty$}. (При $p=\infty$ задача о
сходимости представляется значительно более сложной из-за того,
что вместо (\ref{(*35)}) можно утверждать лишь (\ref{(*36)}).)
Приведем несколько результатов.

\subsection{Некоторые критерии сходимости рядов НД}
\label{paragraph3.2} Пусть $\rho_n$ --- НД с полюсами
$z_k\in{\mathbb C^+}$. При $1<p<3$ имеем следующий критерий
сходимости \cite{[DAN2010]}:
 $$
\|\mu\|_p<\infty  \quad\Leftrightarrow \quad {\rm
Im}\,\left(e^{-i\frac{\pi(p-2)}{2}}\sum\nolimits_{k=1}^{n}
\rho_n^{p-1}(\overline{z_k})\right)\le A<\infty \quad (\forall n),
 $$
где $A=A(p,\{z_k\})$ не зависит от $n$, а суммируются значения
однозначной в ${\mathbb C}^{-}$ аналитической ветви
$$
\rho_n^{p-1}(z)=|\rho_n(z)|^{p-1}e^{i\varphi(p-1)},\quad
\varphi=\arg\rho_n(z) \in (0,\pi),\quad z\in{\mathbb C}^{-}.
$$
При этом справедлива двусторонняя оценка
\begin{equation} \|\rho_n\|_p^p
\cos\frac{\pi(p-2)}{2}\le 2\pi\;{\rm
Im}\,\left(e^{-i\frac{\pi(p-2)}{2}}\sum\nolimits_{k=1}^{n}
\rho_n^{p-1}(\overline{z_k})\right)\le \|\rho_n\|_p^p.
    \label{(*38)}
\end{equation}
Легко проверить, что при $p\in (1,2)$ мнимые части слагаемых в
(\ref{(*38)}) слабо эквивалентны их модулям, откуда получается
весьма простой по форме критерий: {\it для $\|\mu\|_p<\infty$,
$1<p<2$, необходимо и достаточно выполнение} (\ref{(*777)}).

В общем случае с помощью (\ref{(*8)}) получается следующий
критерий \cite{[DAN2010]}.\smallskip

{\it Условие $\|\mu\|_p<\infty$ $($$\rho_n\to \rho$ в $L_p$$)$
выполняется тогда и только тогда, когда найдется конечная величина
$A=A(p,\{z_k\})$ такая, что для любых $n\in{\mathbb N}$ и $g\in
H_q$ выполняется неравенство}
$$
\left|\sum\nolimits_{k=1}^{n} g(z_k)\right|\le A\,\|g\|_q,\qquad
p^{-1}+q^{-1}=1,\qquad 1<p<\infty.
$$

\subsection{Некоторые другие условия сходимости}
\label{paragraph3.3} Пусть $\rho=\rho_\infty$ --- бесконечная НД с
полюсами $z_k\in{\mathbb C^+}$, $p\in(1,\infty)$, причем
$\|\mu\|_p<\infty$. Из двойственности (\ref{(*8)}) получается
аналогичное (\ref{(*440)}) неравенство с заменой величины $A_0(p)$
на величину $A_1=A_1(\|\mu\|_p,p)$. Эта оценка и отмеченный в \S
\ref{paragraph2.1} результат А.\,Е.\,Додонова \cite{[DODONOV]}
уточняют результат В.\,Ю.\,Протасова (\ref{(*30)}).

Приведем еще некоторые результаты о сходимости из
\cite{[DAN2010]}. Если выполнено (\ref{(*31)}), то
последовательность $\rho_n(x)$ сходится к $\rho(x)$ {\it
равномерно} на $\mathbb{R}$, функция $\rho$ аналитична в полосе
$\{z:\;{\rm Im}\,z<A_0(p) \|\rho\|^{-q}_p\}$ и
$$
\|\rho\|_\infty \le A(p)\|\rho\|_p^q, \qquad |\rho'(x)| \le 2
\|\rho\|_\infty \cdot \mu(x),\qquad x\in \mathbb{R}.
$$
Кроме того, справедлива эквиваленция
$$ \|\mu\|_p<\infty\quad
\Leftrightarrow\quad
\sum\nolimits_{k=-\infty}^{\infty}\mu^{p}(c_k)<\infty\qquad
\forall c_k\in [k,k+1],
$$
из которой, в частности, вытекает импликация
$\|\mu\|_p<\infty\Rightarrow \|\mu\|_r<\infty$ при $p<r$. Из
импликации получаем: $S_p\subset S_r$ при $r>p$.

И.~Р.~Каюмов в \cite{[KAY2011],[KAY2012]} получил
необходимое условие для (\ref{(*31)}) в виде:
\begin{equation}
\|\mu\|_p<\infty\quad\Rightarrow \quad \sum_{k=1}^{\infty}
k^{p-1}|z_k^*|^{1-p}<\infty,
     \label{(*41)}
\end{equation}
где $z_k^*$ --– последовательность, полученная из $z_k$ путем
упорядочивания $|z_k|$ в порядке возрастания. Кроме того, им
получено достаточное для (\ref{(*31)}) условие в виде
$\sum_{k=1}^{\infty} k^{p-1}y_k^{1-p}<\infty$, $z_k=x_k+iy_k$,
$y_k>0$, которое в силу (\ref{(*41)}) оказывается также и
необходимым при условии, что $y_k$ упорядочены по возрастанию, а
$z_k\in {\mathbb C}^+$ лежат в некотором угле: $|z_k|<c\, y_k$.
Отметим, что при выполнении последнего свойства из (\ref{(*777)})
получается
$$
\|\mu\|_p<\infty \quad \Leftrightarrow\quad
\sum_{k=1}^{\infty}\left(\sum_{j=1}^{\infty
}\frac{1}{y_k+y_j}\right)^{p-1}<\infty,\qquad p\in (1,2).
 $$

\section{ Модификации и обобщения НД} \label{paragraph9}

Всюду в этом параграфе через $f$ и $h$ обозначены фиксированные
функции, аналитические в окрестности начала координат и имеющие
представления
\begin{equation}
f(z)=\sum_{m=0}^{\infty} f_m z^m, \qquad h(z)=\sum_{m=0}^{\infty}
h_m z^m, \qquad h_{m} \neq 0.
 \label{9.1-f-h}
\end{equation}
(Условие на коэффициенты $h$ накладывается здесь для простоты
изложения; в общем случае достаточно выполнения импликации
$f_m\neq 0 \Rightarrow h_m\neq 0$, см. \cite{[DAN2008]}.) Кроме
того, считаем, что $h$ аналитична в круге $D_r:=\{z:|z|\le r\}$ с
некоторым фиксированным $r>0$.

\subsection{Интерполяция Паде $h$-суммами, применение в численном
анализе}
\label{paragraph9.1}
В работе \cite{[DAN2008]} введены так называемые \textit{$h$-суммы}
вида
\begin{equation}
\label{9.1-h-sum}
{H}_n(z)={H}_n(\{\lambda_k\},h;z)=\sum_{k=1}^n\lambda_k
h(\lambda_k z), \qquad z, \lambda_k\in \mathbb{C},\qquad n\in
\mathbb{N},
\end{equation}
где  $h$  --- \textit{базисная} функция из (\ref{9.1-f-h}),
$|z|<r\,\min_{k=1,\ldots,n} |\lambda_k|^{-1}$. Если в качестве
базисной функции выбрать, например, $h(z)=(z-1)^{-1}$, то
$h$-сумма будет иметь вид НД с полюсами $\lambda_k^{-1}$ --- тем
самым $h$-суммы являются естественным обобщением НД.

В \cite{[DAN2008]} рассматривалась задача Паде-интерполяции
\begin{equation}
\label{9.1-iPade} f(z)-{H}_n(z)=O(z^n),\qquad z\to 0,
\end{equation}
которая, очевидно, обобщает рассмотренную в \S\ref{paragraph7.7}.
С учетом (\ref{9.1-f-h}) легко проверить, что (\ref{9.1-iPade})
равносильно системе уравнений для степенных сумм чисел
$\lambda_k$:
\begin{equation}
\label{9.1-sps}
S_{m+1}(\lambda_1,\ldots,\lambda_n):=\sum\nolimits_{k=1}^n
\lambda_k^{m+1}=s_m,\qquad s_m:=f_m/h_m, \quad m=\overline{0,n-1}.
\end{equation}
Как уже говорилось в \S\ref{paragraph7.1}, эта система всегда
разрешима, $\lambda_k$ являются корнями многочлена, построенного
как и (\ref{7.1-Tn}), но с заменой $f_m$ на $(-s_m)$.

В \cite{[DAN2008],Chu2010} были получены оценки величин
$|\lambda_k|$ и $|O(z^n)|$ при определенных ограничениях на числа
$|s_m|$ и даны оценки остаточных членов интерполяции. Например,
если для степенных сумм вида (\ref{9.1-sps}) имеем $|S_m|\le a^m$,
$m=\overline{1,n}$, то суммы ${H}_n$ сходятся к $f$ равномерно в
любом круге $|z|\le r(1-\delta)\,a^{-1}$, $\delta\in (0,1)$,
причем
$$
|{H}_n(z)-f(z)|\le (\theta\delta)^{-1}(1-\theta\delta)^n, \qquad
\theta\in (0,1),\qquad n\ge n_0.
$$
При этом радиус $a$ круга сходимости не может быть увеличен, а
число $1-\theta\delta$ в оценке не может быть заменено числом
меньшим $1-\delta$ (см. \cite{Chu2010}).

В дальнейшем $h$-суммы применялись как операторы численного
дифференцирования, интегрирования и экстраполяции на определенных
классах аналитических функций
\cite{DanChu2011,Chu2010,Chu2012,[FRN1],[FRN2],[FRN3]}. В этом
случае числа $\lambda_k$ уже не зависят от индивидуальных функций
и имеют универсальный характер. Например, в \cite{[DAN2008]}
найдены следующие точные на многочленах степени $\le n-1$ формулы
численного дифференцирования и интегрирования:
\begin{equation}
 zh'(z)\approx -h(z)+\sum_{k=1}^{n}\lambda_{1,k}
h(\lambda_{1,k} z);\quad \int_{0}^{z} h(t)\,dt\approx z
\sum_{k=1}^{n}\lambda_{2,k} h(\lambda_{2,k} z). \label{9.1-fDan}
\end{equation}
Здесь числа $\lambda_{l,k}$ --- абсолютные постоянные, являющиеся
корнями многочленов $P_{l,n}$ $(l=1,2)$, которые можно определить
рекуррентно следующим образом. Пусть $P_{l,0}=1$, $v_{l,1}=-1$
$(l=1,2)$. Тогда при $k=1,2,\ldots$ имеем
$$
P_{l,k}=\lambda P_{l,k-1}-v_{l,k},\quad
v_{1,k}=1+\sum_{j=1}^{k-1}\left(1-\frac{j}{k}\right)v_{1,j},\quad
v_{2,k}=\frac{1}{k^2}+\sum_{j=1}^{k-1}\frac{v_{2,j}}{k(k-j)}.
$$
Первая формула из (\ref{9.1-fDan}) была обобщена в \cite{Chu2010}
на случай любого порядка дифференцирования. В работах
\cite{[FRN1],[FRN2],[FRN3]} рассматривались аппроксимации
операторов более общего вида --- дифференциальных полиномов. Там
показано, что при фиксированном натуральном $q$ и любом
натуральном $n > q+5$ существуют комплексные числа $\lambda_k$,
$|\lambda_k|<1$, зависящие только от $q$ и $n$, такие, что имеет
место приближенное равенство
$$
\sum_{s=1}^q p_s h_{q-s} z^{q-s}\approx
\sum_{k=1}^{Nq}P(\lambda_k)\;h(\lambda_k z), \qquad N:=[n/q],
$$
где $P=\sum_{s}p_s \lambda^s$ --- произвольный многочлен степени
$\le q$. Погрешность формулы имеет порядок $o(n^{-n/q})$, $n\to
\infty$, при любом фиксированном $z$.

Еще одним примером использования $h$-сумм в численном анализе
является метод экстраполяции, предложенный в
\cite{DanChu2011,Chu2012}.  Показано, что при $a>1$ найдется
зависящее только от $a$ и $n$ множество $\{\lambda_k\}_{k=1}^n$ со
свойством
$$
\max_{k=\overline{1,n}}|\lambda_k|\le a-(a-1)n^{-1}, \qquad
S_m(\lambda_1,\ldots,\lambda_{n})=a^{m-1}, \qquad
m=\overline{1,n},
$$
откуда для функции $h$ получается {\it экстраполяционная} формула
$h(z)\approx {H}_n^{(\mu)}(z)$, где $\mu\in{\mathbb N}$ и
\begin{equation}
{H}_n^{(\mu)}(z)=\sum_{k_1,\ldots,k_\mu=1}^n\lambda_{k_1}
\cdots\lambda_{k_\mu}
h\left(\tfrac{\lambda_{k_1}\cdots\lambda_{k_\mu}}{a^\mu}z\right),
\quad
\left|\tfrac{\lambda_{k_1}\cdots\lambda_{k_\mu}}{a^\mu}z\right|
<\left(1-\tfrac{a-1}{an}\right)^{\mu}|z|.
 \label{9.1-ext1}
\end{equation}
Эта формула является экстраполяционной в том смысле, что значения
функции $h$ выражаются через ее значения в точках с меньшими
модулями.

Показано \cite{DanChu2011,Chu2012}, что при определенном $d_0$,
$d_0<r$, за счет сбалансированного выбора параметров $\mu$ и $n$
можно сколь угодно точно экстраполировать значения $h$ на
окружность $|z|=r$ из круга $|z|\le d_0$ (т.е. все узлы
экстраполяции $\lambda_{k_1}\cdots\lambda_{k_{\mu}}a^{-\mu}z$
лежат в круге $|z|\le d_0$). При этом для остаточного члена
$r_n(z):=h(z)-{H}_n^{(\mu)}(z)$ имеем \cite{Chu2012}:
\begin{equation}
\label{9.1-ext2} r_n(z)= \sum_{m=n}^{\infty}h_m
\left(1-\left(\frac{S_{m+1}}{a^{m}}\right)^{\mu}\right)z^m, \qquad
|r_n(z)| \le \sum_{m=n}^{\infty}|h_m| |z|^m, \qquad z\in D_r.
\end{equation}

Заметим еще, что при увеличении $\mu$ с фиксированным $n$ радиусы
кругов, в которых лежат узлы экстраполяции, стремятся к нулю как
геометрическая прогрессия (число узлов, естественно, возрастает),
см. (\ref{9.1-ext1}). Однако оценка погрешности из
(\ref{9.1-ext2}) не ухудшается, т.к. просто не зависит от $\mu$. В
\cite{Chu2012} приведены примеры, демонстрирующие определенные
преимущества указанной экстраполяции перед традиционной
интерполяцией (или экстра\-по\-ля\-цией) в корнях из единицы.

\subsection{Интерполяция амплитудно-частотными суммами.
Регуляризация} \label{paragraph9.2} В
\cite{[DanChu2014],[DanChu2013],[DanChu2013-2],[DanChu2016]}
предложено естественное обобщение $h$-сумм
--- \textit{амплитудно-час\-тот\-ные суммы} вида
$$
\mathcal{H}_n(\{\mu_k\},\{\lambda_k\},h;z):=\sum_{k=1}^n \mu_k
h(\lambda_k z),\qquad \mu_k,\lambda_k\in \mathbb{C},
$$
где \textit{амплитуды} $\mu_k$ и \textit{частоты} $\lambda_k$ ---
независимые друг от друга параметры. Как и $h$-суммы,
амплитудно-частотные суммы применялись и для аппроксимации
отдельных аналитических функций $f$ (и тогда
$\lambda_k=\lambda_k(f,h,n)$, $\mu_k=\mu_k(f,h,n)$), и как
специальные операторы (интегрирования, дифференцирования и
экстраполяции), действующие на определенном классе (и тогда
$\lambda_k=\lambda_k(n)$, $\mu_k=\mu_k(n)$) \cite{[DanChu2016]}.

Отметим, что амплитудно-частотные суммы являются естественным
обобщением некоторых классических аппаратов приближения, таких как
экспоненциальные суммы, тригонометрические полиномы, дроби Паде
(см. \cite{[DanChu2014],[DanChu2016]}).

В случае амплитудно-частотных сумм разумно ставить задачу уже
$2n$-кратной Паде-интерполяции функции $f$ в точке $z=0$:
\begin{equation}
f(z)=\mathcal{H}_n(\{\mu_k\},\{\lambda_k\},h;z)+O(z^{2n}), \qquad
z\to 0, \label{6-2n-int}
\end{equation}
что равносильно выполнению следующих условий на так называемые
\textit{обобщенные степенные суммы} (\textit{моменты}) (ср. с
(\ref{9.1-sps})):
\begin{equation}
\label{6-SRS}
  \mathcal{S}_m:=\sum_{k=1}^n \mu_k \lambda_k^m=s_m,\qquad s_m:=f_m/h_m,\qquad
m=\overline{0,2n-1}.
\end{equation}
Эту систему относительно неизвестных $\lambda_k,\mu_k$ при
известных правых частях $s_m$ называют \textit{задачей дискретных
моментов}.

Следуя \cite{Lyubich}, назовем совместную систему (\ref{6-SRS}) и
ее решение \textit{регулярными}, если все $\lambda_k$ попарно
различны, а все $\mu_k$ отличны от нуля.  Для решения регулярных
систем существует классический метод Прони (см.
\cite{Sylvester,Lyubich,Lyubich2}). Ключевую роль в этом методе
играет {\it производящий} многочлен
\begin{equation}
\label{6-G_n} G_n(\lambda):=\sum_{m=0}^n g_{m} \lambda^{m}= \left|
\begin{array}{ccccc}
1 & \lambda & \lambda^2 & \ldots & \lambda^n\\
s_0 & s_1 & s_2 & \ldots & s_n\\
s_1 & s_2 & s_3 & \ldots & s_{n+1}\\
\ldots & \ldots & \ldots & \ldots & \ldots\\
s_{n-1} & s_n & s_{n+1} & \ldots & s_{2n-1}\\
\end{array}
\right|,
\end{equation}
построенный по числам $s_m$, $m=\overline{0,2n-1}$, из
$(\ref{6-SRS})$. Регулярный случай характеризуется тем, что
степень многочлена $G_n$ равна $n$ и все его корни попарно
различны. Эти корни и дают нужные {частоты} $\lambda_k$. Затем
находятся и амплитуды $\mu_k$ из линейной относительно них системы
(\ref{6-SRS}).

Рассмотрим пример. Пусть $h$ --- аналитическая в круге $D_r$ функция и
$$
f(x):=\frac{1}{x}\int_{-x}^{x}h(t)\,dt,\qquad 0\le x <r.
$$
Несложно проверить, что в данном случае система (\ref{6-SRS})
регулярна, причем $s_{m-1}=({1-(-1)^{m}})/m$, $m=\overline{1,2n}$.
Ее решение методом Прони (с помощью производящего многочлена)
позволяет построить оператор
$\mathcal{H}_n(\{\mu_k\},\{\lambda_k\},h;x)$ численного
интегрирования --- квадратурную формулу Гаусса:
$$
f(x)\approx \mathcal{H}_n(\{\mu_k\},\{\lambda_k\},h;x),
$$
где частоты $\lambda_k$ вещественны, попарно различны и лежат на
интервале $(-1,1)$, а  амплитуды $\mu_k$ положительны. Хорошо
известно (см., например, \cite[Гл. 7, \S2]{Krylov}), что формулы
Гаусса имеют наивысший порядок точности среди всех квадратурных
формул порядка $n$ и точны на многочленах степени $2n-1$.

Возник вопрос о построении столь же высокоточных интерполяционных
формул для операторов численного дифференцирования и экстраполяции
\cite{[DanChu2014],[DanChu2016]}. Однако в этих случаях условия регулярности не
выполняются и, более того, для численного дифференцирования
соответствующая задача моментов вовсе не имеет решения.

Для преодоления этой трудности в \cite{[DanChu2014],[DanChu2016]} предложен
метод аналитической регуляризации задачи $2n$-кратной интерполяции
путем добавления к интерполируемой функции $f$ специального бинома
вида $b(z):=p\,h_{n-1}z^{n-1}+q\,h_{2n-1}z^{2n-1}$ с тем условием,
чтобы для функции $f+b$ задача (\ref{6-SRS}) была регулярной.
Тогда метод Прони дает <<подправленную>> интерполяционную формулу
\begin{equation} \label{6-f-reg-ex}
f(z)=-p\,h_{n-1}z^{n-1}-q\,h_{2n-1}z^{2n-1} +
\mathcal{H}_n(\{\mu_k\},\{\lambda_k\},h;z)+O(z^{2n}).
\end{equation}

В рассматриваемых в \cite{[DanChu2014],[DanChu2016]} приложениях
важно, чтобы $p$ и $q$ задавались явными формулами, позволяющими
получать оценки амплитуд, частот и остаточных членов интерполяции.
Остановимся на этих вопросах подробней.

\subsection{Численное дифференцирование
посредством амплитудно-частот\-ных сумм} \label{paragraph9.4}

Рассмотрим задачу $2n$-кратной интерполяции функции $zh'(z)$
посредством сумм $\mathcal{H}_n$, где в качестве базисной выбрана
сама функция $h$. Мы приходим к нерегулярной задаче дискретных
моментов (\ref{6-SRS}), где $s_m=m$. Применим метод регуляризации,
взяв
$$
q=q_0(p):=-2\,{\frac {p\left (3\,p+{n}^{2}-1\right )}{\left ({n}-1
\right )\left (n-2\right )}},\qquad p\in {\mathbb C}.
$$
Тогда производящий многочлен для $zh'(z)+b(z)$ примет следующий
вид
 \begin{equation}
  \label{6-FACG}
  {G}_n(\lambda)={g}_n \left({\lambda}^{n}-
\frac{6\lambda\left(\lambda^{n-1}-(n-1)\lambda+n-2\right)}{\left
(n- 1\right)\left(n-2\right)\left (\lambda-1\right)^{2}} +2+{\frac
{6\,p}{\left (n-1\right )\left (n-2\right)}} \right),
 \end{equation}
причем сколь угодно малой вариацией $p$ всегда можно добиться,
чтобы все $n$ корней этого многочлена были попарно различны и мы
имели регулярный случай. Тем самым, получается формула
(\ref{6-f-reg-ex}) с $f(z)=zh'(z)$, где частоты $\lambda_k$ ---
попарно различные корни многочлена $(\ref{6-FACG})$, а амплитуды
$\mu_k$ однозначно определяются по этим корням. Заметим, что
$\mu_k=\mu_k(p,n)$ и $\lambda_k=\lambda_k(p,n)$ не зависят от вида
аналитической функции $h$ и в этом смысле являются универсальными
(в \cite{[DanChu2014],[DanChu2016]} были получены некоторые оценки
$|\lambda_k|$, оценки соответствующих амплитуд --- вопрос
открытый).

Полученная формула численного дифференцирования точна на
многочленах степени не выше $2n-1$. При этом требуется знать
только $n$ значений функции $h$ и два фиксированных значения ее
производных в точке $z=0$. Традиционный интерполяционный подход
при таком количестве известных значений позволяет получить, вообще
говоря, только порядок $O(z^{n+2})$ (см.
\cite{Ash_Janson_Jones,Salzer,Lyness,Lyness2,[DAN2008]}).

\subsection{Экстраполяция
посредством амп\-ли\-туд\-но-час\-тот\-ных сумм}
\label{paragraph9.5} Экстраполяция посредством сумм
$\mathcal{H}_n$ заметно улучшает качество экстраполяции
$h$-сум\-ма\-ми (см. \S\ref{paragraph9.1}). Пусть $a>0$, а $h$ ---
функция, аналитическая в некотором круге $|z|<r$. Если функцию
$h(az)$ интерполировать посредством амплитудно-частотной суммы
${\mathcal{H}}_n(\{\mu_k\},\{\lambda_k\},h;z)$ (в качестве
базисной выбрана сама $h$), то, как и в случае дифференцирования,
мы получаем нерегулярную задачу моментов с $s_m=a^m$. В
\cite{[DanChu2014],[DanChu2016]} показано, что она регуляризуется, например,
добавлением бинома вида $b(z)$ с $p>0$, $q=0$. При этом
производящий многочлен принимает вид
$$
{G}_n(\lambda)={{g}}_n \left(\lambda^n-\frac{a^{n}}{na^{n-1}+p}
\frac{\lambda^n-a^n}{\lambda-a}\right)
$$
и всегда имеет $n$ попарно различных корней
$\{\lambda_k\}_{k=1}^n$. Таким образом, при указанных параметрах
$p$ и $q$ справедлива интерполяционная формула
\begin{equation}
\label{9.1-AF-extr} h(az)=-
p\,h_{n-1}z^{n-1}+\mathcal{H}_n(\{\mu_k\},\{\lambda_k\},
h;z)+r_n(z),\qquad r_n(z)=O(z^{2n}),
\end{equation}
точная на многочленах степени не выше $2n-1$. Здесь, как и в
предыдущей задаче, частоты и амплитуды не зависят от $h$. В
(\ref{9.1-AF-extr}) для всех частот справедливы неравенства
$$
|\lambda_k|<\delta\,a;\qquad
\delta:=\left(1+p/(na^{n-1})\right)^{-1/n}<1,
$$
т.е. формула (\ref{9.1-AF-extr}) является экстраполяционной. При
этом для остаточного члена экстраполяции из (\ref{9.1-AF-extr})
справедлива следующая, не зависящая от $p$, оценка
$$
|r_n(z)|\le \sum_{m=2n}^\infty |f_m||az|^m.
$$
Из этого, в свою очередь, вытекает, что теоретическая погрешность
экстраполяции не возрастает при изменении параметра $p$. С другой
стороны, в силу свойств величины $\delta=\delta(n,a,p)$ при
возрастании $p$ и фиксированных прочих параметрах узлы
экстраполяции стягиваются к точке $z=0$ --- возникает интересное
явление стягивания узлов экстраполяции в одну точку без влияния на
оценку погрешности. Аналогичный эффект уже отмечался в сходных
задачах для $h$-сумм \cite{DanChu2011,Chu2012} (также см. выше~\S\ref{paragraph9.1}).

В заключение отметим, что обычно при $n$-точечной простой или
кратной экстраполяции на основе многочленов Лагранжа и других
сходных аппаратов получаются экстраполяционные формулы, точные на
многочленах порядка не выше $n-1$ --- см. такие формулы, например,
в \cite{Salzer2,DanChu2011,Chu2012}.

Полученные в \cite{[DanChu2014],[DanChu2016]} $n$-точечные
экстраполяционные формулы точны на многочленах порядка $2n-1$,
причем такое удвоение порядка точности достигается за счет
введения лишь одного регуляризующего слагаемого $ph_{n-1}z^{n-1}$.

\subsection{Замечание о вещественных формулах}
\label{paragraph100.5}

При работе с вещест\-вен\-но\-знач\-ны\-ми на $\mathbb{R}$
аналитическими функциями (которых в приложениях большинство) к
определенному неудобству обсуждавшихся в \S\ref{paragraph9}
формул можно отнести то, что среди подходящих значений параметров
~$\lambda_k,\mu_k$ для $h$-сумм и амплитудно-частотных сумм могут
быть невещественные. Оперирование комплексными значениями в таких
задачах весьма неудобно, например, если известны лишь
табулированные значения функций на действительной оси. Возникает
задача выделения достаточно широкого подкласса функций $h$, в
котором возможно построение {\it вещественных} интерполяционных
формул того же типа, и самого построения таких формул.

Эта задача непроста и весьма далека от полного решения. Наибольшим
продвижением в этом направлении на данный момент является работа
Ю.\,М.~Нигматяновой \cite{[Nigm2017]}, в которой предложен численный
метод пробных алгебраических многочленов, с помощью которого
строятся вещественные интерполяционные формулы рассматривавшихся в
\S\ref{paragraph9.1} типов. Например, для четных аналитических в
окрестности начала функций $h(z)$ построены операторы вида
$\sum\nolimits_{k=1}^n\lambda_kh(\lambda_kz)$ нечетного порядка
$n$ с вещественными параметрами $\lambda_k$, аппроксимирующие
дифференциальный оператор $(zh(z))'$ с локальной погрешностью
$O(z^{n+2})$ $(z\to 0)$ при $n\le 51$.

\subsection{Вещественные амплитудно-фазовые суммы}
\label{paragraphAFO}

Сформулируем основной результат из \cite{[Dan-Vas1],[Dan-Vas2],[Dan-Vas3],[Dan-Vas4]} о
выделении гармоник из тригонометрического многочлена
$T_{n}(t)=a_0+\sum_{k=1}^{n} \tau_k(t)$.\smallskip

{\it Пусть натуральные $s\ge 2$, $\mu\ge 1$, $n=s\mu -2$,
$m=\mu(s-1)$. Тогда
\begin{equation} \tau_{\mu}(t)+ a_0\,\omega=
\sum_{j=1}^{m} X_j\cdot T_{n}\left(t-\lambda_j\right),\quad
m=\mu(s-1),
  \label{AFO2}
 \end{equation}
где
$$\omega=-2\cos \varphi_{\alpha},\quad
\varphi_{\alpha}:=\frac{\pi {\alpha}}{s+1},\quad
{\alpha}=\overline{1,s},\quad \lambda_k=-\arg z_k,
$$
а значения $z_k$ $( k=\overline{1,m})$ попарно различны  и
составляют множество
 $$
\{z_k\}=\{\sqrt[(s+1)\mu]{(-1)^{\alpha}}\}\setminus
\{\sqrt[\mu]{e^{i\varphi_{\alpha}}},
\sqrt[\mu]{e^{-i\varphi_{\alpha}}}\};
$$
$$
X_k=\frac{1}{(s+1)\mu}\left(\omega+2\,{\rm
Re}\,z_k^{\mu}\right),\quad \sum_{k=1}^{m} X_k=\omega;
$$
где $ \hbox{все } X_k<0$, если $\alpha=1$ и все $X_k>0$,
 если $\alpha=s$.
Кроме того, формула $(\ref{AFO2})$ верна и для сходящихся
тригонометрических рядов, в которых отсутствуют гармоники с
номерами}
$$
\beta=\mu+{\theta}\,k,\quad n+1+{\theta}\,(k-1)\le \beta\le
\theta+{\theta}(k-1),\quad \hbox{где}\quad {\theta}=n+\mu+2,\quad
k\in {\mathbb N}.
$$

Выделение гармоник наложением сигналов (без использования
промежуточных спектральных замеров) позволяет эффективно применять
амплитудно-фа\-зо\-вые суммы для оценок гармоник
тригонометрических полиномов. Отметим, что оценкам коэффициентов
тригонометрических многочленов посвящено много работ (см.,
например, работы А.\,С. Белова и С.\,В. Конягина \cite{BELOV1,BELOV2} и обширную библиографию в них).

Особенно важным для оценок является построение амплитудно-фазовых
сумм, у которых все $X_k$ одного знака (это соответствует
значениям $\alpha=1$, $\alpha=s$ в предыдущей теореме). С помощью
таких сумм, например, получаются точные оценки $L_p$-норм гармоник
через $L_p$-норму самого многочлена, а в случае неотрицательных
полиномов $T_{n}\not\equiv 0$ --- точные оценки их коэффициентов
через свободный член (обобщение неравенства Фейера для первой
гармоники). Из предыдущей теоремы вытекают следующие оценки
\cite{[Dan-Vas1]}:\smallskip

{\it В условиях теоремы о выделении гармоник при $1\le p\le\infty$
имеем оценку
$$
\|a_0\pm\omega^{-1}\tau_{\mu} (t)\|_{L_p[0,2\pi]}\le
 \|T_{n}\|_{L_p[0,2\pi]},\qquad \hbox{где}\quad
 \omega=2\cos\frac{\pi}{s+1}.
$$
В частности, для вещественного многочлена имеем
 $$
 |a_0|+\omega^{-1}\sqrt{a_{\mu}^2+b_{\mu}^2}\le
 \|T_{n}\|_{\infty}
 $$
 $($при $s=2$ достигается равенство на многочлене
$a_0+\tau_{\mu}(t)$$)$. Кроме того, для вещественного
неотрицательного многочлена $T_{n}\not\equiv 0$ имеем
$$
a_0 - \omega^{-1} \sqrt{a_{\mu}^2+b_{\mu}^2}\ge \min_t
T_{n}(t),\qquad \sqrt{a_{\mu}^2+b_{\mu}^2}\le \omega\,a_0< 2\,a_0
 $$
 $($при $\mu=1$ последнее неравенство принадлежит Л.~Фейеру.$)$
Множитель $2$ в неравенстве нельзя заменить меньшим. Например, для
последовательности неотрицательных четных многочленов
$T^*_{n}(t)=(1+\cos t)^{n}$ при всех $\mu$ имеем}
$$
\frac{a_{\mu}(T^*_{n})}{2a_0(T^*_{n})}
=\frac{\prod_{k=0}^{\mu-1}(n-k)}{\prod_{k=1}^{\mu}(n+k)}\to
  1,\quad n\to \infty.
 $$

Отметим, что в \cite{[Dan-Dan],[Dan-Dan+]} рассматривалась задача
выделения суммы двух гармоник из тригонометрических многочленов
амплитудно-фазовой суммой; в этом случае также получены точные
оценки типа Фейера.


\begin{thebibliography}{100}


\bibitem{[AND-EID2]} {Anderson J.\,M., Eiderman V.\,Ya.} {\it Cauchy transforms of
point masses: the logarithmic derivative of polynomials}. Ann. of
Math., {\bf 163}:3 (2006), 1057--1076.

\bibitem{Ash_Janson_Jones}
{Ash J.\,M., Janson S., Jones R.\,L.} {\it Optimal numerical
differentiation using $N$ function evaluations}. Calcolo,
\textbf{21}(2) (1984), 151--169.

\bibitem{[BOOL]} {Boole G.} {\it On the comparison of transcendents, with certain
applications to the theory of definite integrals}. Philos. Trans.
R. Soc., {\bf 147} (1857), 745–-803.

\bibitem{[cartan]} {Cartan H.} {\it Sur les syst\`emes de fonctions holomorphes \`a vari\'et\'es lin\'eaires
lacunaires et leurs applications}. Ann. Sci. \'Ecole Norm. Sup.,
{\bf 45}:3 (1928), 255–-346.

\bibitem{[Chui]} {Chui C.\,K.} {\it On approximation in the Bers spaces}. Proc.
Amer. Math. Soc., {\bf 40} (1973), 438--442.

\bibitem{[ChuiShen]} {Chui C.\,K., Shen X.\,C.} {\it Order of approximation by
electrostatic fields due to electrons}. Constr. Approx., {\bf 1}:1
(1985), 121--135.

\bibitem{[Chu2014]} {Chunaev P.} {\it Least deviation of logarithmic
derivatives of algebraic polynomials from zero}. J. Approx. Theory, {\bf 185} (2014), 98–-106.

\bibitem{[DanChu2014]} {Chunaev P.\,V., Danchenko V.\,I.} {\it Approximation
by the amplitude and frequency operators}. arXiv:1409.4188 (2014).

\bibitem{[DanChu2016]} {Chunaev P., Danchenko V.}
 {\it Approximation by amplitude and frequency operators.}  J. Approx Theory. 2016. V. 207. p. 1-30.

\bibitem{[DanChu2016+]} {Chunaev P.\,V., Danchenko V.\,I.} {\it Sharp inequalities of
Jackson-Nikolskii type for rational functions.} arXiv:1611.03485 (2016).

\bibitem{[D-D3]}
{Danchenko V.\,I., Danchenko D.\,Ya.} {\it Nikolskii type inequalities for
simple partial fractions}. Комплексный анализ и его приложения:
материалы VII Петрозаводской международной конференции (29 июня -
5 июля 2014 г.) (2014), 33-37.

\bibitem{[DanDod2013]} {Danchenko V.\,I., Dodonov A.\,E.} {\it Estimates for
exponential sums. Applications}. J. Math. Sci., {\bf 188}:3
(2013), 197--206.

\bibitem{DanChu2011}
{Danchenko V.\,I., Chunaev P.\,V.} {\it Approximation by
simple partial fractions and their generalizations}. J. Math.
Sci., {\bf 176}:6 (2011), 844--859.

\bibitem{[DanChu2013]}
{Danchenko V.\,I., Chunaev P.\,V.} \textit{Approximation
by amplitude and frequency sums}, Joint CRM-ISAAC Conference on
Fourier Analysis and Approximation Theory: Abstracts (Bellaterra,
2013) (2013), 12.

\bibitem{[Dan-Vas1]} {Danchenko V.\,I., Vasilchenkova D.\,G.}
{\it Extraction of harmonics from trigonometric polynomials by
amplitude and phase operators},  arXiv:1606.08716 (2016).

\bibitem{[DODONOV]} {Dodonov A.\,E.} {\it
On Convergence of Series of Simple Partial Fractions in
$L_p(\mathbb R)$}. J. Math. Sci., {\bf 210}:5 (2015), 648–653.

\bibitem{[FRN1]}
{Fryantsev A.\,V.} {\it On numerical approximation of differential polynomials}.
J. Math. Sci., {\bf 157}:3 (2009), 395--399.

\bibitem{[KOR64]} {Korevaar J.} {\it Asymptotically neutral distributions of electrons and
polynomial approximation}. Ann. of Math. (2), {\bf 80}:2 (1964),
403--410.

\bibitem{[KOR]} {Korevaar J.} {\it Limits of polynomials whose zeros lie in a given
set}. Proc. Sympos. Pure Math. Amer. Math. Soc., {\bf 11} (1968), 261--272.

\bibitem{[Kom2011]}
{Komarov M.\,A.} {\it Uniqueness of a simple partial fraction
of best approximation}. J. Math. Sci., {\bf 175:}3 (2011),
284–308.

\bibitem{[Kom2012]}
{Komarov M.\,A.} {\it Interpolation of rational functions by
simple partial fractions}. J. Math. Sci., {\bf 181}:5 (2012),
600--612.

\bibitem{[Kom2012-2]}
{Komarov M.\,A.} {\it Examples related to best approximation
by simple partial fractions}. J. Math. Sci., {\bf 184}:4 (2012),
509--523.

\bibitem{[Kom2013]}
{Komarov M.\,A.} {\it Sufficient condition for the best
uniform approximation by simple partial fractions}. J. Math. Sci.,
{\bf 189}:3 (2013), 482--489.

\bibitem{[Kom2014-1]}
{Komarov M.\,A.} \textit{On the analog of Chebyshev theorem
for simple partial fractions}. Комплексный анализ и его
приложения: материалы VII Петрозаводской международной конференции
(Петрозаводск, 29 июня --- 5 июля 2014 года) (2014), 65--69.

\bibitem{[Kom2015-PMA]}
{Komarov M.\,A.} {\it An Analog of the Haar Condition for
Simple Partial Fractions}. J. Math. Sci., {\bf 208}:2 (2015),
174--180.

\bibitem{[Kom-2016]}
{Komarov M.\,A.} \textit{Alternance and best uniform
approximation of odd functions by simple partial fractions}.
Комплексный анализ и его приложения: материалы VIII Петрозаводской
международной конференции (Петрозаводск, 3
--- 9 июля 2016) (2016), 41--46.

\bibitem{Kung1}
{Kung J.\,P.\,S.} {\it Canonical forms of binary forms:
Variations on a theme of Sylvester. Invariant theory and tableaux
(Minneapolis, MN, 1988)}. IMA Vol. Math. Appl., \textbf{19}
(1990), 46–58.

\bibitem{Lyness}
{Lyness J.\,N.} {\it Differentiation formulas for analytic
functions}. Math. Comp., \textbf{22} (1968), 352--362.

\bibitem{Lyness2}
{Lyness J.\,N., Moler C.\,B.} {\it Numerical differentiation
of analytic functions}. SIAM J. Numer. Anal., \textbf{4} (1967),
202--210.

\bibitem{Lyubich2}
{Lyubich Y.\,I.} {\it Gauss type complex quadrature formulae, power moment problem and elliptic curves}.
Матем. физ., анал., геом., \textbf{9}(2) (2002), 128–145.

\bibitem{Lyubich}
{Lyubich Y.\,I.} {\it The Sylvester-Ramanujan system of equations and
the complex power moment problem}. Ramanujan J., \textbf{8}
(2004), 23–45.

\bibitem{[M-F]} {Macintyre A.\,J., Fuchs W.\,H.\,J.} {\it Inequalities for the
logarithmic derivatives of a polynomial}. J. London Math. Soc.,
{\bf 15}:3 (1940), 162--168.

\bibitem{[MAR]} {Marstrand J.\,M.} {\it The distribution of the logarithmic
derivative of a polynomial}. J. London Math. Soc., {\bf 38}
(1963), 495--500.

\bibitem{[Nigm2017]} {Nigmatyanova Yu.\,M.} {\it Numerical Analysis
of the Method of Differentiation by Means of Real h-Sums}. J.
Math. Sci., {\bf 224}:5 (2017), 735--743.

\bibitem{Prony}
{Prony R.} {\it Sur les lois de la Dilatabilit\'{e} des
fluides \'{e}lastiques et sur celles de la Force expansive de la
vapeur de l'eau et de la vapeur de l'alkool, \`{a} diff\'{e}rentes
temp\'{e}ratures}. J. de l'Ecole Polytech., \textbf{2}(4) (1795),
28--35.

\bibitem{Ramanujan} {Ramanujan S.} {\it Note on a set of simultaneous
equations}. J. Indian Math. Soc., \textbf{4} (1912), 94–96.

\bibitem{Salzer}
{Salzer H.\,E.} {\it Optimal points for numerical
differentiation}. Num. Math., \textbf{2}:1 (1960), 214--227.

\bibitem{Salzer2}
{Salzer H.\,E.} {\it Formulas for best extrapolation}. Num.
Math., \textbf{18} (1971), 144--153.

\bibitem{Sylvester} {Sylvester J.\,J.}
{\it On a remarkable discovery in the theory of canonical forms
and of hyperdeterminants}. Phil. Magazine, \textbf{2} (1851),
391–410.

\bibitem{[Dan-Vas3]} {Vasilchenkova D. G., Danchenko V. I.} {\it Extraction
of harmonics from trigonometric polynomials  by amplitude and
phase transformation} Комплексный анализ и его приложения:
материалы VIII Петрозаводской международной конференции (3 - 9
июля 2016 г.). Петрозаводск:  Изд-во ПетрГУ, 2016. 93-95.

\bibitem{[WARSH]}
{Warschawski S.\,E.} {\it On differentiability at the boundary
in conformal mapping}. Proc. Amer. Math. Soc., {\bf 12}:4 (1961),
615--620.



\bibitem{[AND-EID1]} {Андерсон Дж.\,M., Эйдерман В.\,Я.} {\it Оценки преобразования
Коши точечных масс (логарифмической производной многочлена)}.
Докл. РАН, {\bf 401}:5 (2005), 583--586.

\bibitem{[Akhiezer]} {Ахиезер Н.\,И.} {\it Лекции по теории аппроксимации}. Наука, М. 1965.

\bibitem{Bari}
{Бари Н.\,К.} {\it Обобщение неравенства А.~А.~Маркова и
С.~Н.~Бернштейна}. Изв. АН СССР. Сер. матем., {\bf 18} (1954),
159--176.

\bibitem{BELOV1}
{ Белов А.С., Конягин С.В.} {\it Об оценке свободного члена
неотрицательного тригонометрического полинома с целыми
коэффициентами.} Изв. РАН. Сер. матем. 1996. 60:6. 31–90.

\bibitem{BELOV2}
{Белов А.С.} {\it Некоторые свойства и оценки для
неотрицательных тригонометрических полиномов.} Изв. РАН. Сер.
матем. 2003. 67:4. 3–20.

\bibitem{[BOR2007]} {Бородин П.\,А.} {\it Оценки расстояний до прямых и лучей от полюсов
наипростейших дробей, ограниченных по норме $L_p$ на этих
множествах}. Матем. заметки, {\bf 82}:6 (2007), 803–-810.

\bibitem{[BOR2009]}
{Бородин П.\,А.} {\it Приближение наипростейшими дробями на
полуоси}. Матем. сб., {\bf 200}:8 (2009), 25--44.

\bibitem{[BOR-KOS]} {Бородин П.\,А., Косухин О.\,Н.}
{\it О приближении наипростейшими дробями на действи\-тельной
оси}. Вестн. МГУ. Сер.~1. Матем., мех., {\bf 1} (2005), 3--8.

\bibitem{[BOR2012]} {Бородин П.\,А.}
{\it  Приближение наипростейшими дробями с ограничением на
полюсы}. Матем. сб., {\bf 203}:11 (2012),  23–40.

\bibitem{[BOR2016]} {Бородин П.\,А.}
{\it Приближение наипростейшими дробями с ограничением на
полюсы. II}. Матем. сб., {\bf 207}:3 (2016), 19--30.

\bibitem{[BULANOV]} {Буланов А.\,П.} {\it Асимптотика для наименьших уклонений функции ${\rm sign}\,x$ от рациональных
функций}.  Матем. сб., {\bf 96(138)}:2 (1975),  171--188.

\bibitem{[Dan-Vas2]} {Васильченкова\,Д.\,Г., Данченко\,
В.\,И.} {\it Фильтрация тригонометрических многочленов
амплитудно-фазовыми операторами.}  Междунар. конф. по
дифференциальным уравнениям и динамическим системам. Тезисы
докладов. Суздаль, 8-12 июль 2016.~ 40-41.

\bibitem{[Dan-Vas4]} {Васильченкова\,Д.\,Г., Данченко\,
В.\,И.} {\it Оценки гармоник тригонометрических многочленов.}
Междун. конф. по математической теории управления и
механике.  Тезисы   докладов. Суздаль,  7--11 июля 2017. 48-49.

\bibitem{[GARNET]} {Гарнетт Дж.}
{\it Ограниченные аналитические функции}. Мир, М., 1984.

\bibitem{[GEL]} {Гельфонд А.\,О.} {\it Об оценке мнимых частей корней многочленов с
ограниченными производными от логарифмов на действительной оси}.
Матем. сб., {\bf 71(113)}:3 (1966), ~~~~~~~289--296.

\bibitem{[GOVOR-LAP]} {Говоров Н.\,В., Лапенко Ю.\,П.} {\it Оценки снизу модуля логарифмической производной
многочлена}. Матем. заметки, {\bf 23}:4 (1978), 527--535.

\bibitem{[GOR]} {Горин Е.\,А.} {\it Частично гипоэллиптические дифференциальные
уравнения в частных производных с постоянными коэффициентами}.
Сиб. матем. журн., \textbf{3}:4 (1962), 500--526.

\bibitem{[GON3]} {Гончар А.\,А.} {\it О наилучших приближениях рациональными
функциями}. Докл. АН СССР, {\bf 100}:2 (1955), 205–-208.

\bibitem{[GON2]} {Гончар А.\,А.} {\it Обратные теоремы о наилучших приближениях на
замкнутых множествах}. Докл. АН СССР, {\bf 128}:1 (1959), 25–-28.

\bibitem{[GON1]} {Гончар А.\,А.} {\it Обратные теоремы о наилучших приближениях
рациональными функциями}. Изв. АН СССР. Сер. матем., {\bf 25}:3
(1961), 347--356.

\bibitem{[GON4]} {Гончар А.\,А.} {\it О задачах Е.~И.~Золотарева,
связанных с рациональными функциями}. Матем. сб., {\bf 78(120)}:4
(1969), 640–-654.

\bibitem{[GON5]} {Гончар А.\,А., Григорян Л.\,Д.} {\it Об оценках норм голоморфной
составляющей мероморфной функции}. Матем. сб., {\bf 99} (1976),
634--638.

\bibitem{[GON6]} {Гончар А.\,А., Григорян Л.\,Д.} {\it Об оценке компонент ограниченных
аналитических функций}. Матем. сб., {\bf 132} (1987), 299--303.

\bibitem{[GRIGOR]} {Григорян Л.\,Д.} {\it Оценка нормы голоморфных составляющих мероморфных функций в
областях с гладкой границей}. Матем. сб., {\bf 100(142)}:1(5)
(1976), 156–-164.

\bibitem{[DAN-razdel]}
{Данченко В.\,И.}
{\it О разделении особенностей мероморфных функций}. Матем. сб.,
{\bf 125(167)}:2(10) (1984), 181–198.

\bibitem{[DAN-viniti-1]}
{Данченко В.\,И.} {\it Оценки отклонения от действительной оси нулей
многочлена с нормированной логарифмической производной}. Деп.
ВИНИТИ, № 2990-91 (1991), 1--21.

\bibitem{[DAN-viniti-2]}
{Данченко В.\,И.} {\it Оценки мнимых частей полюсов логарифмических
производных многочленов}. Деп. ВИНИТИ, № 1695-92 (1992), 1-19.

\bibitem{[DAN-DAN1]}
{Данченко В.\,И.} {\it О скорости приближения к действительной оси
полюсов нормированных логарифмических производных полиномов}.
Докл. РАН, {\bf 330}:1 (1993), 15--16.

\bibitem{[DAN1]}
{Данченко В.\,И.} {\it Оценки расстояний от полюсов логарифмических
производных многочленов до прямых и окружностей}. Матем. сб., {\bf
185}:8 (1994), 63–80.

\bibitem{[DAN-DISS]}
{Данченко В.\,И.} {\it Метрические свойства мероморфных
функций. Дисс. . . . докт. физ.-матем. наук}. МГУ, М., 1999.

 \bibitem{[DAN2006]}
{Данченко В.\,И.}
{\it Оценки производных наипростейших дробей и другие вопросы}.
Матем. сб., {\bf 197}:4 (2006), 33–52.

\bibitem{[DAN2008]}
{Данченко В.\,И.} {\it Об аппроксимативных свойствах
сумм вида $\sum_k\lambda_kh(\lambda_k z)$}. Матем. заметки, {\bf
83}:5 (2008), 643–649.

\bibitem{[DAN2010]}
{Данченко В.\,И.} {\it О сходимости наипростейших дробей в
$L_p(\mathbb R)$}. Матем. сб., {\bf 201}:7 (2010), 53--66.

\bibitem{[DAN2013]}
{Данченко В.\,И.}
{\it Интегральные оценки длин линий уровня рациональных функций и
задача Е. И. Золотарева}. Матем. заметки, {\bf 94}:3 (2013),
331–337.

\bibitem{[D-D1]}
{Данченко В.\,И., Данченко Д.\,Я.} {\it О равномерном приближении
логарифмическими производными многочленов}. Теория функций, ее
приложения и смежные вопросы: Материалы школы-конференции,
посвященной 130-летию Д.\,Ф.~Егорова
$($Казань, 1999$)$ (1999), 74-77.

\bibitem{[D-D2]}  {Данченко В.\,И., Данченко Д.\,Я.} {\it О~приближении
наипростейшими дробями}. Матем. заметки, {\bf 70}:4 (2001),
553--559.

\bibitem{[D-D.Abstr2010]}
{Данченко В.\,И., Данченко Д.\,Я.} {\it О единственности
наипростейшей дроби наилучшего приближения}. Международная
конференция по дифференциальным уравнениям и динамическим системам
(Суздаль, 2010): Тезисы докладов (2010), 71-72.

\bibitem{[Dan-Dan+]} {Данченко В.И., Данченко Д.Я.}
 {\it  Оценки сумм двух гармоник тригонометрических многочленов.}
Международная конференция по математической теории управления и
механике. Тезисы   докладов. Суздаль,  7--11 июля 2017. 62-63.

\bibitem{[Dan-Dan]} {Данченко В.И., Данченко Д.Я.}
 {\it  Выделение пар гармоник  из
тригонометрических многочленов амплитудно-фазовым оператором.}
XIII Международная Казанская летняя школа-конференция <<Теория
функций, ее приложения и смежные вопросы>>.  Казань.  21-27
августа 2017, 143-145.

\bibitem{DAN-DOD}
{Данченко В.\,И., Додонов А.\,Е.} {\it Оценки $L_p$-норм наипростейших
дробей}. Изв. вузов. Матем., {\bf 6} (2014), 9–19.

\bibitem{[DAN-KOND2]}
{Данченко В.\,И., Кондакова Е.\,Н.}
{\it Чебышевский альтернанс при аппроксимации констант
наипростейшими дробями}. Тр. МИАН, {\bf 270} (2010), 86–96.

\bibitem{[DAN-KOND1]}
{Данченко В.\,И., Кондакова Е.\,Н.} {\it Критерий возникновения особых
узлов при интерполяции наипростейшими дробями}. Тр. МИАН, {\bf
278} (2012), 49–58.

\bibitem{[DAN-SEMIN]}
{Данченко В.\,И., Семин Л.\,А.} {\it  Точные квадратурные
формулы и неравенства разных метрик для рациональных функций.}
Сиб. матем. журн., {\bf 57}:2 (2016), 282-296.

\bibitem{[DANChu2013+]}
{Данченко В.\,И., Чунаев П.\,В.}
{\it Об оценке мнимых частей полюсов наипростейших дробей с
нормированной производной на действительной оси}. Современные
методы теории функций и смежные проблемы: Материалы Воронежской
зимней математической школы (2013), 78.

\bibitem{[DanChu2013-2]} {Данченко В.\,И., Чунаев П.\,В.} {\it Об аппроксимации посредством
частотно-амплитудных сумм}. Международная Казанская летняя
школа-конференция (Казань, 22-28 август 2013 г.). Теория функций, ее приложения и смежные
вопросы \textbf{46}  (2013), 174-175.

\bibitem{DDJ}
{Данченко Д.\,Я.} {\it Некоторые вопросы аппроксимации и интерполяции
ра\-ци\-она\-ль\-ными функциями. Приложение к уравнениям
эллиптического типа. Дисс. . . . канд. физ.-матем. наук}.  Вла\-ди\-мир\-ский государственный педагогический ун-т (2001).

\bibitem{[DEMIDOVICH]} {Демидович Б.\,П., Марон И.\, А} {\it
Основы вычислительной математики}. Государственное издательство
Физ.-мат. лит. Москва. 1960.

\bibitem{[Dodonov]} {Додонов\,
А.\,Е.} {\it Неравенства разных метрик для наипростейших дробей.}
Международная конференция по математической теории управления и
механике. Тезисы   докладов. Суздаль,  7--11 июля 2017. С. 65.

\bibitem{[DOL2]} {Долженко Е.\,П.} {\it Скорость приближения рациональными дробями и
свойства функций}. Матем. сб., {\bf 56}:4 (1962), 403--432.

\bibitem{[DOL3]} {Долженко Е.\,П.} {\it Рациональные аппроксимации и граничые свойства
аналитических функций}. Матем. сб., {\bf 69(111)} (1966),
497--524.

\bibitem{[DOL4]} {Долженко Е.\,П.} {\it О зависимости граничных свойств аналитической
функции от скорости ее приближения рациональными функциями}.
Матем. сб., {\bf 103(145)} (1977), ~~~~~~~131--142.

\bibitem{[DOL1]} {Долженко Е.\,П.} {\it Некоторые точные интегральные оценки
производных рациональных и алгебраических функций. Приложения}.
Anal. Math., {\bf 4}:4 (1978), 247–-268.

\bibitem{[ZOLOT]} {Золотарев Е.\,И.} {\it Полное собрание сочинений. Вып.
2}. Физ.-мат. ин-т им.~В.~А.~Стеклова АН СССР., Изд-во АН СССР,
Л., 1932.

\bibitem{[KAC]} {Кацнельсон В.\,Э.} {\it О некоторых операторах, действующих в
пространствах, порожденных функциями $\frac{1}{z-z_k}$}. Теория
функций, функциональный анализ и их приложения, {\bf 4} (1967),
58--66.

\bibitem{[KAY2011]}
{Каюмов И.\,Р.} {\it Сходимость рядов наипростейших дробей в $L_p(\mathbb
R)$}. Матем. сб., {\bf 202}:10 (2011), 87--98.

\bibitem{[KAY2012]}
{Каюмов И.\,Р.} {\it Необходимое условие сходимости наипростейших дробей в
$L_p(\mathbb R)$}. Матем. заметки, {\bf 92}:1 (2012), 149–152.

\bibitem{[KAY2012+]}
{Каюмов И.\,Р.} {\it Интегральные оценки наипростейших дробей}. Изв.
вузов. Матем., {\bf 4} (2012), 33--45.


\bibitem{[KAYva]}
{Каюмова А.\,В.} {\it Сходимость рядов простых дробей в $L_p$}. Учeн. зап.
Казан. гос. ун-та. Сер. Физ.-матем. науки, {\bf 154}:1 (2012),
208--213.

\bibitem{[Kolmogorov48]}
{Колмогоров А.\,Н.} {\it Замечание по поводу многочленов
П.\,Л. Чебышева, наимение уклоняющихся от заданой функции}. Успехи
матем. наук, {\bf 3}:1 (1948), 216–-221.

\bibitem{[KOM2]}
{Комаров М.\,А.} {\it О неединственности наипростейшей дроби наилучшего
равномерного приближения}. Изв. вузов. Матем., {\bf 9} (2013),
28–37.

\bibitem{[KOM3]}  {Комаров М.\,А.} {\it Критерий наилучшего приближения констант
наипростейшими дробями}. Матем. заметки, {\bf 93}:2 (2013),
209–215.

\bibitem{[KOM1]}
{Комаров М.\,А.}
{\it Критерий разрешимости задачи кратной интерполяции посредством
наипростейших дробей}. Сиб. матем. журн., {\bf 55}:4 (2014),
750–763.

\bibitem{[Kom2015-IzvRAN]}
{Комаров М.\,А.}
{\it Критерий наилучшего равномерного приближения наипростейшими
дробями в терминах альтернанса}. Изв. РАН. Сер. матем., {\bf 79}:3
(2015), 3--22.

\bibitem{[Kom2015-MZ]}
{Комаров М.\,А.}
{\it Скорость наилучшего приближения констант наипростейшими
дробями и альтернанс}. Матем. заметки, {\bf 97}:5 (2015),
718--732.

\bibitem{[Kom2017-IzvRAN]}
{Комаров М.\,А.} {\it Критерий наилучшего равномерного
приближения наипростейшими дробями в терминах альтернанса. II}.
Изв. РАН. Сер. матем., {\bf 81}:3 (2017), 109--133.

\bibitem{[Kom2017-IzVuz]}
{Комаров М.\,А.} \textit{Аппроксимация посредством
дробно--линейных преобразований наипростейших дробей и их
разностей}. Изв. вузов. Матем. (статья принята к публикации)

\bibitem{[KOND]} {Кондакова Е.\,Н.}
{\it Интерполяция наипростейшими дробями}. Изв. Сарат. ун-та. Нов.
сер. Сер. Математика. Механика. Информатика, {\bf 9}:2 (2009),
30–-37.

\bibitem{[KOND.Abstr2010]}
{Кондакова Е.\,Н.} \textit{Особые случаи интерполяции
посредством наипростейших дробей}. Международная конференция по
дифференциальным уравнениям и динамическим системам (Суздаль,
2010): Тезисы докладов (2010), 105-106.


\bibitem{[KOND-diss]}
{Кондакова Е.\,Н.} \textit{Интерполяция и аппроксимация
наипростейшими дробями: Автореф. дис. ... канд. физ.-матем. наук},
Саратов (2012).

\bibitem{[KOS1]}
{Косухин О.\,Н.} {\it Об аппроксимационных свойствах наипростейших
дробей}. Вестн. МГУ. Сер.~1. Матем., мех., {\bf 4} (2001), 54--58.

\bibitem{[KOS2]}
{Косухин~О.\,Н.} {\it О некоторых нетрадиционных методах
приближения, связанных с~комплексными полиномами: Дисс....канд.
физ.-мат. наук}. МГУ, М., 2005.

\bibitem{[KOS3]}
{Косухин О.\,Н.} {\it Об оценках расстояний от полюсов наипростейших
дробей до компактов}. Материалы Международной научной
конф., посв. 105-летию акад. С.~М.~Никольского (МГУ им. М.~В.~Ломоносова (17-19 мая 2010 г.)), МГУ, М., (2010) 25.

\bibitem{Krylov}
{Крылов В.\,И.} {\it Приближенное вычисление интегралов}. Наука, М., 1967.

\bibitem{Kurosh} {Курош А.\,Г.} {\it Курс высшей алгебры}. Физматлит, М., 1963.

\bibitem{[KUSIS]}
{Кусис П.} {\it Введение в теорию пространств
$H^{p}$}. Мир, М., 1984.

\bibitem{[Minc]} {Минк Х.} {\it Перманенты}. Мир, М. 1982.

\bibitem{[MOROZOV]} {Морозов А.\,Н.} {\it Об одном описании пространств
дифференцируемых функций}. Матем. заметки, {\bf 70}:5 (2001),
758–768.

\bibitem{[NIK]} {Николаев Е.\,Г.} {\it Геометрическое свойство корней
многочленов}. Вестн. МГУ. Сер. 1. Матем., мех., {\bf 5} (1965),
23--26.

\bibitem{Nikolskiy}
{Никольский С.\,М.} {\it Неравенства для целых функций
конечной степени и их применение в теории дифференцируемых функций
многих переменных}. Тр. МИАН, {\bf 38} (1951), 244--278.

\bibitem{[NOV2]}
{Новак Я.\,В.} {\it О наилучшем локальном приближении наипростейшими
дробями}. Матем. заметки, {\bf 84}:6 (2008), 882–887.

\bibitem{[NOV-diss]}
{Новак Я.\,В.} {\it Апроксимацiйнi та iнтерполяцiйнi властивостi найпростiших
дробiв. Дисс. . . . канд. физ.-матем. наук}. ИМ НАН Украины, Киев,
2009.

\bibitem{[NOV*]} {Новак Я.\,В.} {\it Критерiй iснування неперервних похiдних у
функцiй з класу $L_p$ на вiдрiзку в термiнах локальних наближень
найпростiшими дробами}.  Доповiдi НАН України, {\bf 5} (2009),
36--40.

\bibitem{[NOV2010]} {Новак Я.\,В.} {\it Критерiй типу Колмогорова для найпростiших
дробiв}. Збiрник праць Iн-ту математики НАН України, {\bf 7}:2
(2010), 385–392.

\bibitem{[PRO]}
{Протасов В.\,Ю.} {\it Приближения наипростейшими~дробями и
преобразование Гильберта}. Изв. РАН. Сер. матем., {\bf 73}:2
(2009), 123--140.

\bibitem{Timan}
{Тиман А.\,Ф.} {\it Теория приближения функций действительного
переменного}. Физматгиз, М., 1960.

\bibitem{WOLSH} {Уолш Д.\,Л.} {\it Интерполяция и аппроксимация рациональными
функциями в комплексной области}. М., 1961.

\bibitem{[FRN2]}
{Фрянцев A.\,В.} {\it О~численной
аппроксимации дифференциальных полиномов}. Изв. Сарат. ун-та. Нов.
сер. Сер. Математика. Механика. Информатика, {\bf 7}:2 (2007),
39--43.

\bibitem{[FRN3]}
{Фрянцев A.\,В.}
{\it О полиномиальных решениях линейных дифференциальных
уравнений}. Успехи математ. наук, {\bf 63}:3 (2008), 149–150.

\bibitem{Chu2010} {Чунаев П.\,В.} {\it Об одном нетрадиционном методе
аппроксимации}. Тр. МИАН, {\bf 270} (2010), 281--287.

\bibitem{Chu2012} {Чунаев П.\,В.} {\it Об экстраполяции аналитических
функций суммами вида $\sum_k \lambda_k h (\lambda_k z)$}. Матем.
заметки, {\bf 92}:5 (2012), 794--797.

\bibitem{[EID1]}  {Эйдерман В.\,Я.} {\it Оценки картановского типа для потенциалов с ядром Коши и с
действительными ядрами}. Матем. сб., {\bf 198}:8 (2007), 115–-160.
\end{thebibliography}
\end{document}